\documentclass{article}
\usepackage[utf8]{inputenc}
\usepackage{amssymb,amsmath,amsthm,amsfonts,latexsym}
\usepackage{bbm}
\usepackage[toc,page]{appendix}
\usepackage{tikz}
\usetikzlibrary{decorations.pathreplacing,angles,quotes}
\usepackage{stackrel}
\usepackage{fancyvrb}
\usepackage{float}
\usepackage{epsfig}
\usepackage{graphicx, color, subfigure}
\usepackage{siunitx}
\usepgflibrary{arrows}
\newtheorem{theorem}{Theorem}[section]
\newtheorem{lemma}[theorem]{Lemma}

\newtheorem{remark}{Remark}[section]
\usepackage{caption}
\usepackage{afterpage}
\usepackage{fullpage}
\usepackage{mathtools}
\usepackage[mathscr]{eucal}
\usepackage{pdfpages}
\usepackage[T1]{fontenc}
\usepackage[english]{babel}
\usepackage{xspace}
\usepackage{amscd,bm}
\usepackage{xcolor,soul}
\usepackage{verbatim}
\usepackage{csvsimple}
\usepackage{hyperref}
\numberwithin{equation}{section}

\newtheorem{ex}{Example}[section]

\DeclareMathAlphabet{\mathpzc}{OT1}{pzc}{m}{it}

\title{{The Hybrid-dimensional Darcy's Law: A Reinterpreted Discrete Fracture Model for Fracture and Barrier Networks on Non-conforming Meshes} \footnote{\baselineskip 1.2pc The second author was funded by Funded by the National Nature Science Foundation of China 52074336, the Major Science and Technology Projects of China National Petroleum Corporation ZD2019-183-008, and the Fundamental Research Funds for the Central Universities (Grant No.18CX05029A). The last author was funded by the NSF grant DMS-1818467}}
\author{Ziyao Xu\footnote{Division of Applied Mathematics, Brown University, Providence, RI 02906. E-mail: ziyao\_xu1@brown.edu}, \quad Zhaoqin Huang\footnote{Department of Petroleum engineering, China University of Petroleum, Qingdao, Shangdong, China. E-mail: huangzhqin@upc.edu.cn}, \quad Yang Yang\footnote{Department of Mathematical Sciences, Michigan Technological University, Houghton, MI 49931. E-mail: yyang7@mtu.edu}}
\date{}
\begin{document}
\baselineskip 1.2pc
\maketitle
\begin{abstract}
\baselineskip 1.2pc
In this paper, we extend the reinterpreted discrete fracture model for flow simulation of fractured porous media containing flow blocking barriers on non-conforming meshes. The methodology of the approach is to modify the traditional Darcy's law into the hybrid-dimensional Darcy's law where fractures and barriers are represented as Dirac-$\delta$ functions contained in the permeability tensor and resistance tensor, respectively.
As a natural extension of the reinterpreted discrete fracture model \cite{paper} for highly conductive fractures, this model is able to account for the influence of both highly conductive fractures and blocking barriers accurately on non-conforming meshes.
The local discontinuous Galerkin (LDG) method is employed to accommodate the form of the hybrid-dimensional Darcy's law and the nature of the pressure/flux discontinuity.
The performance of the model is demonstrated by several numerical tests.

\end{abstract}

\section{Introduction}
Fractures are ubiquitous in crustal rocks as a result of geological processes such as jointing and faulting, which have a substantial impact on the hydraulic properties of the rock mass. Fractures can serve as high-conductivity conduits or barriers for fluid and solute transport depending on whether fractures are filled with impermeable minerals or not \cite{Caine1996, Berkowitz2002, Obeysekara2018}.
Understanding and modeling fluid flow in fractured media is of great importance in various engineering applications, such as oil and gas exploitation, $\text{CO}_2$ sequestration, geothermal extraction, and radioactive waste disposal. As an efficient approach, discrete fracture models (DFM) are usually used in the flow modeling of fractured media, in which the fractures are represented individually. To accurately capture the complexity of a fractured media, unstructured grids are commonly required to explicitly discretize the complex fracture geometries in DFM. Based on DFM, a series of corresponding numerical methods, including FVM \cite{Granet2001, Karimi-Fard2003,Sandve2012,Ahmed2015, Glaser2017}, FEM \cite{Kim2000, Firoozabadi2003,Zhang2013}, mixed FEM \cite{Martin2005,Hoteit2008,Limiter2,Hoteit2006,Moortgat2013,Zidane2014,Zidane2018}, and mimetic finite difference method \cite{Huang2014} have been applied to the flow simulations in fractured media. However, great difficulties can be brought into the mesh generation processes when the complex fracture geometries are highly developed. The use of conforming grids in real geological models is still limited due to the complex gridding and high computational cost.

As a compromise, the embedded discrete fracture model (EDFM), a kind of method based on nonconforming meshes, is proposed \cite{Lee2001,Li2008}. In this approach, the structured grids are usually used to discretize the reservoir domain, and the fractures are directly incorporated into the background grids. Fracture grids are generated by segmenting fractures with the structured background grids, and additional non-neighboring connections and transmissivities are introduced to explicitly consider the influence of fractures. The EDFM has been widely applied to reservoirs containing highly-conductive fractures, and implemented in various complex problems. Moinfar et al.\cite{Moinfar2014} implemented the EDFM in an in-house compositional reservoir simulator for nonvertical fractures. Panfili et al.\cite{Panfili2014} and Fumagalli et al.\cite{Fumagalli2016} applied the EDFM to perform case studies with corner-point grids. Xu and Sepehrnoori \cite{Xu2019} gave the detailed algorithms and workflow for the implementation of EDFM in corner-point grids. Liu et al.\cite{Liu2020} proposed a modified EDFM to improve the fracture discretization process by using two sets of independent grids for matrix and fracture systems, which promotes the modeling of 3D complex fracture geometries in real-field geological structures. However, the conventional EDFM is not suitable in cases when fractures behave as barriers, i.e. the fracture permeability lies far below that of the matrix. To resolve this limitation, the projection-based EDFM (pEDFM) is proposed \cite{Tene2017,Jiang2017}. The effective flow area between adjacent matrix grids is calculated using the original interface area minus the projected area of the fracture segment, which will become zero if the fracture fully penetrates through the matrix cell. Olorode et al.\cite{Olorode2020} extended the pEDFM into three-dimensional compositional simulation of fractured reservoirs. However the pEDFM still cannot describe the complex multiphase flow behavior in the matrix blocks within barrier fractures.

In this paper, we extend the idea in DFM, which accounts for the influence of fractures by superposing the stiffness matrices (FEM-DFM) or fluxes (FVM-DFM) contributed by fractures on that of porous matrix, to model barriers in a similar fashion.
Recall that in the classical FEM-DFM, the fractures are aligned with edges of meshes in order to superpose the stiffness matrices of fractures on that porous matrix conveniently. We have revealed \cite{paper} that such an approach is equivalent to employing the FEM on the same mesh to the modified Darcy's law, in which the permeability tensor of fractures represented by Dirac-$\delta$ functions are superposed on the permeability tensor of porous matrix. The new approach is named as reinterpreted DFM (RDFM) as we did not introduce a new model, yet to write out another representation of the DFM by using partial differential equations (PDEs). It is not difficult to see that the equivalence also holds for FVM-DFM if we apply FVM on the same mesh to the modified Darcy's law. The approach of adding barriers in the model is drawn from the simple observation that the Darcy's law in one space dimension can be written as ${k}^{-1} {u} = -\frac{ d p}{d x}$, in which ${k}^{-1}$ can be explained as the resistance as in the Ohm's law. The barriers thus can be viewed as the region with high resistance, in contrast with the high permeability region of fractures.
Therefore, we can also represent the barriers by Dirac-$\delta$ functions in a similar way as we did in modeling fractures and the only issue remains is to quantify the blocking ability of barriers. Though the idea of the model inherits from DFM, its implementation is more similar to the pEDFM, in the sense that it can be naturally applied on structured meshes, in particular on rectangular meshes, since the model is built on the PDEs level thereby independent of meshes. However, unlike the pEDFM, which handles barriers in a totally different way to fractures, this approach models fractures and barriers in an equal and symmetric way.

The rest of the paper is organized as follows. In section 2, we model barriers in $1$D porous media. In section 3, we extend the model to two space dimensions and fracture barrier networks. In section 4, we establish the numerical scheme of the novel discrete fracture model based on local discontinuous Galerkin methods. The performance of the approach for steady-state single-phase flow in $2$D porous media is demonstrated by several well-known benchmarks in section 5. Finally, we show an application of the model in contaminant transportation in porous media in section 6. We will ends in section 7 with concluding remarks.

\section{Modeling barriers in 1D porous media}

To rationalize the formulation of the hybrid-dimensional Darcy's law, this section explores the hybrid-dimensional model for barriers in $1$D porous media.

\subsection{Equi-dimensional model for barriers in 1D porous media}
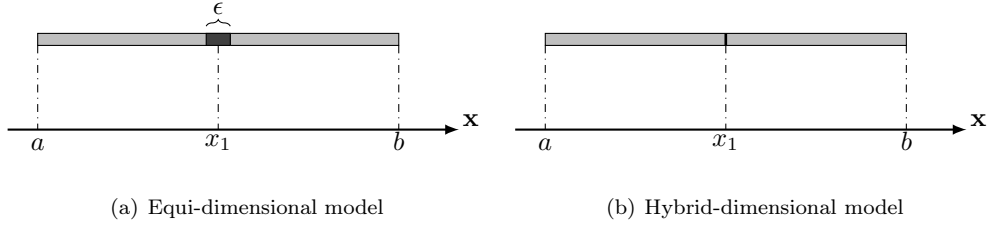
\begin{figure}[!htbp]
\centering
\subfigure[Equi-dimensional model]
{
\begin{tikzpicture}[scale=0.8]
\draw [fill=lightgray] (-3,-0.1) -- (3,-0.1) -- (3,0.1) --(-3,0.1) -- (-3,-0.1) -- cycle;
\draw [fill=darkgray] (-0.2,-0.1) rectangle (0.2,0.1);

\draw[decoration={brace,raise=5pt},decorate]
  (-0.2,0) -- node[above=6pt] {$\epsilon$} (0.2,0);

\draw[thick, -latex] (-3.5,-1.5) -- (4.0,-1.5);
\coordinate [label=$\mathbf{x}$] (h) at(4.2,-1.6) ;

\draw [thin,dash dot] (-3,-0.1) -- (-3,-1.5);
\draw [thin,dash dot] (3,-0.1) -- (3,-1.5);
\draw [thin,dash dot] (0,-0.1) -- (0,-1.5);

\coordinate [label=$a$] (h) at(-3,-2.0) ;
\coordinate [label=$b$] (h) at(3,-2.0) ;
\coordinate [label=$x_1$] (h) at(0,-2.0) ;
\end{tikzpicture}
}
\subfigure[Hybrid-dimensional model]
{
\begin{tikzpicture}[scale=0.8]
\draw [fill=lightgray] (-3,-0.1) -- (3,-0.1) -- (3,0.1) --(-3,0.1) -- (-3,-0.1) -- cycle;
\draw [line width=1] (-0.000,-0.1) -- (-0.000,0.1);

\draw[thick, -latex] (-3.5,-1.5) -- (4.0,-1.5);
\coordinate [label=$\mathbf{x}$] (h) at(4.2,-1.6) ;

\draw [thin,dash dot] (-3,-0.1) -- (-3,-1.5);
\draw [thin,dash dot] (3,-0.1) -- (3,-1.5);
\draw [thin,dash dot] (0,-0.1) -- (0,-1.5);

\coordinate [label=$a$] (h) at(-3,-2.0) ;
\coordinate [label=$b$] (h) at(3,-2.0) ;
\coordinate [label=$x_1$] (h) at(0,-2.0) ;
\end{tikzpicture}
}
\caption{Barrier model in 1D porous media}
\label{fig:Barriers1D}
\end{figure}

For the steady-state single-phase flow in 1D porous media, the Darcy's velocity $u$ and pressure $p$ are governed by the following differential equations,
\begin{eqnarray}
&& u=-k(x)\frac{d p}{d x},\quad x\in(a,b),\label{eq:Darcy1D}\\
&& \frac{d u}{d x} = f, \quad x\in(a,b),\label{Cont1DEqui}
\end{eqnarray}
where $k(x)$ is the permeability of the porous media and $f$ is the source term. The equation (\ref{eq:Darcy1D}) is known as the Darcy's law. When the porous media contains low permeable barrier regions, the permeability $k(x)$ can be expressed as follows,
\begin{equation}\label{Permeability1DbarrierEqui}
k(x) =
 \begin{cases}
      k_m & x\in \Omega_m \\
      k_{\epsilon} & x\in \Omega_{\epsilon}  ,
   \end{cases}
\end{equation}
where $k_m$ is the permeability of porous matrix, $k_{\epsilon}$ is the permeability of barrier, $\Omega_m$ and $\Omega_{\epsilon}$ denote the matrix region and barrier region, respectively. See Figure \ref{fig:Barriers1D}(a) for an illustration, in which the barrier region $\Omega_\epsilon=(x_1-\frac{\epsilon}{2},x_1+\frac{\epsilon}{2})$ is colored dark gray and the matrix region $\Omega_m=(a,b)\setminus\Omega_\epsilon$ is colored light gray. Note that the thickness of barrier region is exaggerated for the sake of visibility here.

We consider the Dirichlet boundary conditions for this problem, i.e.
\begin{equation}\label{Bound1DEqui}
    p(a)=p_a, \quad p(b)=p_b.
\end{equation}

The equations \eqref{eq:Darcy1D}, \eqref{Cont1DEqui} and \eqref{Bound1DEqui} are called the equi-dimensional model for barriers in 1D porous media.

\subsection{Hybrid-dimensional model for barriers in 1D porous media}

Following the idea in our previous work \cite{paper}, which uses Dirac-$\delta$ functions to describe highly conductive fractures, this section explores the hybrid-dimensional representation for barriers in 1D porous media.

If there is no source term, the Darcy's velocity $u$ can be solved from the equations \eqref{eq:Darcy1D}, \eqref{Cont1DEqui} and \eqref{Bound1DEqui}:
\begin{equation} \label{eq:Ohm0}
    u=\frac{p_a-p_b}{\displaystyle \int_{a}^{b}\frac{dx}{k(x)}}.
\end{equation}
In analogy to the Ohm's law $RI=U$, we can rewrite equation (\ref{eq:Ohm0}) as
\begin{equation*}
    {\displaystyle \left(\int_{a}^{b}\frac{dx}{k(x)}\right)u}={p_a-p_b},
\end{equation*}
where $\displaystyle \int_{a}^{b}\frac{dx}{k(x)}$ can be interpreted as the resistance of the porous media.

Replace $b$ by the variable $x\in (a,b)$, the resistance of porous media on interval $(a,x)$ can be approximated as $\displaystyle \int_{a}^{x}\frac{dx}{k(x)}\approx\displaystyle \int_{a}^{x}\frac{dx}{k_m}+\frac{\epsilon}{k_\epsilon}H(x-x_1)$, where $H(x)$ is the Heaviside function defined as $H(x)=1$ when $x \geq 0$ while $H(x)=0$ when $x<0$.
Therefore, the following governing equation of integration form holds:
\begin{equation*}
    {\displaystyle \left(\displaystyle \int_{a}^{x}\frac{dx}{k_m}+\frac{\epsilon}{k_\epsilon}H(x-x_1)\right)u}={p_a-p(x)}, \quad x\in (a,b).
\end{equation*}
The corresponding differentiation form is the hybrid-dimensional barrier model, which can be viewed as a modified version of Darcy's law:
\begin{equation}\label{eq:Darcy1DHybrid1}
    \left(k_m^{-1}+\frac{\epsilon}{k_{\epsilon}}\delta(x-x_1)\right)u=-\frac{d p}{d x},\quad x\in(a,b),
\end{equation}
where $\delta$ is the Dirac-$\delta$ function.
If there are multiple barriers in the media, the equation becomes
\begin{equation}\label{eq:Darcy1DHybrid2}
\left(k_m^{-1}+\sum^{M}_{i=1}\frac{\epsilon_i}{k_{\epsilon_i}}\delta(x-x_i)\right)u=-\frac{d p}{d x},\quad x\in(a,b)
\end{equation}
where $M$ is the number of barriers and $x_i's$ are the positions of barrier's centers.

Equations  \eqref{eq:Darcy1DHybrid2}, \eqref{Cont1DEqui} and \eqref{Bound1DEqui} are called the hybrid-dimensional model for barriers in 1D porous media.

\section{Hybrid-dimensional Darcy's law for 2D porous media}
This section extends the hybrid-dimensional model for barriers from $1$D to $2$D and then presents the complete formulation of the hybrid-dimensional Darcy's law by combining it with the fracture model \cite{paper}.

\subsection{Hybrid-dimensional model for barriers in 2D porous media}
For the steady-state single-phase flow in 2D porous media, the Darcy's velocity $\bm{u}$ and pressure $p$ are governed by the following partial differential equations,
\begin{eqnarray}
{\bf K}^{-1}\bm{u}&=&-\nabla p, \quad x\in\Omega,\\
\nabla\cdot\bm{u}&=&f,\quad x\in\Omega, \label{eq:ContEq}
\end{eqnarray}
where ${\bf K}$ (${\bf K}^{-1}$) is the permeability (resistance) tensor of the porous media.

We consider the mixed boundary condition
\begin{equation}\label{BVP}
    p = p_D,\quad \text{on}~ \Gamma_D\in\partial\Omega,\quad \text{and} \quad \bm{u} \cdot {\bm n}=q_N, \quad \text{on}~ \Gamma_N = \partial\Omega\setminus\Gamma_D,
\end{equation}
where ${\bm n}$ is the unit outer normal vector of the boundary $\partial \Omega$.
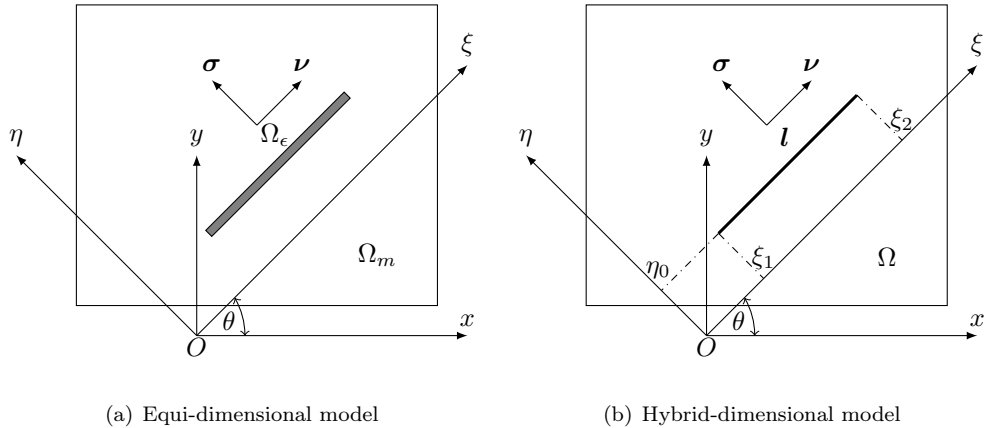
\begin{figure}[!htbp]
\centering
\subfigure[Equi-dimensional model]
{
\begin{tikzpicture}[scale=0.8]
\draw (-3,-2.5) -- (3,-2.5) -- (3,2.5) --(-3,2.5) -- (-3,-2.5) -- cycle;
\draw [fill=gray] (-0.8+0.05,-1.3-0.05) -- (1.5+0.05,1.0-0.05) -- (1.5-0.05,1.0+0.05) -- (-0.8-0.05,-1.3+0.05) --(-0.8+0.05,-1.3-0.05) -- cycle;

\coordinate [label=${\Omega}_m$] (h) at(2,-2) ;

\draw[-latex] (-1,-3) -- (3.5,1.5);
\coordinate [label=$\xi$] (h) at(3.5,1.5) ;
\draw[-latex] (-1,-3) -- (-4,0);
\coordinate [label=$\eta$] (h) at(-4,0) ;
\draw[-latex] (-1,-3) -- (3.5,-3);
\coordinate [label=$x$] (h) at(3.5,-3) ;
\draw[-latex] (-1,-3) -- (-1,0);
\coordinate [label=$y$] (h) at(-1,0) ;
\coordinate [label=$\Omega_\epsilon$] (h) at(0.3,0.0) ;

\draw[-latex] (0,0.5) -- (0.75,1.25);
\coordinate [label=$\bm{\nu}$] (h) at(0.75,1.25) ;
\draw[-latex] (0,0.5) -- (-0.75,1.25);
\coordinate [label=$\bm{\sigma}$] (h) at(-0.75,1.25) ;

\coordinate [label=$O$] (h) at(-1,-3.5) ;
\draw[<->] (-0.2,-3) arc (0:30:1.25);
\coordinate [label=left:$\theta$] (a) at (-0.2,-2.75);

\end{tikzpicture}
}
\subfigure[Hybrid-dimensional model]
{
\begin{tikzpicture}[scale=0.8]
\draw (-3,-2.5) -- (3,-2.5) -- (3,2.5) --(-3,2.5) -- (-3,-2.5) -- cycle;

\coordinate [label=${\Omega}$] (h) at(2,-2) ;

\draw[black, line width=0.4mm] (-0.8,-1.3) -- (1.5,1.0);
\coordinate [label=$\bm{l}$] (h) at(0.3,0.0) ;

\draw[-latex] (0,0.5) -- (0.75,1.25);
\coordinate [label=$\bm{\nu}$] (h) at(0.75,1.25) ;
\draw[-latex] (0,0.5) -- (-0.75,1.25);
\coordinate [label=$\bm{\sigma}$] (h) at(-0.75,1.25) ;

\draw[-latex] (-1,-3) -- (3.5,1.5);
\coordinate [label=$\xi$] (h) at(3.5,1.5) ;
\draw[-latex] (-1,-3) -- (-4,0);
\coordinate [label=$\eta$] (h) at(-4,0) ;

\draw[-latex] (-1,-3) -- (3.5,-3);
\coordinate [label=$x$] (h) at(3.5,-3) ;
\draw[-latex] (-1,-3) -- (-1,0);
\coordinate [label=$y$] (h) at(-1,0) ;

\coordinate [label=$O$] (h) at(-1,-3.5) ;
\draw[<->] (-0.2,-3) arc (0:30:1.25);
\coordinate [label=left:$\theta$] (a) at (-0.2,-2.75);

\draw [dash dot] (-0.8,-1.3) -- (-0.05,-2.05);
\coordinate [label=$\xi_1$] (h) at(-0.25+0.2,-2.25+0.2) ;
\draw [dash dot] (1.5,1.0) -- (2.25,0.25);
\coordinate [label=$\xi_2$] (h) at(2.25,0.25) ;
\draw [dash dot] (-0.8,-1.3) -- (-1.75,-2.25);
\coordinate [label=$\eta_0$] (h) at(-1.75-0.05,-2.25+0.05) ;
\end{tikzpicture}
}
\caption{Barrier model in 2D porous media}
\label{fig:Barriers2D}
\end{figure}

When the porous media contains low permeable barrier strips, see Figure \ref{fig:Barriers2D}(a), the permeability ${\bf K}$ can be expressed as follows,
\begin{equation}\label{EquiDimensionK}
{\bf{K}} =
 \begin{cases}
      {\bf K}_m & x\in \Omega_m \\
      {\bf K}_{\epsilon} & x\in \Omega_{\epsilon},
   \end{cases}
\end{equation}
where $\Omega_m$ is the matrix region and $\Omega_\epsilon$ is the barrier region, and ${\bf K}_m$ and ${\bf K}_\epsilon$ are the permeability tensors of these two regions, respectively.

Analogous to (\ref{eq:Darcy1DHybrid1}) and note that the effect of a barrier on the flow in its tangential direction $\bm{\nu}$ is negligible, we propose the hybrid-dimensional barrier model for $2$D porous media:
\begin{equation}\label{eq:barrier2d}
\left({\bf K}_m^{-1}+\frac{\epsilon}{k_{\epsilon}}\delta(\cdot)\mathbbm{1}(\cdot)\bm{\sigma}\bm{\sigma}^{T}\right)\bm{u}=-\nabla p,
\end{equation}
where $\epsilon$ is the thickness of the barrier, $k_\epsilon$ is the permeability of the barrier in its normal direction, $\bm{\sigma}$ is the unit normal vector of the barrier $\bm{l}$, $\delta(\cdot)$ and $\mathbbm{1}(\cdot)$ are the shorthand of two functions containing the information of the position of the barrier $\bm{l}$.

The full expressions of $\delta(\cdot)$ and $\mathbbm{1}(\cdot)$ under the local coordinate system $\xi O \eta$ and global coordinate system $xOy$ are given as
\begin{equation}\label{eq:shorthand1}
\delta(\cdot) =\delta(\eta-\eta_0)=\delta(-\sin(\theta)x+\cos(\theta)y-\eta_0),
\end{equation}
and
\begin{equation}\label{eq:shorthand2}
\mathbbm{1}(\cdot)=\mathbbm{1}(\xi_1\leq \xi \leq \xi_2)=\mathbbm{1}(\xi_1\leq \cos(\theta)x+\sin(\theta)y\leq \xi_2),
\end{equation}
where $\delta$ is the Dirac-$\delta$ function and $\mathbbm{1}$ is the indicator function defined as $\mathbbm{1}(\text{expr})$ equals $1$ if expr is true while equals $0$ otherwise. For the coordinate systems $\xi O\eta$ and $xOy$, see Figure \ref{fig:Barriers2D} as an illustration.

\subsection{Hybrid-dimensional Darcy’s law for fracture and barrier networks in 2D porous media}
In the previous work \cite{paper}, we proposed the hybrid-dimensional fracture model in 2D porous media:
\begin{equation}\label{eq:fracture2d}
\bm{u}=-\left({\bf K}_m+\epsilon k_{f}\delta(\cdot)\mathbbm{1}(\cdot)\bm{\nu}\bm{\nu}^{T}\right)\nabla p.
\end{equation}
where $\epsilon$ is the thickness of the fracture, $k_f$ is the permeability of the fracture in its tangential direction, $\bm{\nu}$ is the unit tangential vector of the fracture, and $\delta(\cdot)$, $\mathbbm{1}(\cdot)$ are the shorthand of their full expressions as in (\ref{eq:shorthand1}) and (\ref{eq:shorthand2}) that contain the geometric information of the fracture.

Combining the barrier model (\ref{eq:barrier2d}) and fracture model (\ref{eq:fracture2d}) together, we obtain the hybrid-dimensional Darcy's law for single-phase flow in fractured media containing fracture and barrier networks:
\begin{equation}\label{eq:HybridDimDarcy}
\left({\bf{I}}+{\bf K}_m\sum^{M}_{i=1}\frac{\epsilon_i}{k_{\epsilon_ i}}\delta_{i}(\cdot)\mathbbm{1}_i(\cdot)\bm{\sigma}_i\bm{\sigma}_{i}^{T}\right)\bm{u}=-\left({\bf K}_m+\sum^{M+N}_{i=M+1}\epsilon_i k_{f_i}\delta_{i}(\cdot)\mathbbm{1}_i(\cdot)\bm{\nu}_i\bm{\nu}_{i}^{T}\right)\nabla p,
\end{equation}
where $i=1,\ldots,M$ are barriers and $i=M+1,\ldots,M+N$ are fractures.

Note that in $3$D model, we can follow \cite{paper} and replace $\bm{\nu}_i\bm{\nu}_{i}^{T}$ by ${\bf{I}}-\bm{\sigma}_i\bm{\sigma}_{i}^{T}$ on the right hand side of \eqref{eq:HybridDimDarcy}.

\section{Numerical algorithm of the reinterpreted discrete fracture model}
In this section, we establish the numerical algorithm of the RDFM \eqref{eq:HybridDimDarcy}, \eqref{eq:ContEq} and \eqref{BVP} by applying the local discontinuous Galerkin (LDG) methods to it, see \cite{Beatrice2008,Shu1998} for more details of the method. For brevity, we only consider the $2$D case.

\subsection{Basic notations}
For the computational domain $\Omega$, consider the partition $\mathpzc{T}$, i.e. $\Omega=\bigcup_{T\in \mathpzc{T}} T$, where $T$'s are elements in the partition.
We denote $\Gamma=\bigcup_{T\in \mathpzc{T}}\{e| e\in\partial T\}$ to be the collection of cell interfaces and $\Gamma_0=\Gamma\setminus\partial\Omega$ for all interior ones.

For a scalar-valued function $v$, define its jump and average on $e=\partial T_{1}\cap \partial T_{2}\in \Gamma_0$ by
$$[v]_e=v_1\bm{n}_1+v_2\bm{n}_2,\qquad \{v\}_e=\frac{1}{2}(v_1+v_2),$$
where $v_i=v|_{T_i}$ and $\bm{n}_i$ is the unit outer normal vector of edge $\partial T_i$ in $T_i$.
For convenience, if $e\in\partial\Omega$, we denote $[v]_e=v\bm{n}$ and $ \{v\}_e=v,$ where $\bm{n}$ is the unit outer normal of $\partial\Omega$.
Similarly,
the jump and average of a vector-valued function $\bm{w}$ are defined as
$$[\bm{w}]_e=\bm{w}_1\cdot\bm{n}_1+\bm{w}_2\cdot\bm{n}_2,\qquad \{\bm{w}\}_e=\frac{1}{2}(\bm{w}_1+\bm{w}_2).$$
for $e=\partial T_{1}\cap \partial T_{2}\in \Gamma_0$, and $[\bm{w}]_e=\bm{w}\cdot\bm{n}$, $\{\bm{w}\}_e=\bm{w}$ for $e\in\partial\Omega$.

Moreover, define the broken inner product for scalar-valued functions as
$$\left(u,v\right)=\sum_{T\in \mathpzc{T}}\int_{T}u v~ dxdy,$$
and for vector-valued functions as
$$\left(\bm{u},\bm{v}\right)=\sum_{T\in \mathpzc{T}}\int_{T}\bm{u}\cdot \bm{v}~ dxdy.$$

\subsection{LDG scheme}
We take the finite element spaces as
\[V_h=\{v\in L^{2}(\Omega):v|_T\in P^k(T),\:\:\forall T\in\mathpzc{T} \} \text{ and } \bm{W}_h = V_h \times V_h , \]
where $P^k(T)$ denotes the space of polynomials of degree at most $k$ on triangular element $T$,
or
\[V_h=\{v\in L^{2}(\Omega):v|_T\in Q^k(T),\:\:\forall T\in\mathpzc{T} \} \text{ and } \bm{W}_h = V_h \times V_h , \]
where $Q^k(T)$ denotes the space of tensor products of $1$D polynomials of degree at most $k$ on rectangular element $T$.

Lemma \ref{lem:delta} will help us reduce the $2$D integration of the line Dirac-$\delta$ terms in \eqref{eq:HybridDimDarcy} to line integrals on fractures and barriers. One can find the detailed proof in \cite{paper}.
\begin{lemma}\label{lem:delta}
Let $\delta(\cdot)$ and $\mathbbm{1}(\cdot)$ be the shorthand of their full expressions in \eqref{eq:shorthand1} and \eqref{eq:shorthand2}, respectively. For all continuous function $g$ on $\Omega$, the following identity holds:
\begin{equation}
\int_{\Omega} \delta(\cdot)\mathbbm{1}(\cdot)g(x,y)~dx dy =\int_{l} g(x,y)~ ds,
\end{equation}
where ${\bm{l}}$ is the support of $\delta(\cdot)\mathbbm{1}(\cdot)$ as shown in Figure \ref{fig:Barriers2D}(b).
\end{lemma}

Now we introduce the auxiliary variable $\bm{s}=-\nabla p$ and establish the LDG scheme for the model problem \eqref{eq:HybridDimDarcy},\eqref{eq:ContEq} and \eqref{BVP} as follows:

Find $\bm{s},\bm{u}\in\bm{W}_{h}$ and $p\in V_{h}$, such that $\forall \bm{\xi},\bm{\eta}\in\bm{W}_h,\zeta\in V_{h}$,
\begin{equation}\label{scheme:part1}
\left(\bm{s},\bm{\xi}\right)=(p,\nabla\cdot\bm{\xi})-\int_{\Gamma_0\cup\Gamma_{N}}\hat{p}[\bm{\xi}]~ds-\int_{\Gamma_D}p_{D}\bm{\xi}~ds,
\end{equation}
\begin{equation}\label{scheme:part2}
\left(\bm{u},\bm{\eta}\right)+\sum_{i=1}^{M}\int_{\bm{l}_i}\left(\frac{\epsilon_i}{k_{\epsilon_i}}{\bf{K}_{m}}\bm{\sigma}_{i}\bm{\sigma}_{i}^{T}\bm{u}\right)\cdot\bm{\eta}~ds=\left({\bf K}_{m}\bm{s},\bm{\eta}\right)+\sum_{i=M+1}^{M+N}\int_{\bm{l}_{i}}\left(\epsilon_i k_{f_i}\bm{\nu}_i \bm{\nu}_{i}^{T}\bm{s}\right)\cdot\bm{\eta}~ds,
\end{equation}
\begin{equation}\label{scheme:part3}
-\left(\bm{u},\nabla \zeta\right)+\int_{\Gamma_0\cup\Gamma_{D}}\hat{\bm{u}}\cdot[\zeta]~ds+\int_{\Gamma_N}q_N \zeta~ds=\left(f,\zeta\right),
\end{equation}
where ${\hat{\bm{u}}}$ and $\hat{p}$ are numerical fluxes defined by
\begin{equation}\label{NumFlux_u}
\hat{\bm{u}}|_{e} =
 \begin{cases}
      \{\bm{u}\}, & \text{ if } T_1\cup T_2 \text{ is crossed by some barriers } \bm{l}_i, i = 1,2,\cdots,M   \\
      \{\bm{u}\} + \alpha[p],& \text{ otherwise}.
   \end{cases}
\end{equation}
\begin{equation}\label{NumFlux_p}
\hat{p}|_{e} =
 \begin{cases}
      \{p\}+\beta[{\bm u}], & \text{ if } T_1\cup T_2 \text{ is crossed by some barriers } \bm{l}_i, i = 1,2,\cdots,M \\
      \{p\}, & \text{ otherwise}.
   \end{cases}
\end{equation}
on $e=\partial T_{1}\cap \partial T_{2}\in \Gamma_0$, and $\hat{\bm{u}}|_{e}=\bm{u}+\alpha\left(p-p_D\right)\bm{n}$ on $e=\partial T\cap\Gamma_D$, and $\hat{p}|_{e}=p$ on $e=\partial T\cap\Gamma_N$, where $\alpha,\beta>0$ are penalty parameters, and $\bm{n}$ is the unit outer normal vector of $\partial\Omega$.

\begin{remark}\label{rmkflux}
It turns out to be crucial to take the correct penalty on the jump of $\bm{u}$ and $p$ for a successful simulation when defining the numerical fluxes.
The reason why we enforce the penalty in this way is, physically, the normal component of the Darcy's velocity $\bm{u}$ is continuous across a barrier but might be discontinuous across a fracture, while the pressure $p$ is discontinuous when crossing a barrier but is continuous in all other cases.
In the special case where some barriers happen to be aligned (or very close to be aligned) with some cell interfaces, we can take a subtle adaptation on the penalties in the numerical fluxes \eqref{NumFlux_u} and \eqref{NumFlux_p} to improve the accuracy of the simulation, by enforcing the penalty on jump of $\bm{u}$ only on the aligned edges and enforcing the penalty on jump of $p$ on other edges of that cell if no other barrier crosses it.
\end{remark}
\begin{remark}
Unlike other models that generate system of equations only for pressure $p$, the LDG methods give system of equations for the unknown pressure $p$, Darcy's velocity $\bm{u}$ and the pressure gradient $\bm{s}$. This difference makes our methods to be computationally more costly than other methods since the degrees of freedom increases a lot.
However, it's not a total disadvantage since in many applications what really useful is the Darcy's velocity and gradient of pressure rather than the pressure itself, and by computing these variables simultaneously, the accuracy of Darcy's velocity and gradient of pressure is higher than those obtained from the post processing of pressure.
\end{remark}

Two issues in implementation are worthy to be mentioned.
First, as one may notice, if the fracture or barrier $\bm{l}$ happens to be aligned with cell interfaces of a mesh, the ambiguity of the line integral on $\bm{l}$ in scheme \eqref{scheme:part2} needs to be resolved.
Rigorously speaking, we should use a well-defined trace of $\bm{u}$ or $\bm{s}$ in \eqref{scheme:part2} on the cell interfaces to replace of the double-valued $\bm{u}$ or $\bm{s}$, e.g. $\textbf{tr}(\bm{u})|_{\bm{l}}=\{\bm{u}\}$, $\textbf{tr}(\bm{s})|_{\bm{l}}=\{\bm{s}\}$, which is equivalent to divide the fracture or barrier $\bm{l}$ equally into its two neighboring cells.
However, in practice, we can take the advantage of the non-conforming nature of our algorithm to add very tiny random perturbations on the position of $\bm{l}$ and only code for the non-conforming case, which is what we shall do in the numerical experiments in the following section.
Second, to treat the special case where fractures and barriers have intersections, say $\bm{l}_1$ and $\bm{l}_2$ intersect on the cell $T$, we only consider either $\bm{l}_1\cap T$ or $\bm{l}_2\cap T$ on $T$, and remove the other one, depending on whether the fracture $\bm{l}_1$ penetrates the barrier $\bm{l}_2$ (fracture-dominated) or the barrier $\bm{l}_2$ blocks the fracture $\bm{l}_1$ (barrier-dominated).

\subsection{Slope limiter}
Due to the discontinuity of pressure across barriers and numerical features of LDG, the numerical approximation of pressure may have strong oscillations on cells containing barriers.
Therefore, a properly designed slope limiter to eliminate spurious oscillations on those cells is essential in the post-processing of numerical approximation of pressure.

One such limiter was introduced in \cite{Limiter0} and modified in \cite{Limiter1, Limiter2}.
The key idea of the limiter is to achieve an appropriate maximum-principle of pressure on barrier cells by imposing the constraint that the value of pressure on a vertex $o$ lies between the minimum and the maximum of the cell averages of neighboring cells of this vertex while make the change as small as possible.

Since we only apply the slope limiter on cells containing barriers rather than on all cells, we adapt the limiter in \cite{Limiter1, Limiter2} in a way that might be less accurate but is easier to implement.
The slope limiting of $p$ on a target cell $T$ with vertices $o_1, o_2, \ldots, o_n$ ($n=4$ for rectangular meshes) is shown as follows:
\begin{equation}\label{SlopeLimiter}
\tilde{p}_T(x) = \bar{p}_T +\theta \left(p_T(x)-\bar{p}_T\right),
\end{equation}
where $\bar{p}_T=\frac{1}{|T|}\int_{T}p_T(x)dx$ is the cell average, and $\theta\in[0,1]$ is chosen to be the largest possible value such that $\tilde{p}_T(o_i)\in[p^{\min}_i, p^{\max}_i]$, for $i=1,\ldots, n$, where $p^{\min}_i=\min_{o_i\in K}\{\bar{p}_K\}$ and $p^{\max}_i=\max_{o_i\in K}\{\bar{p}_K\}$. It's clear from \eqref{SlopeLimiter} that this limiter doesn't change the cell average.

The calculation of $\theta$ in \eqref{SlopeLimiter} can be explicitly given as follows:
\begin{equation}
\theta = \min_{i=1,2,\ldots, n}{\theta_i}, \quad \text{where}~ \theta_i= \begin{cases}
\min\{\frac{p^{\min}_{i}-\bar{p}_T}{p_T(o_i)-\bar{p}_T}, 1\}, & \text{if}~ p_T(o_i)<\bar{p}_T,\\
1 & \text{if}~ p_T(o_i)=\bar{p}_T ,\\
\min\{\frac{p^{\max}_{i}-\bar{p}_T}{p_T(o_i)-\bar{p}_T}, 1\}, & \text{if}~ p_T(o_i)>\bar{p}_T,
\end{cases}
\end{equation}

For the original slope limiter, which is of course applicable in our algorithm, one can refer to the appendix of \cite{Limiter2}, where an intuitive graphical explanation is given.

On cells where the penalty on jump of $p$ is imposed on some edges, the pressure oscillation is greatly depressed by the scheme itself, thereby we do not apply the slope limiter on them.
These include the cells that are not crossed by any barrier, and the cells whose only partial edges are aligned with barriers as is discussed in Remark \ref{rmkflux}.

\section{Numerical experiments}\label{NumExp}
In this section, we provide four numerical experiments, in an order of increasing geometric complexity, to demonstrate the performance of the RDFM proposed in the paper.
The implementations use the piecewise linear finite element spaces, i.e. $k=1$ in $V_h$, unless otherwise stated.
The Examples \ref{ex3} and \ref{ex4} are benchmarks taken from \cite{Benchmark,benchmark2,benchmark3,ThesisXFEM,LMFEM2D}, so one can refer to the solutions in these articles for reference.
For simplicity, the flow in all examples is driven by boundary conditions thus the source term $f=0$ in \eqref{eq:ContEq}.

Flemisch et al. \cite{Benchmark} published the data of grids and solutions of several benchmarks computed by a number of DFM algorithms on \cite{Website}, which enables us to evaluate and compare our model with the existing ones. We declare that all the data and figures used in Example \ref{ex3} and \ref{ex4} for references, evaluations and comparison come from them.


\begin{ex}\label{ex1}
\textbf{Convergence test}

In this example, we test the convergence of our algorithm.
We provide the exact solutions of two particular cases, which are the single fracture case and the single barrier case described in case (a) and case (b), respectively. In both cases, the domain is set to be $\Omega=[-1,1]\times[-1,1]$ and the permeability of the matrix ${\bf{K}}_m={\bf{I}}$.

\textbf{Case (a): Single-fracture}
Consider the fluid flow in the porous media containing a single fracture across the domain.
Suppose the fracture passes the origin with its tangential direction $\bm{\nu}=(\cos(\theta),\sin(\theta))^{T}$, the thickness $\epsilon=10^{-4}$, and the tangential permeability $k_f=2\times 10^{4}$, then the corresponding governing equation is
\begin{equation*}
\bm{u}=-\left({\bf I}+2\delta(-\sin(\theta)x+\cos(\theta)y)\begin{bmatrix}
\cos^2(\theta) & \sin(\theta)\cos(\theta) \\
\sin(\theta)\cos(\theta) & \sin^2(\theta)
\end{bmatrix}\right)\nabla p,\quad \nabla\cdot\bm{u}=0
\end{equation*}

One can verify that $p(x,y)=\sin(\cos(\theta)x+\sin(\theta)y)e^{|-\sin(\theta)x+\cos(\theta)y|}$ is the exact solution of the equation under the corresponding Dirichlet boundary conditions.

\textbf{Case (b): Single-barrier}
Consider the fluid flow in the porous media containing a single barrier across the domain.
Suppose the barrier passes the origin with its tangential direction $\bm{\nu}=(\cos(\theta),\sin(\theta))^{T}$, the thickness $\epsilon=10^{-4}$, and the normal permeability $k_\epsilon= 10^{-4}$, then the corresponding governing equation is
\begin{equation*}
\left({\bf I}+\delta(-\sin(\theta)x+\cos(\theta)y)\begin{bmatrix}
\sin^2(\theta) & -\sin(\theta)\cos(\theta) \\
-\sin(\theta)\cos(\theta) & \cos^2(\theta)
\end{bmatrix}\right)\bm{u}=-\nabla p,\quad \nabla\cdot\bm{u}=0
\end{equation*}

The exact solution is $p(x,y)=(\sin(\theta)-\cos(\theta))x-(\sin(\theta)+\cos(\theta))y+H(\sin(\theta)x-\cos(\theta)y)$, where $H(x)$ is the Heaviside function, under the corresponding Dirichlet boundary conditions.
\end{ex}

\begin{figure}[!htbp]
\subfigure[$\theta=0$, conforming mesh]{\includegraphics[width = 3in]{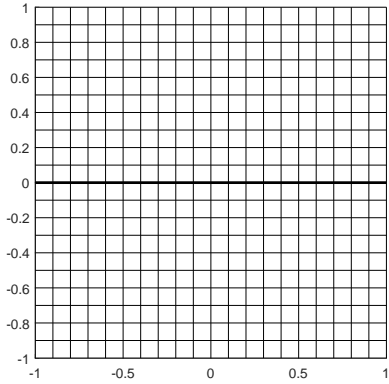}}
\subfigure[$\theta=1$, non-conforming mesh]{\includegraphics[width = 3in]{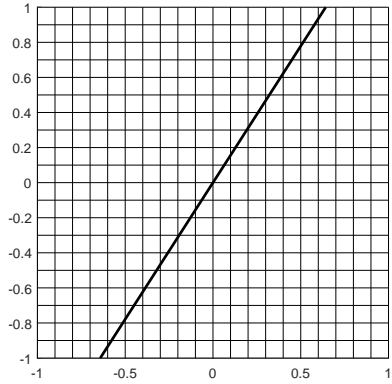}} \\
\caption{The single fracture/barrier with different angles in Example \ref{ex1} }\label{fig:ex1meshes}
\end{figure}

\begin{table}[!htbp]
\begin{center}
\begin{tabular}{c|c|c|c|c|c|c|c|c}
\hline
&\multicolumn{4}{c|}{Case (a) }&\multicolumn{4}{c}{Case (b) }   \\
\hline
&\multicolumn{2}{c|}{$\theta=0$ }&\multicolumn{2}{c|}{$\theta=1$ }&\multicolumn{2}{c|}{$\theta=0$ }&\multicolumn{2}{c}{$\theta=1$}\\\hline
$N\times N$ & $L^1$ error & order& $L^1$ error & order& $L^1$ error & order& $L^1$ error & order\\
\hline
$20\times20$  &1.82e-03& -- &1.26e-02& -- &1.89e-04& -- &6.10e-02& -- \\
$40\times40$ &4.58e-04&1.99&9.02e-03&0.48&4.73e-05&2.00&3.08e-02&0.98\\
$80\times80$ &1.14e-04&2.00&4.35e-03&1.05&1.18e-05&2.00&1.54e-02&1.00\\
$160\times160$ &2.85e-05&2.00&2.64e-03&0.72&2.96e-06&2.00&8.39e-03&0.88\\
$320\times320$ &7.10e-06&2.00&1.41e-03&0.91&7.39e-07&2.00&3.72e-03&1.18\\
\hline
$N\times N$ & $ L^2 $ error & order& $L^2$ error & order& $L^2$ error & order& $L^2$ error & order\\
\hline
$20\times20$ &1.34e-03& -- &9.74e-03& -- &1.23e-04& -- &1.47e-01& -- \\
$40\times40$ &3.39e-04&1.98&7.32e-03&0.41&3.07e-05&2.00&1.08e-01&0.44\\
$80\times80$ &8.47e-05&2.00&3.49e-03&1.07&7.68e-06&2.00&7.60e-02&0.51\\
$160\times160$ &2.11e-05&2.01&2.10e-03&0.73&1.92e-06&2.00&5.38e-02&0.50\\
$320\times320$ &5.26e-06&2.00&1.12e-03&0.91&4.80e-07&2.00&3.83e-02&0.49\\
\hline
\end{tabular}
\caption{\label{table:ex1_k1} Error table of Example \ref{ex1} with order $k=1$}
\end{center}
\end{table}

\begin{table}[!htbp]
\begin{center}
\begin{tabular}{c|c|c|c|c|c|c|c|c}
\hline
&\multicolumn{4}{c|}{Case (a) }&\multicolumn{4}{c}{Case (b) }   \\
\hline
&\multicolumn{2}{c|}{$\theta=0$ }&\multicolumn{2}{c|}{$\theta=1$ }&\multicolumn{2}{c|}{$\theta=0$ }&\multicolumn{2}{c}{$\theta=1$}\\\hline
$N\times N$ & $L^1$ error & order& $L^1$ error & order& $L^1$ error & order& $L^1$ error & order\\
\hline
$20\times20$  &1.04e-04& -- &2.12e-03& -- &1.89e-05& -- &5.73e-02& -- \\
$40\times40$ &1.31e-05&2.98&5.95e-04&1.83&2.37e-06&3.00&3.12e-02&0.88\\
$80\times80$ &1.65e-06&3.00&2.46e-04&1.27&2.96e-07&3.00&1.52e-02&1.04\\
$160\times160$ &2.06e-07&3.00&1.32e-04&0.90&3.69e-08&3.00&8.44e-03&0.85\\
\hline
$N\times N$ & $ L^2 $ error & order& $L^2$ error & order& $L^2$ error & order& $L^2$ error & order\\
\hline
$20\times20$ &7.94e-05& -- &2.30e-03& -- &1.23e-05& -- &1.41e-01& -- \\
$40\times40$ &9.71e-06&3.03&8.88e-04&1.37&1.54e-06&3.00&1.08e-01&0.39\\
$80\times80$ &1.20e-06&3.02&3.18e-04&1.48&1.92e-07&3.00&7.40e-02&0.54\\
$160\times160$ &1.48e-07&3.01&1.74e-04&0.87&2.40e-08&3.00&5.45e-02&0.44\\
\hline
\end{tabular}
\caption{\label{table:ex1_k2} Error table of Example \ref{ex1} with order $k=2$ }
\end{center}
\end{table}

By taking different $\theta$, we can test the errors and rates of convergence in both conforming and non-conforming scenarios, see Figure \ref{fig:ex1meshes}.
The error tables with the order of finite element spaces $k=1$ and $k=2$ are given in Table \ref{table:ex1_k1} and Table \ref{table:ex1_k2}\footnote{Due to the restriction of the memory of the computer, we are unable to compute the solution on $320\times320$ mesh for $k=2$.}, respectively.

A brief summery is given as follows.

For tests of case(a), the rates of convergence are optimal when meshes are aligned with the fracture (as aforementioned, the meshes are not rigorously conforming but almost conforming up to a tiny perturbation on the position of the fracture for convenience of coding, so this condition is not so strong.)
Further tests show that the rates of convergence remain optimal on conforming meshes for higher order $k$'s, under the setting of penalty parameter $\alpha=O(1/h^{k+1})$.
As for the non-conforming meshes, the rate of convergence remains first order for $k\geq1$, due to the low regularity of the exact solution ($p\in H^{1}(\Omega)$). But the errors are still reduced significantly when $k$ increases.
It's notable that we found $\alpha=O(1/h^{3})$ is needed for the algorithm to be convergent for $k\geq 1$ on non-conforming meshes.

For tests of case(b), the rates of convergence are optimal when meshes are aligned with the barrier.
Further tests show that the rates of convergence remain optimal on conforming meshes for higher order $k$'s, under the setting of penalty parameter $\alpha=O(1/h), \beta = O(1/h^{k+1})$.
As for the non-conforming meshes, the rate of convergence remains first order in $L^1$ norm and half order in $L^2$ norm for $k\geq1$, due to the low regularity of the exact solution ($p\in H^{s}(\Omega), s<\frac12$). And the errors are almost unchanged when $k$ increases.
It's notable that we found setting $\alpha=O(1/h), \beta=O(1/h)$ is sufficient for the algorithm to be convergent for $k\geq 1$ on non-conforming meshes.

\begin{figure}[!htbp]
\centering
{\includegraphics[width = 3in]{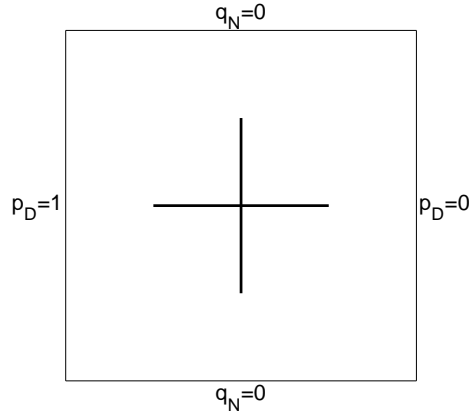}} \\
\caption{Domain and boundary conditions of Example \ref{ex2}}\label{fig:ex2setting}
\end{figure}

\begin{ex}\label{ex2}
\textbf{Cross-shaped networks}

In this example, we test the cross-shaped fracture network and barrier network in case (a) and case (b), respectively.
In both cases, the domain is the unit square $\Omega=[0,1]\times[0,1]$ with permeability of porous matrix ${\bf K}_m={\bf I}$.
The region of fractures/barriers is $[0.25,0.75]\times[0.4995,0.5005]\cup[0.4995,0.5005]\times[0.25,0.75]$, see Figure \ref{fig:ex2setting}.
Moreover, the Dirichlet boundary conditions $p_D=1$ and $p_D=0$ are imposed on the left and right boundaries respectively, and the top and bottom boundaries are set to be impermeable, i.e. $q_N=0$.
This example is the same with the test case given by \cite{Tene2017}, so one can refer to the solution of it for a comparison.

\textbf{Case (a): Cross-shaped fractures}  The permeability of fractures is $k_f=10^{8}$.

\textbf{Case (b): Cross-shaped barriers} The permeability of barriers is $k_\epsilon=10^{-8}$.
\end{ex}

\begin{figure}[!htbp]
\subfigure[Reference solution of case (a)]{\includegraphics[width = 3in]{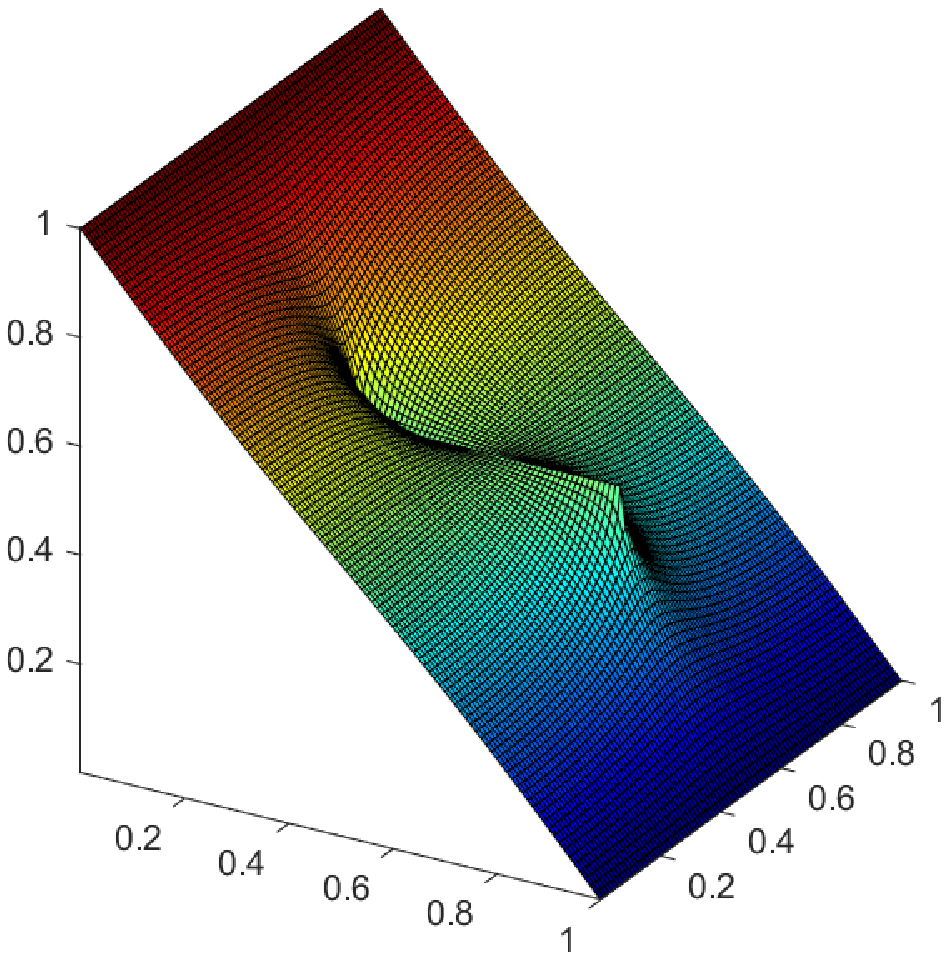}}
\subfigure[Reference solution of case (b)]{\includegraphics[width = 3in]{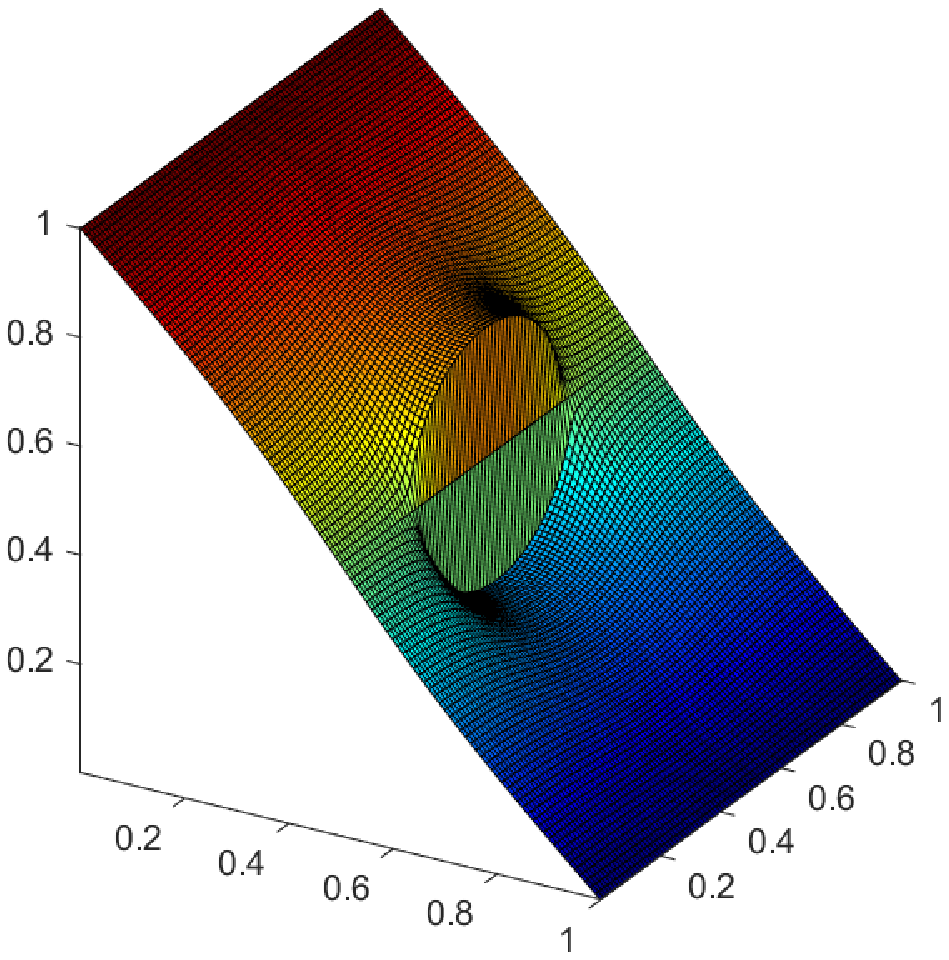}} \\
\caption{The reference solutions of Example \ref{ex2} }\label{fig:ex2ref}
\end{figure}

\begin{figure}[!htbp]
\subfigure[Solution of case (a) on $10\times10$ mesh]{\includegraphics[width = 3in]{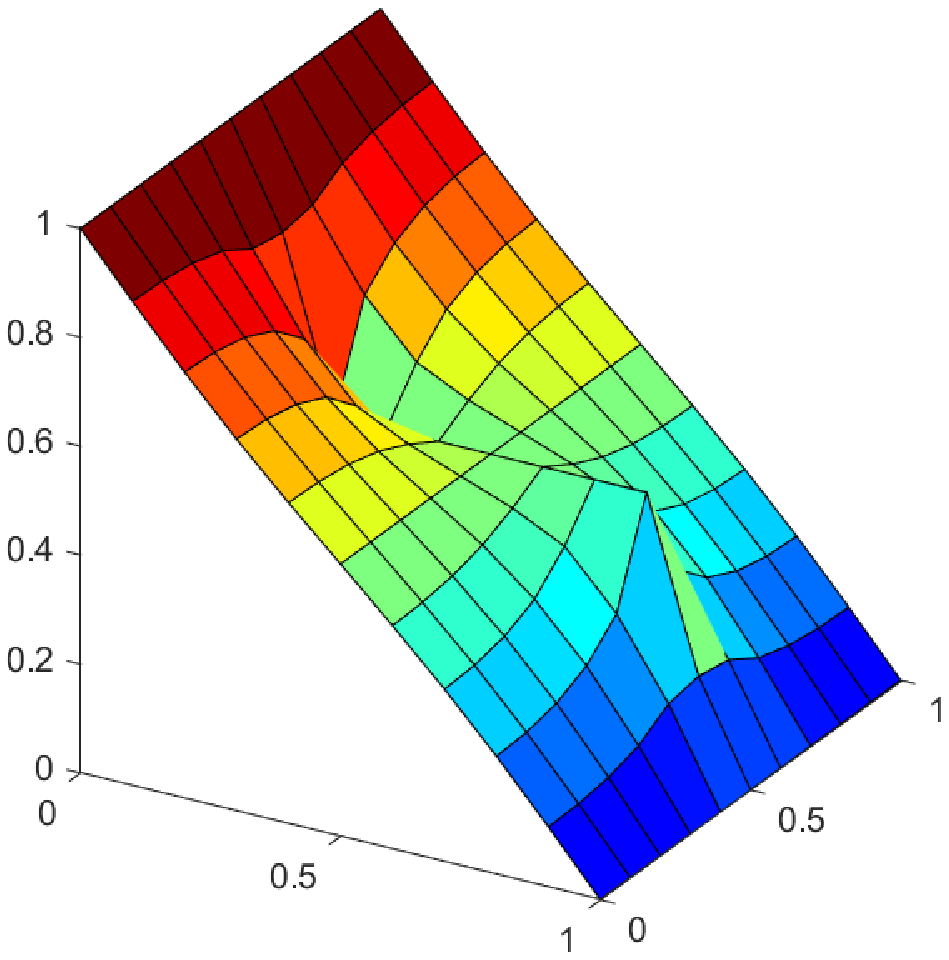}}
\subfigure[Solution of case (a) on $11\times11$ mesh]{\includegraphics[width = 3in]{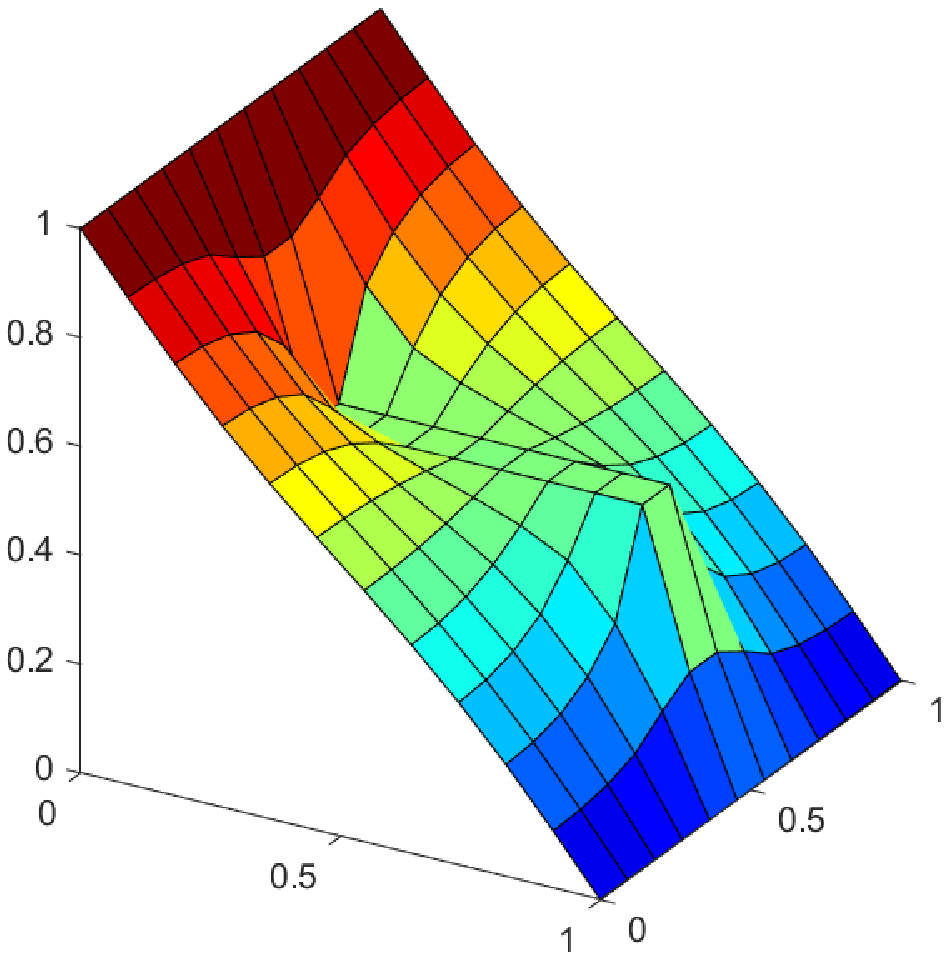}} \\
\subfigure[Solution of case (a) on $20\times20$ mesh]{\includegraphics[width = 3in]{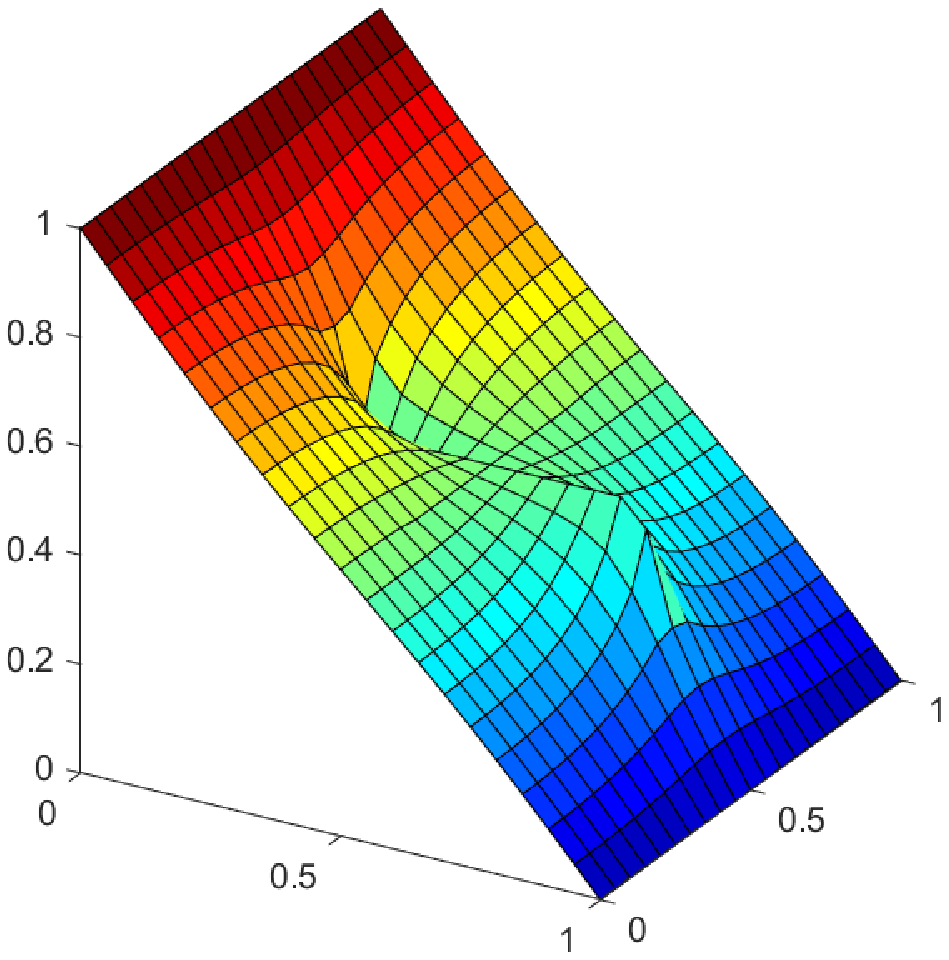}}
\subfigure[Solution of case (a) on $21\times21$ mesh]{\includegraphics[width = 3in]{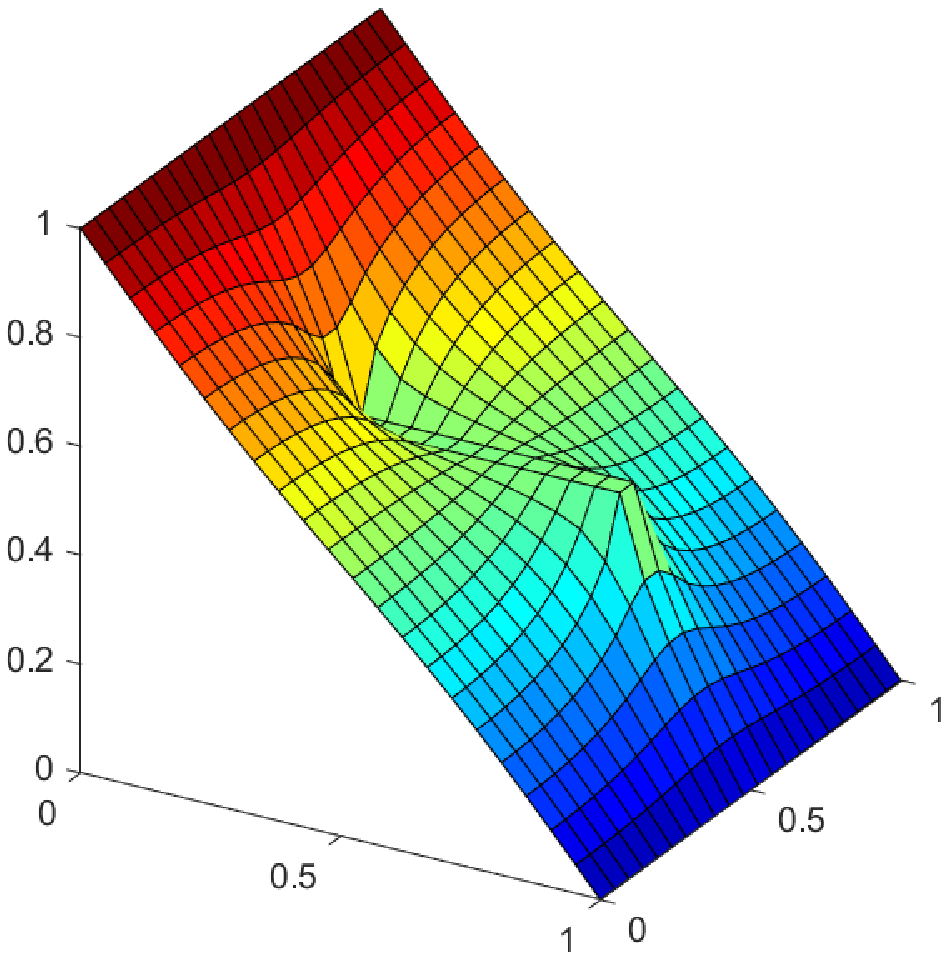}} \\
\caption{Simulation results of case (a) of Example \ref{ex2} }\label{fig:ex2solutionsCaseA}
\end{figure}

\begin{figure}[!htbp]
\subfigure[Solution of case (b) on $10\times10$ mesh]{\includegraphics[width = 3in]{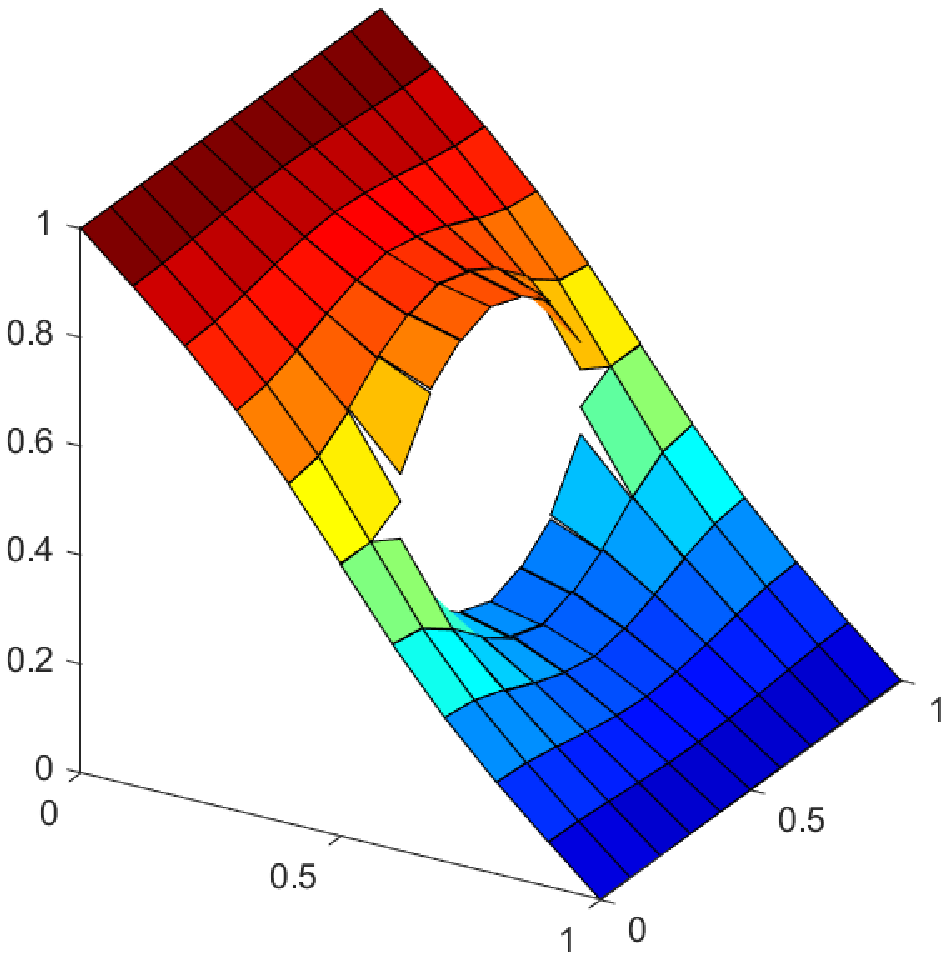}}
\subfigure[Solution of case (b) on $11\times11$ mesh]{\includegraphics[width = 3in]{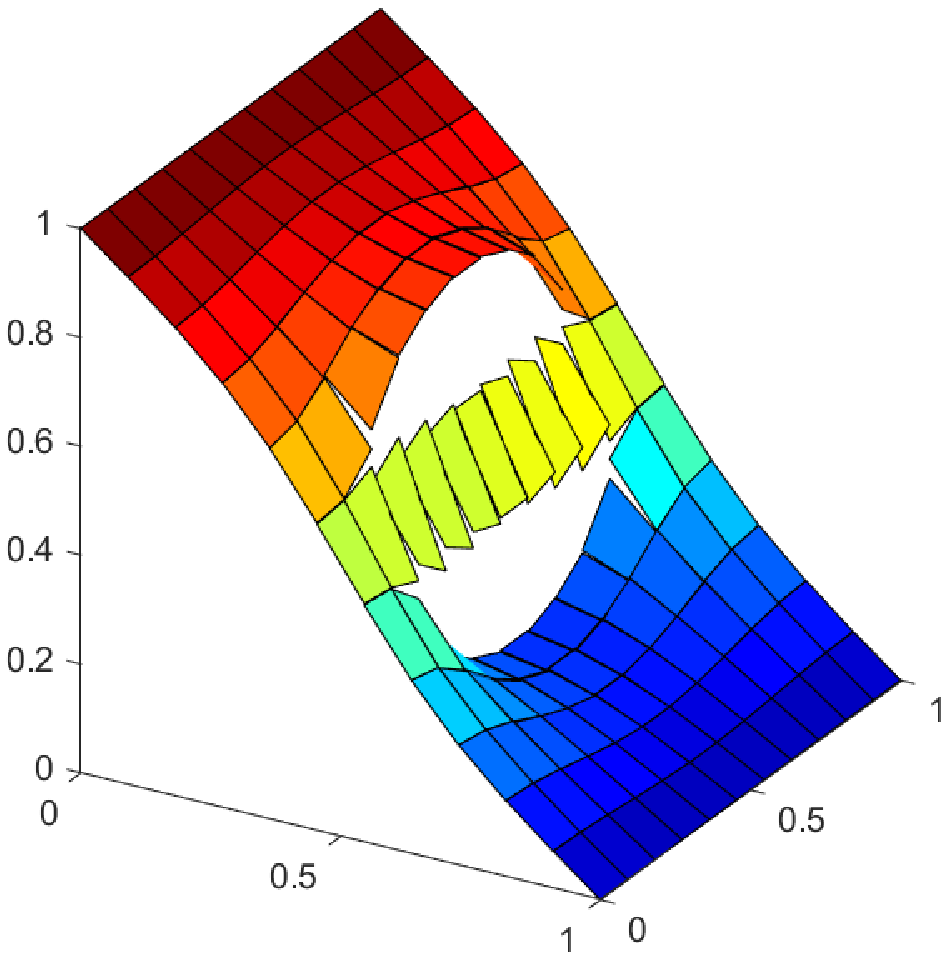}} \\
\subfigure[Solution of case (b) on $20\times20$ mesh]{\includegraphics[width = 3in]{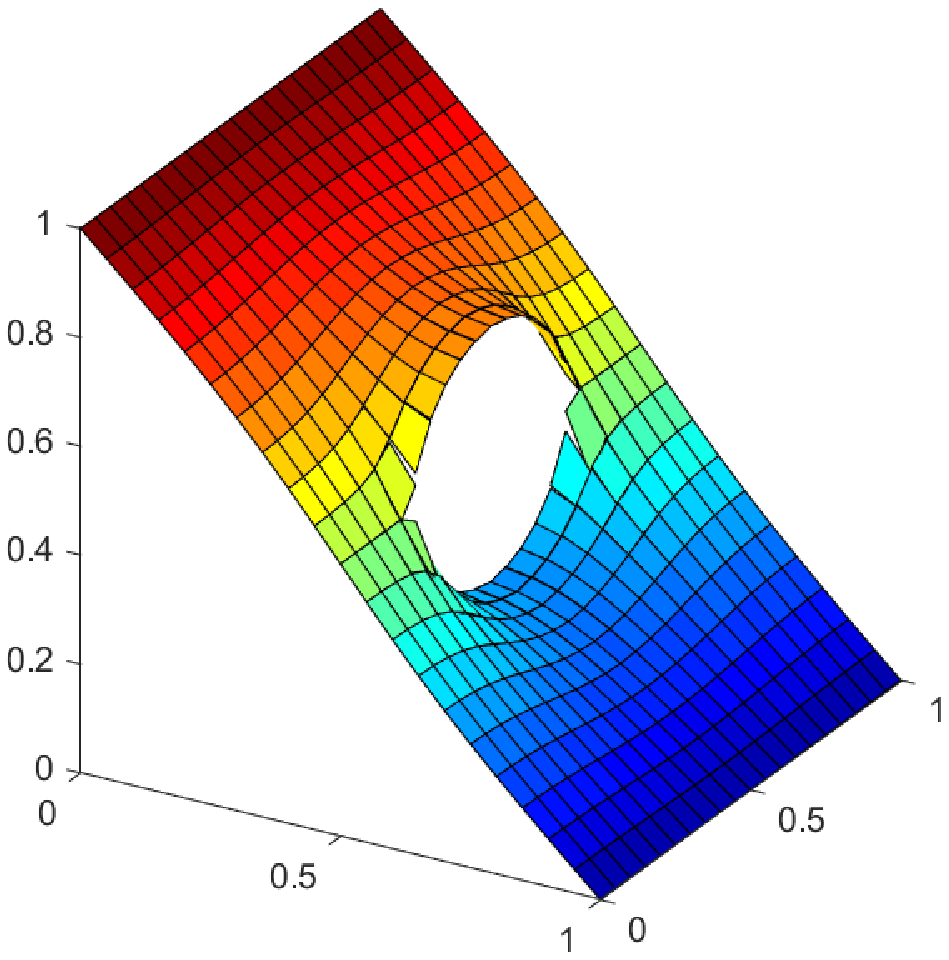}}
\subfigure[Solution of case (b) on $21\times21$ mesh]{\includegraphics[width = 3in]{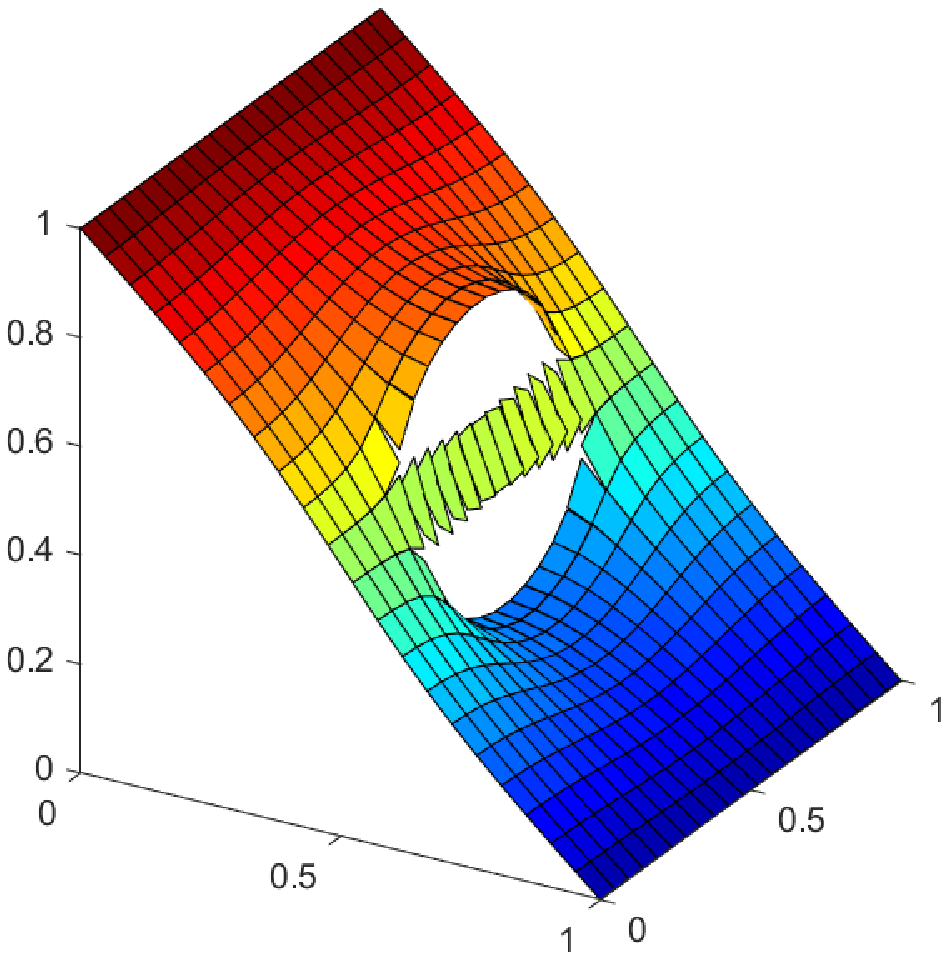}} \\
\caption{Simulation results of case (b) of Example \ref{ex2} }\label{fig:ex2solutionsCaseB}
\end{figure}

The reference solutions are computed by finite element method on the $1001\times1001$ fully resolved mesh and plotted in Figure \ref{fig:ex2ref}.
The solutions of RDFM on different meshes for case (a) and case (b) are shown in Figure \ref{fig:ex2solutionsCaseA} and Figure \ref{fig:ex2solutionsCaseB}, respectively.

Comparing the numerical solutions of RDFM with the reference solutions, we can see the effect of conductive fractures and blocking barriers are captured well on coarse meshes.
Moreover, the pressure along fractures on $11\times11$ and $21\times21$ meshes of case (a) are flat. This phenomenon is reasonable since the domain, meshes and boundary conditions are symmetric.

\begin{figure}[!htbp]
\centering
{\includegraphics[width = 3in]{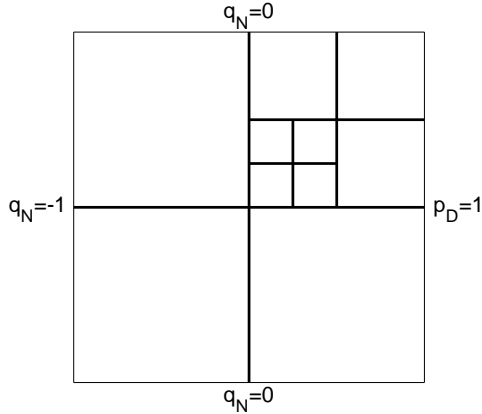}} \\
\caption{Domain and boundary conditions of Example \ref{ex3}}\label{fig:ex3setting}
\end{figure}

\begin{ex} \label{ex3}
\textbf{Regular networks}

In this example, we test regular networks of fractures and barriers.
This example was originally from \cite{benchmark3} and modified by \cite{Benchmark}, which contains two subcases.
The case (a) is a fracture network and the case (b) is a barrier network.
In both cases, the domain is the unit square $\Omega=[0,1]\times[0,1]$ with permeability ${\bf K}_m={\bf I}$.
The networks are composed of six fractures/barriers with central axis $x=0.5, y=0.5, x=0.75, y=0.75, x=0.625, y=0.625$ and uniform thickness $\epsilon = 10^{-4}$, see Figure \ref{fig:ex3setting}.
The left boundary is imposed by Neumann boundary condition $q_N=-1$ and the right boundary is imposed by Dirichlet boundary condition $p_D=1$.
The top and bottom boundaries are impermeable.

\textbf{Case (a): Regular fracture network} The permeability of fractures is $k_f=10^{4}$.

\textbf{Case (b): Regular barrier network} The permeability of barriers is $k_\epsilon=10^{-4}$.
\end{ex}

\begin{figure}[!htbp]
\subfigure[Solution of case (a) on $25\times25$ mesh]{\includegraphics[width = 3in]{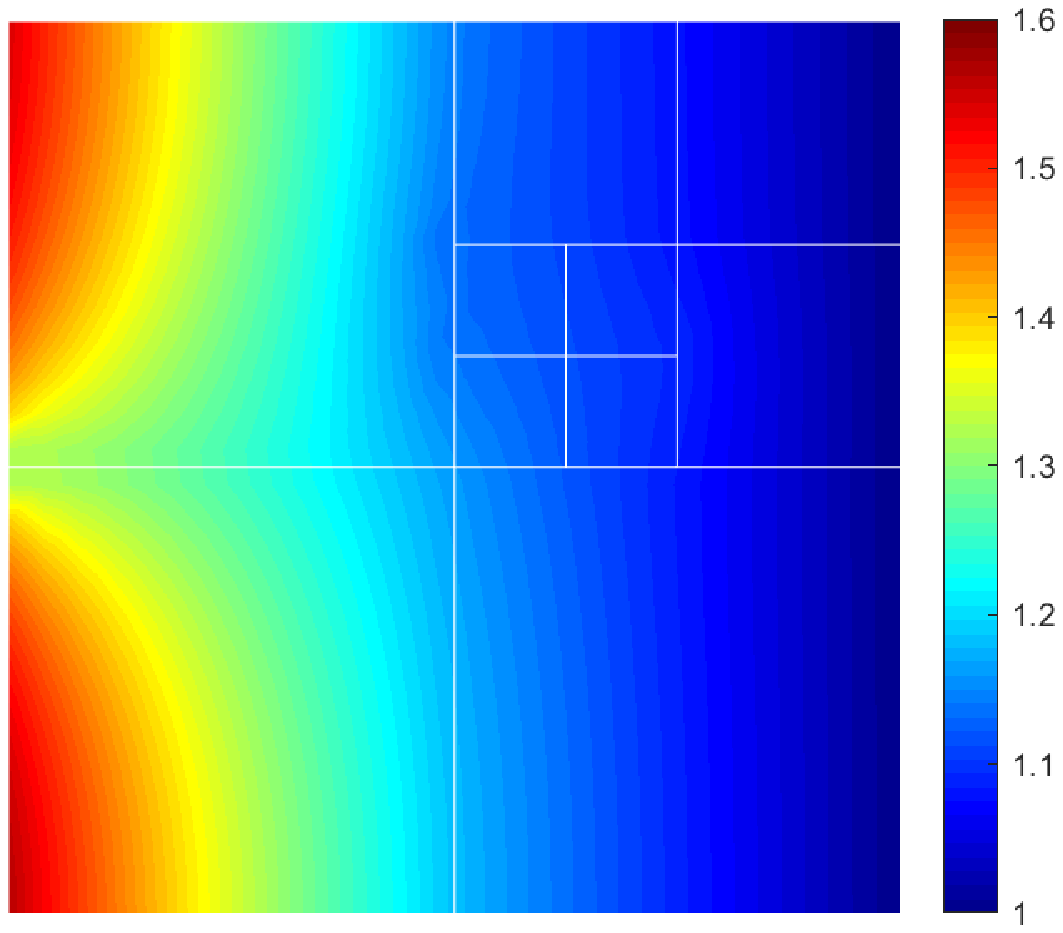}}
\subfigure[Solution of case (a) on $35\times35$ mesh]{\includegraphics[width = 3in]{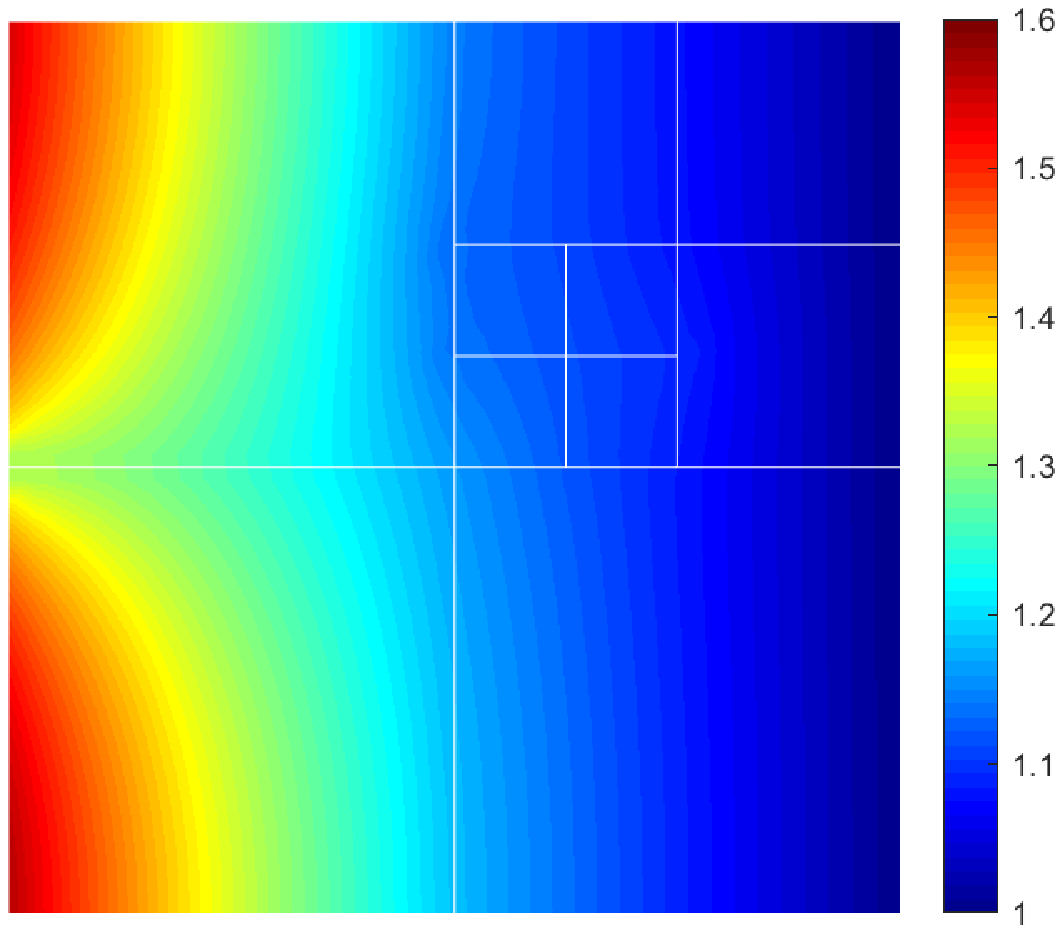}} \\
\subfigure[Solution of case (b) on $25\times25$ mesh]{\includegraphics[width = 3in]{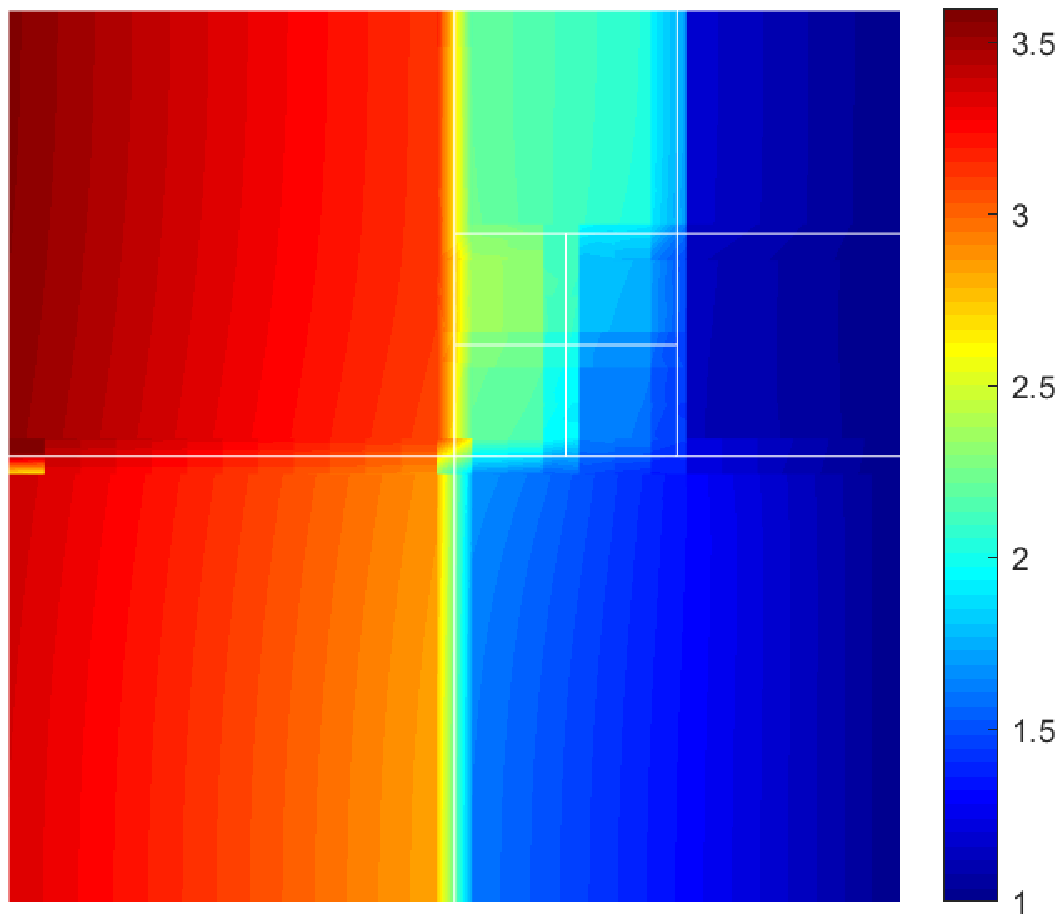}}
\subfigure[Solution of case (b) on $35\times35$ mesh]{\includegraphics[width = 3in]{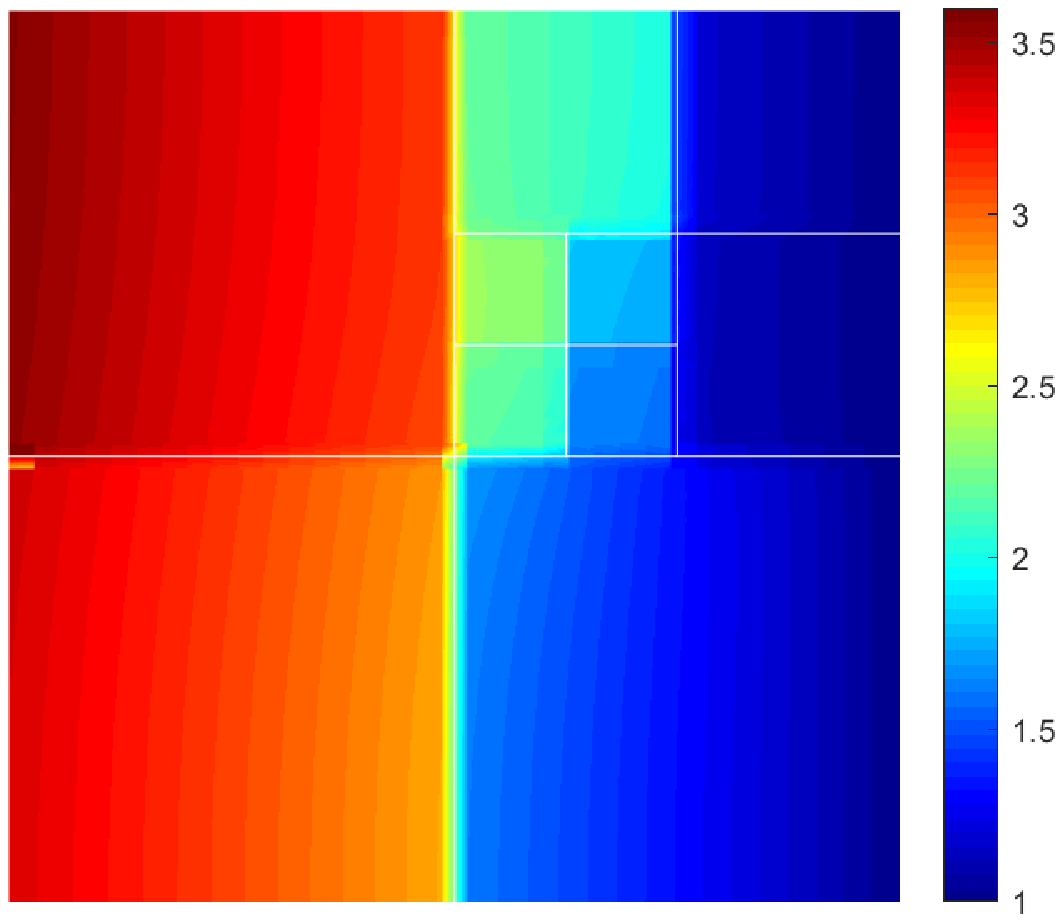}} \\
\caption{Simulation results of Example \ref{ex3} on non-conforming meshes}\label{fig:ex3contour}
\end{figure}

\begin{figure}[!htbp]
\subfigure[Pressure along $y=0.7$ on $25\times25$ mesh]{\includegraphics[width = 3in]{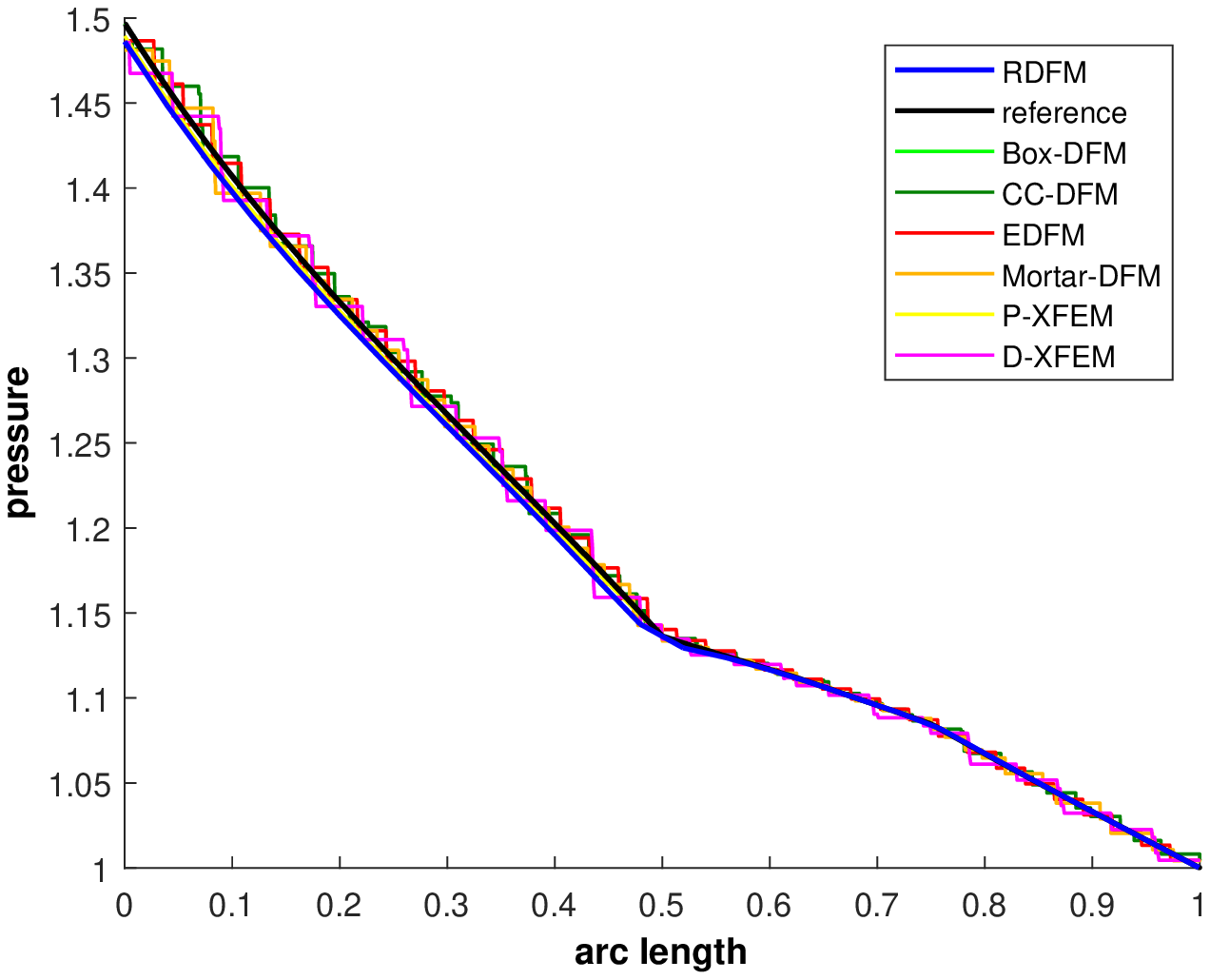}}
\subfigure[Pressure along $y=0.7$ on $35\times35$ mesh]{\includegraphics[width = 3in]{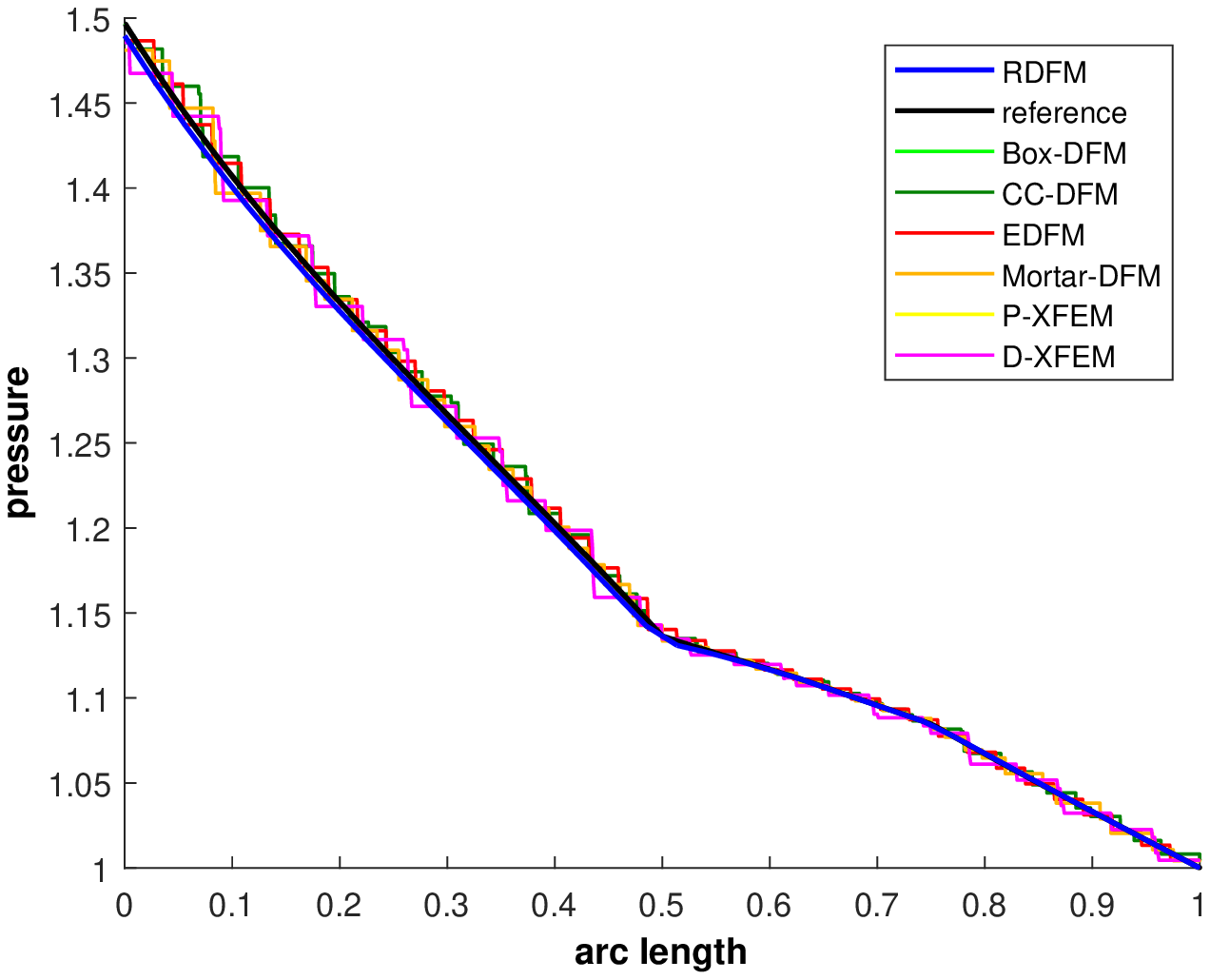}} \\
\subfigure[Pressure along $x=0.5$ on $25\times25$ mesh]{\includegraphics[width = 3in]{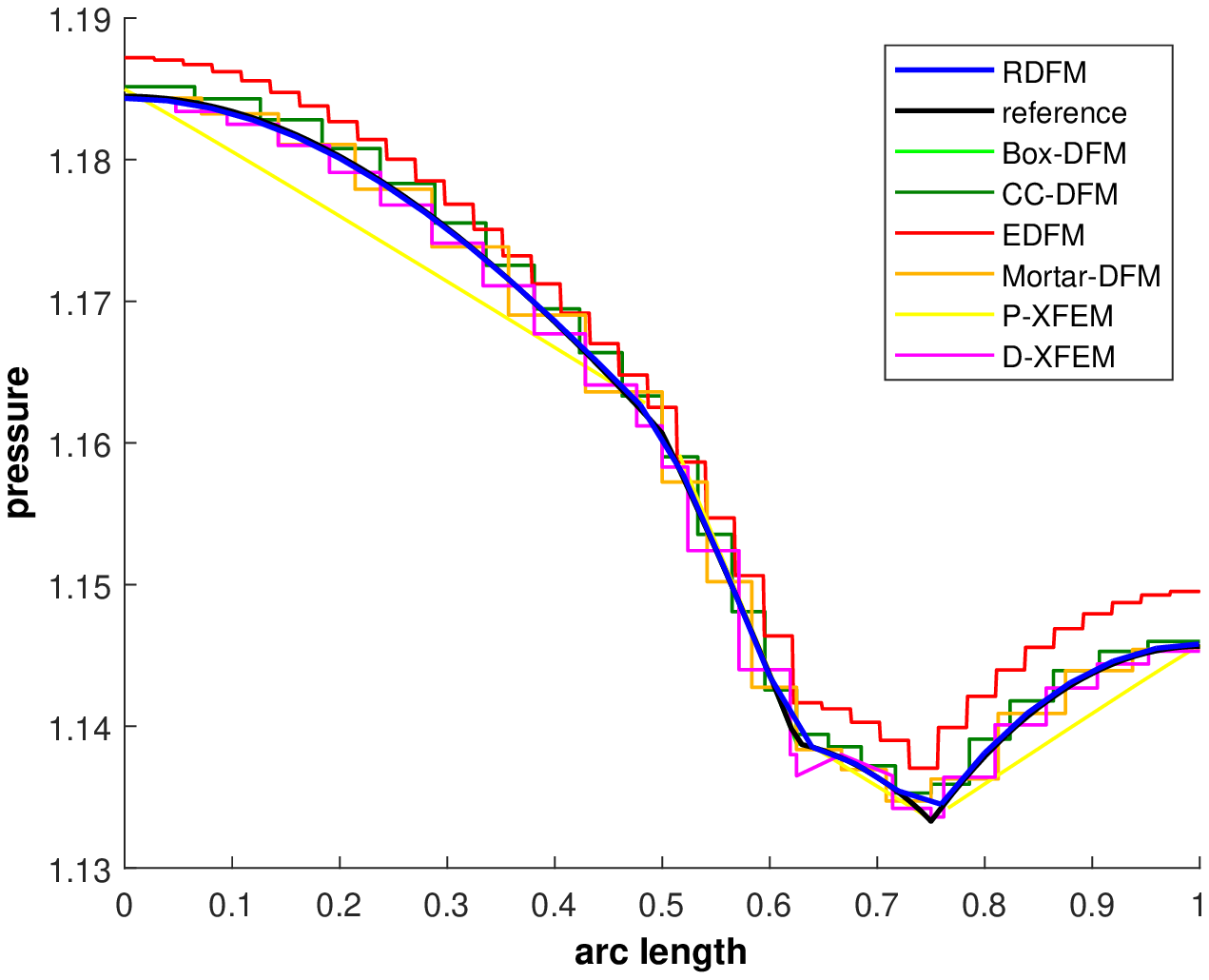}}
\subfigure[Pressure along $x=0.5$ on $35\times35$ mesh]{\includegraphics[width = 3in]{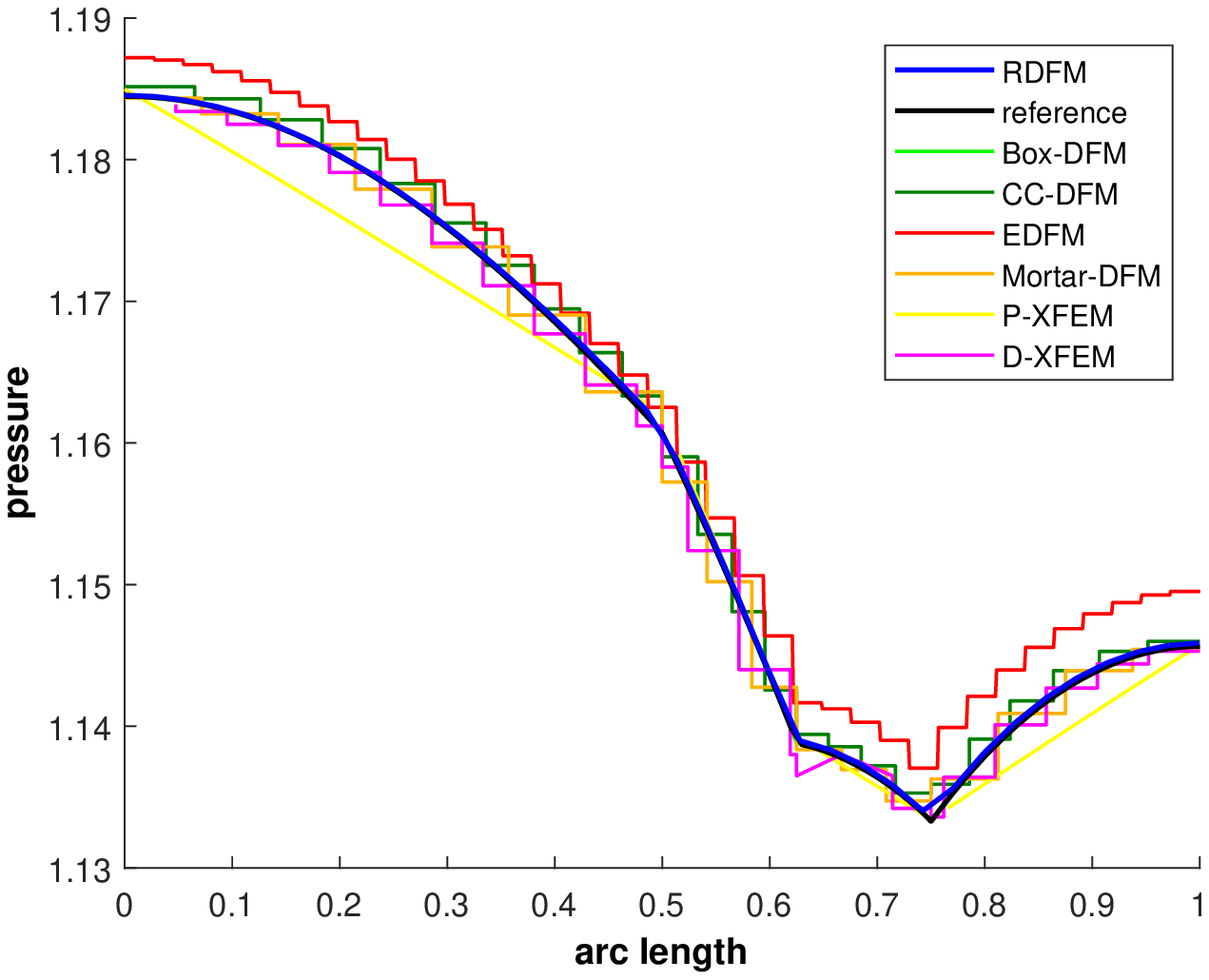}} \\
\caption{Slices of pressure in case (a) of Example \ref{ex3} on non-conforming meshes}\label{fig:ex3SlicesCaseA}
\end{figure}

\begin{figure}[!htbp]
\subfigure[Pressure along $(0.0, 0.1)-(0.9, 1.0)$ on $25\times25$ mesh]{\includegraphics[width = 3in]{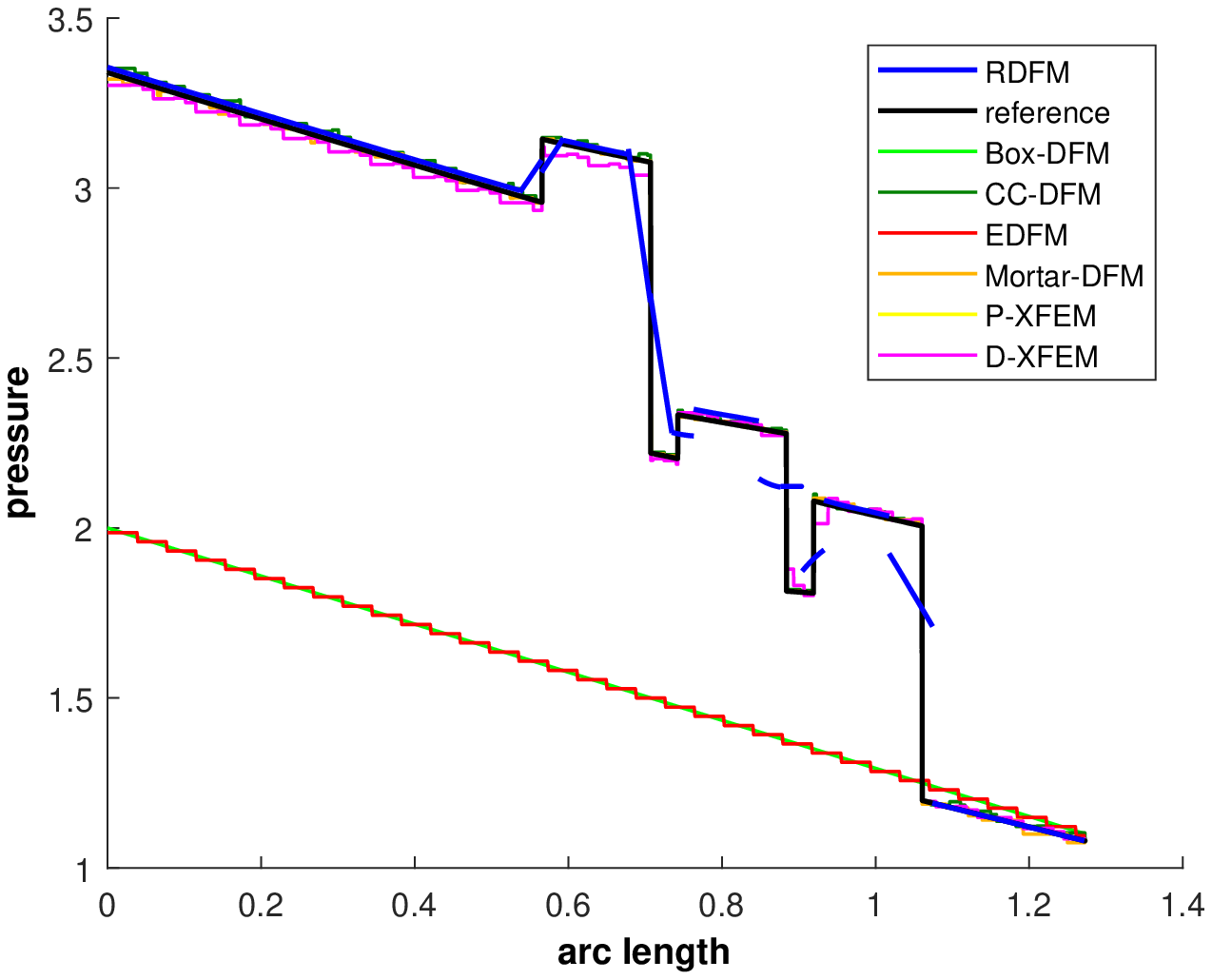}}
\subfigure[Pressure along $(0.0, 0.1)-(0.9, 1.0)$ on $35\times35$ mesh]{\includegraphics[width = 3in]{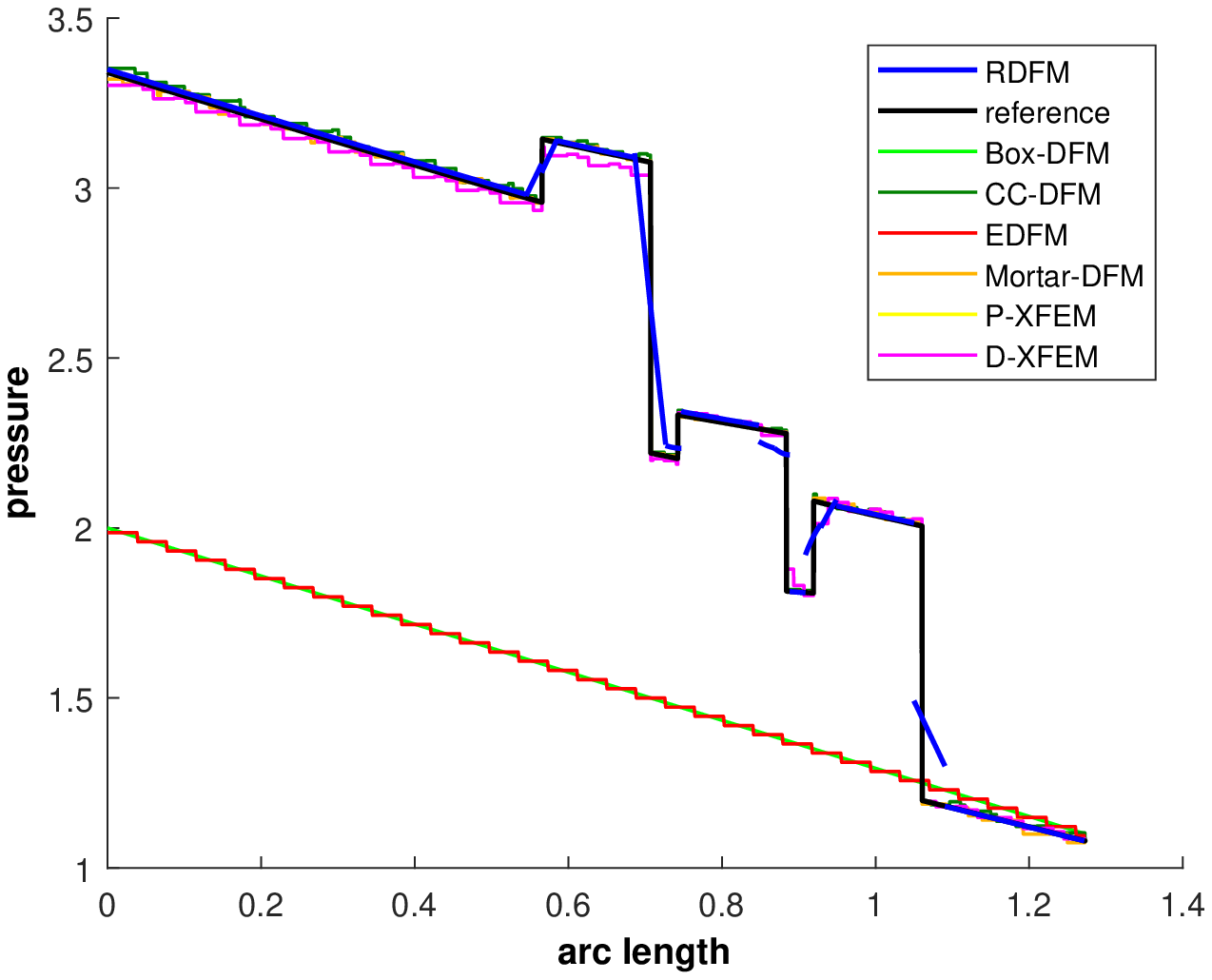}} \\
\caption{Slices of pressure in case (b) of Example \ref{ex3} on non-conforming meshes}\label{fig:ex3SlicesCaseB}
\end{figure}

\begin{figure}[!htbp]
\subfigure[Solution of case (a) on $8\times8$ mesh]{\includegraphics[width = 3in]{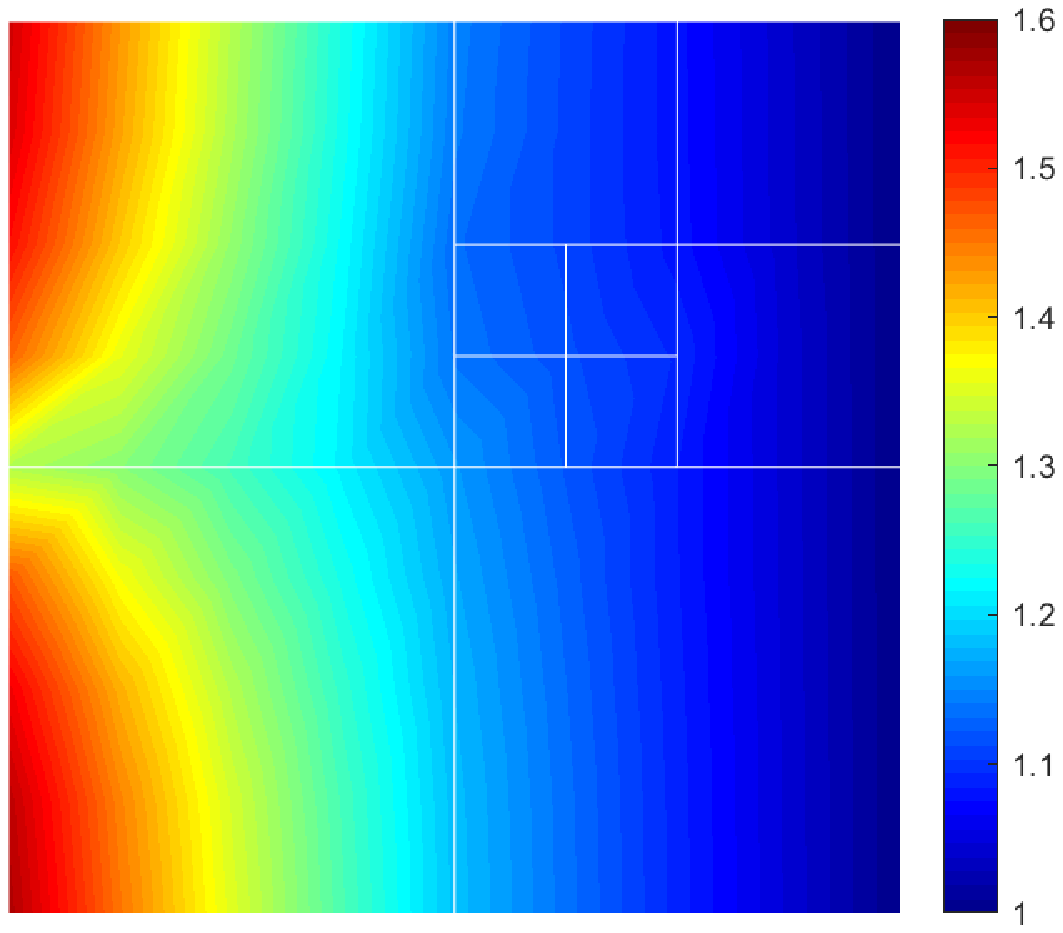}}
\subfigure[Solution of case (a) on $16\times16$ mesh]{\includegraphics[width = 3in]{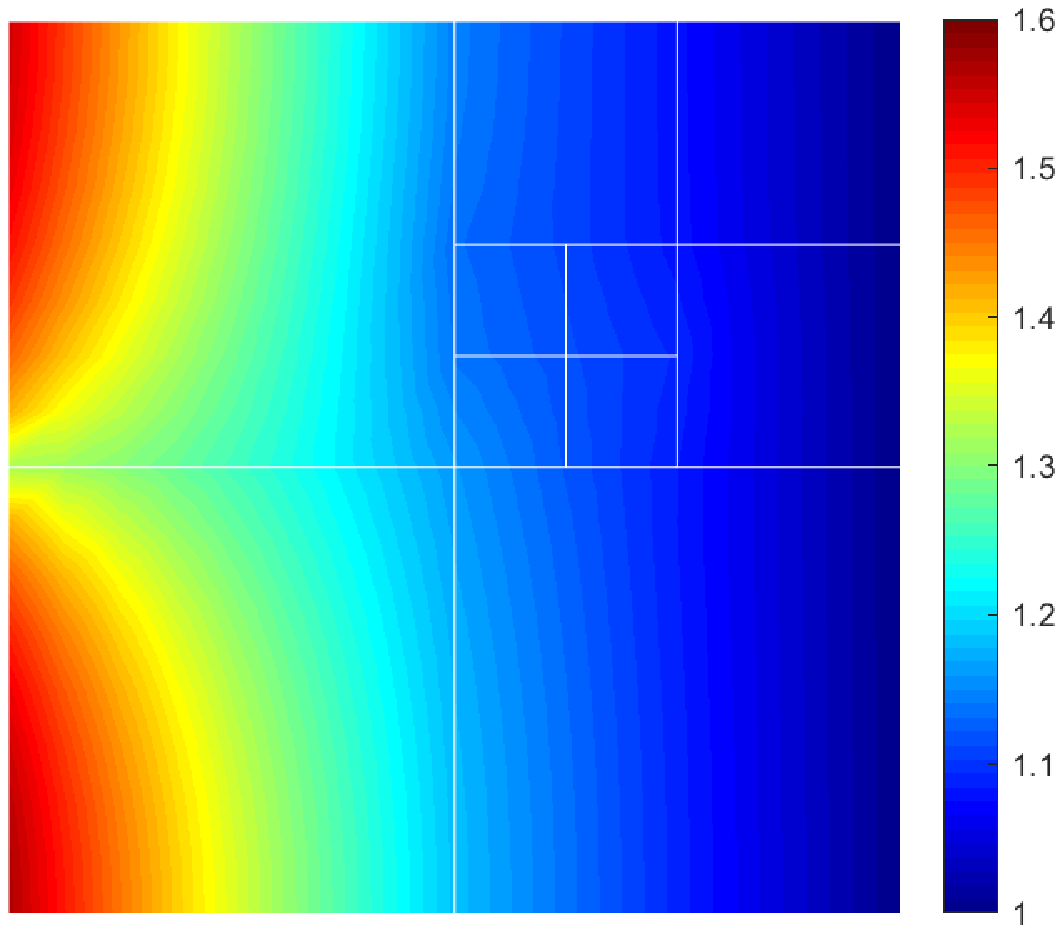}} \\
\subfigure[Solution of case (b) on $8\times8$ mesh]{\includegraphics[width = 3in]{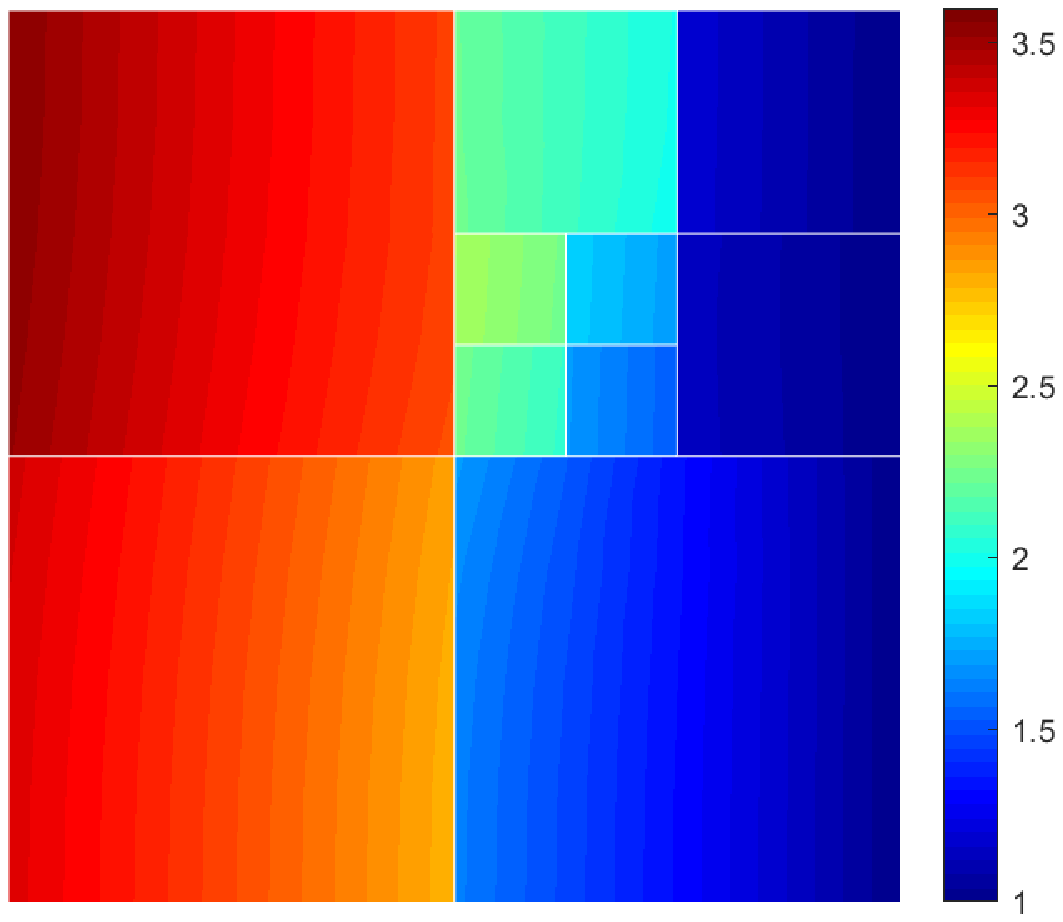}}
\subfigure[Solution of case (b) on $16\times16$ mesh]{\includegraphics[width = 3in]{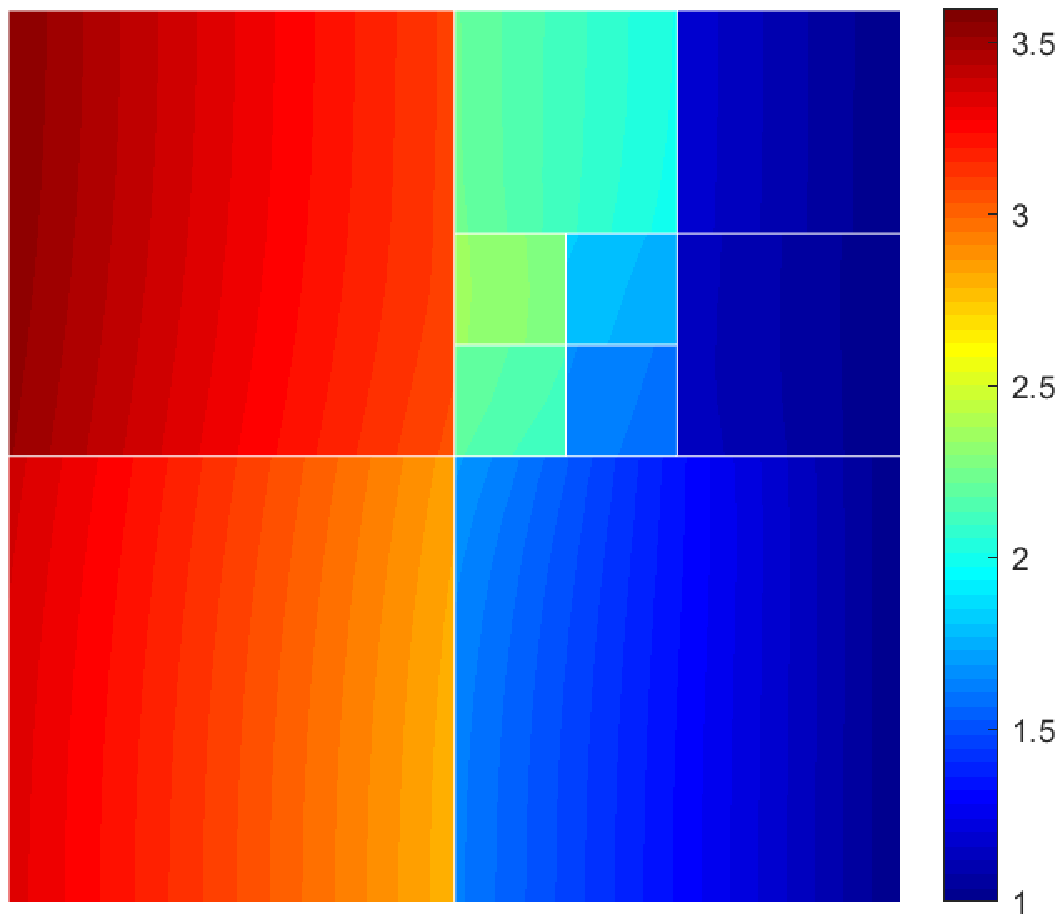}} \\
\caption{Simulation results of Example \ref{ex3} on conforming meshes}\label{fig:ex3contour2}
\end{figure}

\begin{figure}[!htbp]
\subfigure[Pressure along $y=0.7$ on $8\times8$ mesh]{\includegraphics[width = 3in]{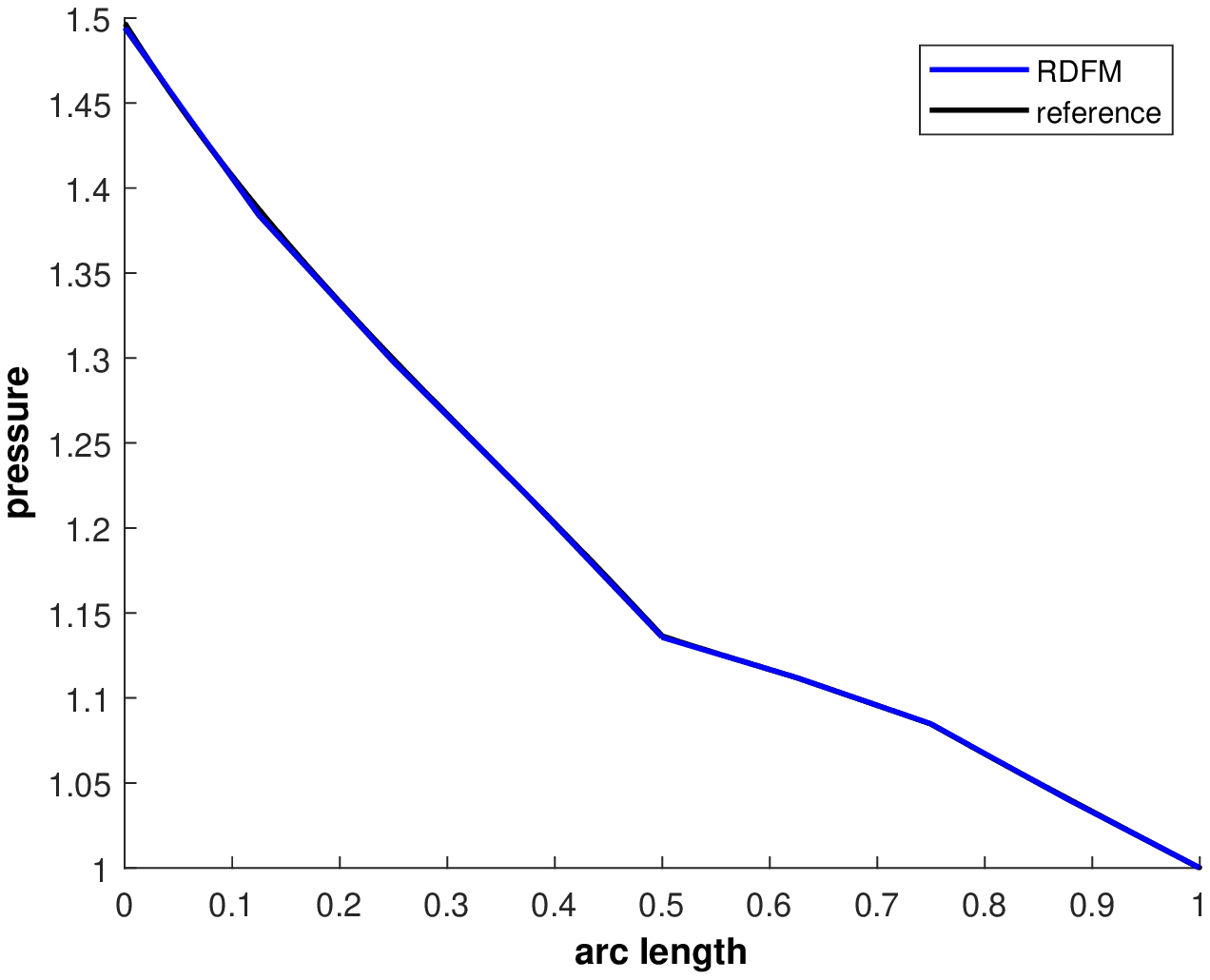}}
\subfigure[Pressure along $y=0.7$ on $16\times16$ mesh]{\includegraphics[width = 3in]{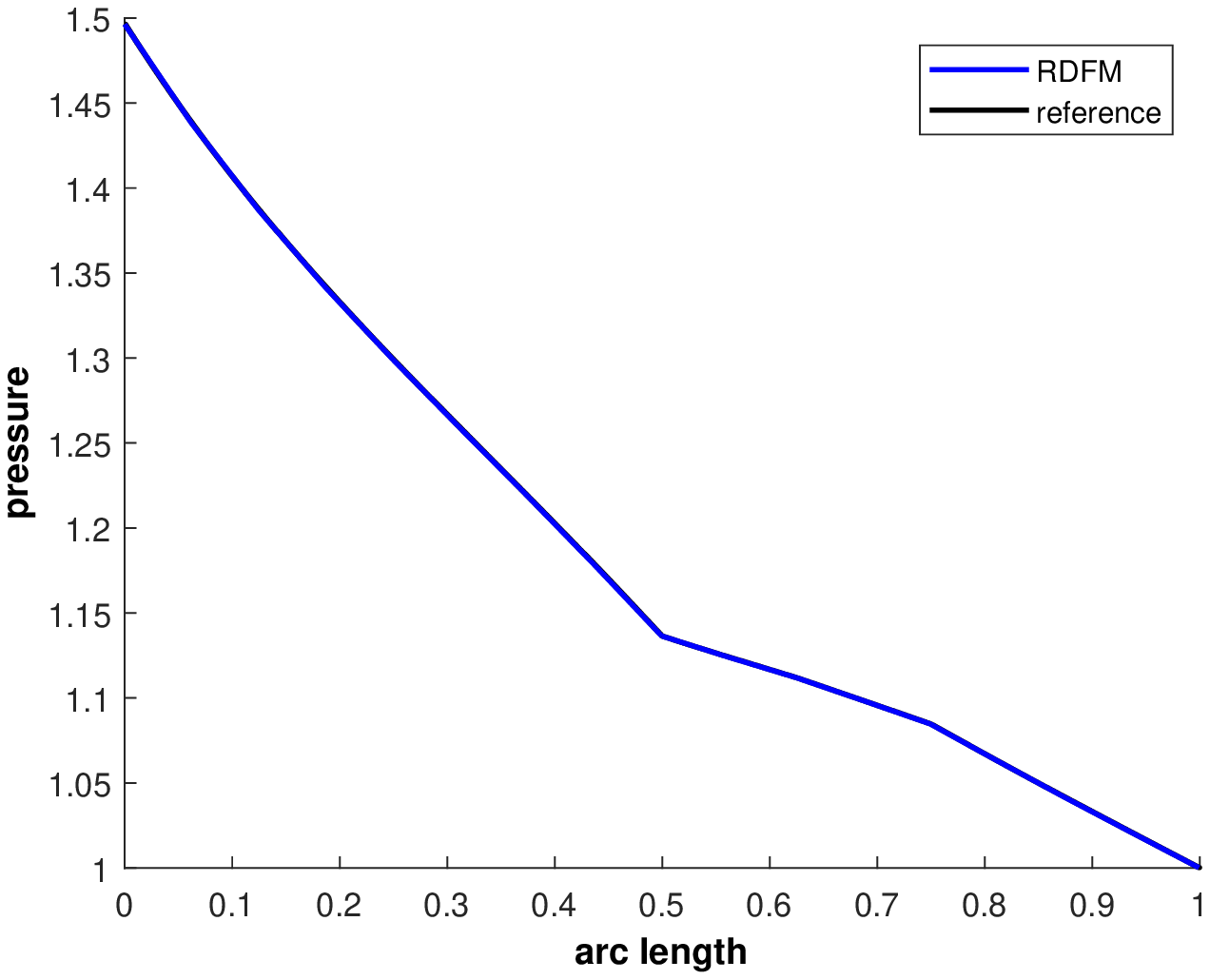}} \\
\subfigure[Pressure along $x=0.5$ on $8\times8$ mesh]{\includegraphics[width = 3in]{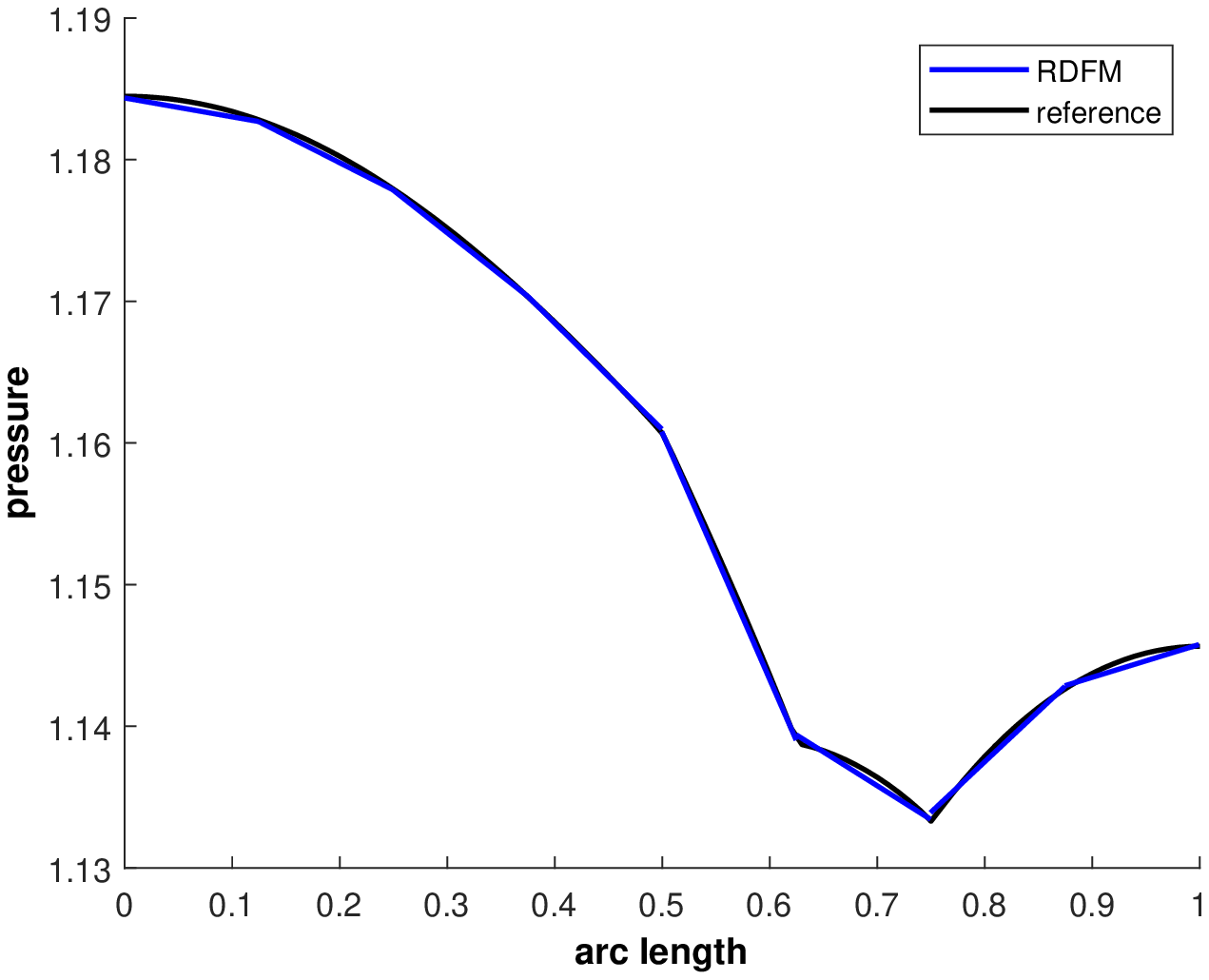}}
\subfigure[Pressure along $x=0.5$ on $16\times16$ mesh]{\includegraphics[width = 3in]{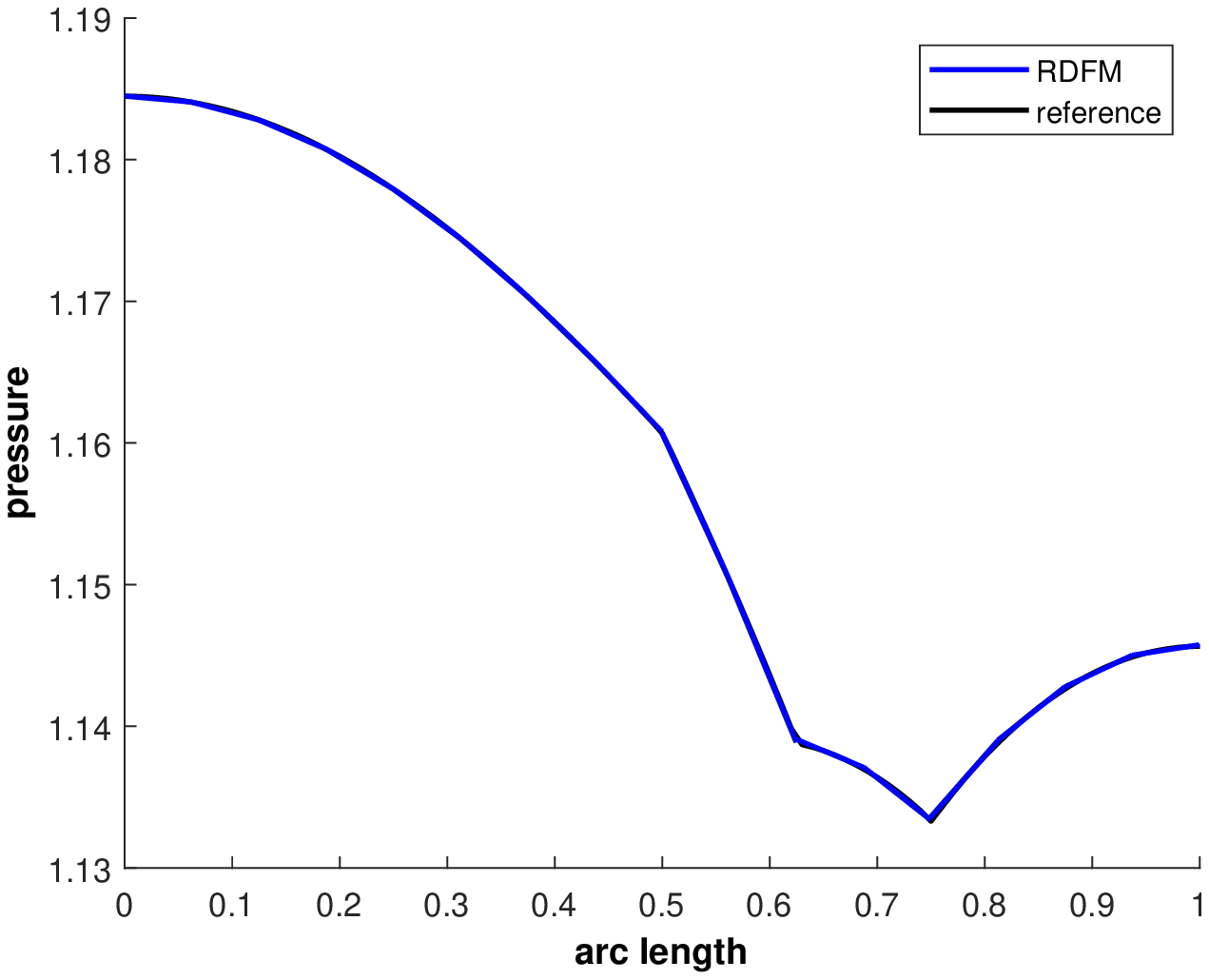}} \\
\caption{Slices of pressure in case (a) of Example \ref{ex3} on conforming meshes}\label{fig:ex3SlicesCaseAConforming}
\end{figure}

\begin{figure}[!htbp]
\subfigure[Pressure along $(0.0, 0.1)-(0.9, 1.0)$ on $8\times8$ mesh]{\includegraphics[width = 3in]{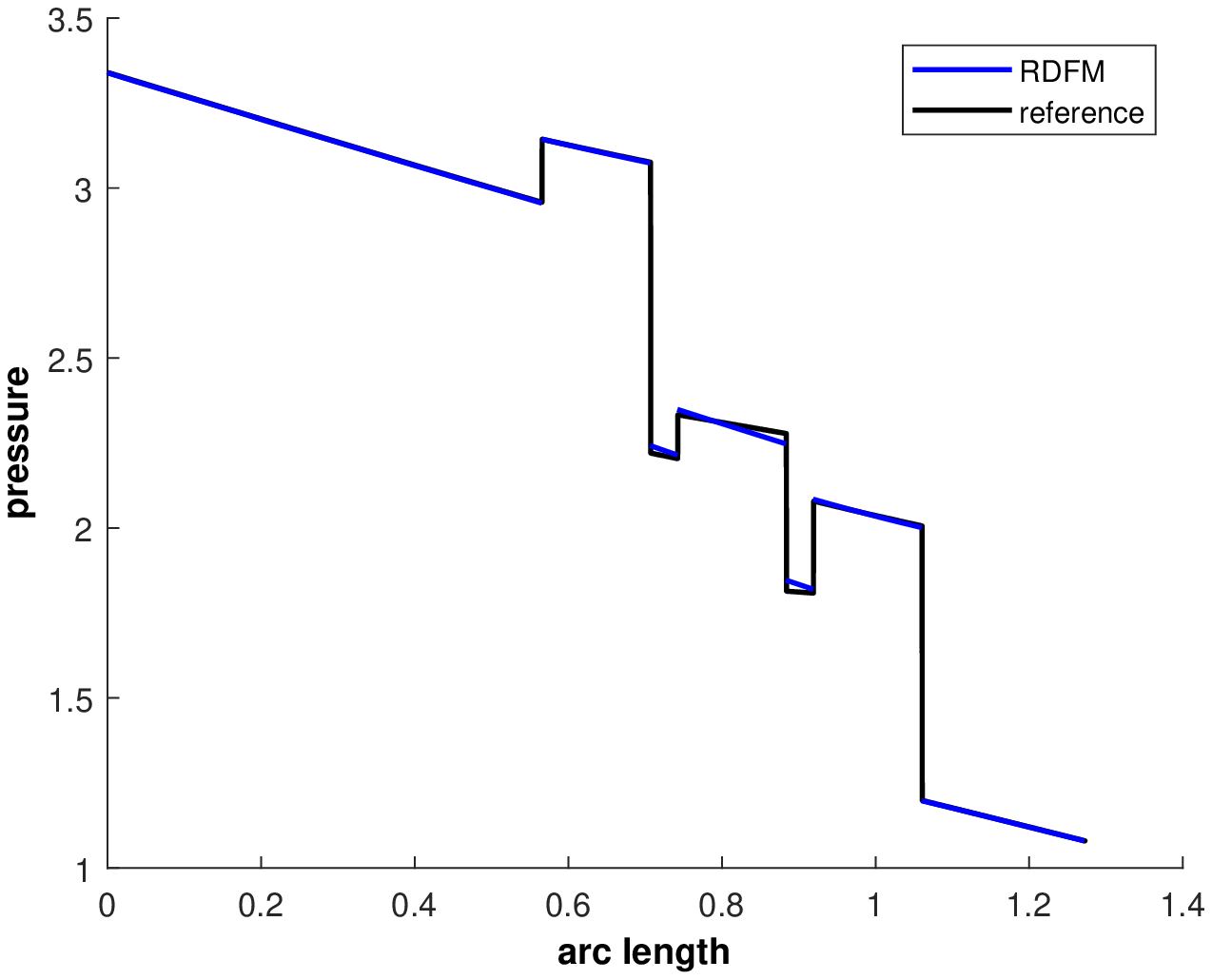}}
\subfigure[Pressure along $(0.0, 0.1)-(0.9, 1.0)$ on $16\times16$ mesh]{\includegraphics[width = 3in]{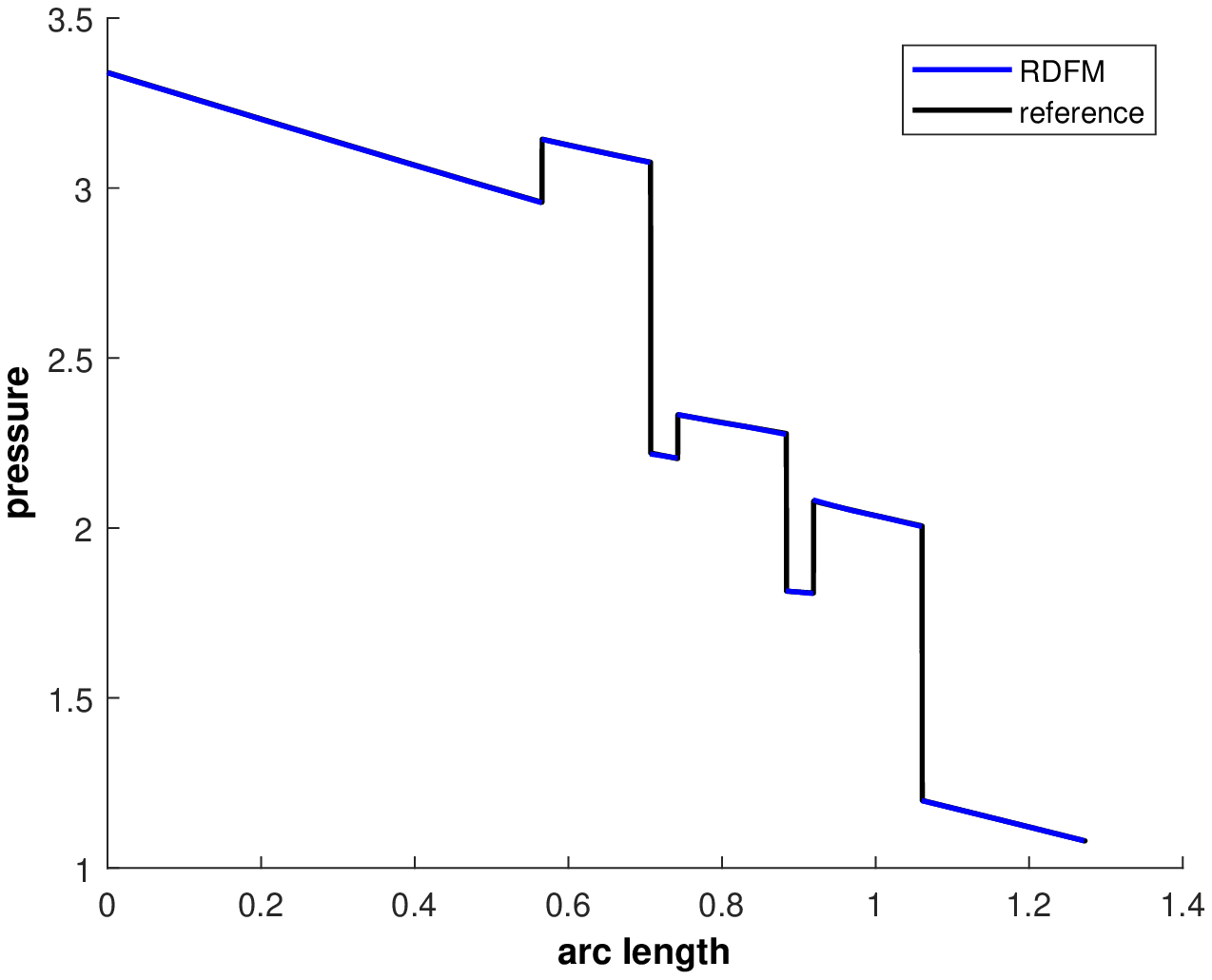}} \\
\caption{Slices of pressure in case (b) of Example \ref{ex3} on conforming meshes}\label{fig:ex3SlicesCaseBConforming}
\end{figure}

\begin{table}[!htbp]
\centering
\caption{\label{tab:ex3_CaseA} Evaluation data of different methods of case (a) of Example \ref{ex3}}
\begin{tabular}{|c|c|c|c|c|c|}
  \hline
method & matrix elements & fracture elements & mesh & err$_m$ &err$_f$ \\\hline
Box-DFM&$2691$ triangles & $130$ & conforming       &6.7e-03 & 1.1e-03  \\\hline
CC-DFM&$1386$ triangles & $95$ & conforming & 1.1e-02 & 5.0e-03 \\\hline
EDFM&$1369$ rectangles & $132$ & non-conforming & 6.5e-03 & 4.0e-03 \\\hline
Mortar-DFM&$1280$ triangles & $75$ & conforming  & 1.0e-02 & 7.4e-03 \\ \hline
P-XFEM&$961$ rectangles & $318$ & non-conforming  & 1.7e-02 & 6.0e-03 \\ \hline
D-XFEM&$1250$ triangles & $126$ & non-conforming  & 9.6e-03 & 8.9e-03 \\ \hline
RDFM& $625$ rectangles & $90$ &non-conforming & 8.6e-03 & 8.8e-03 \\ \hline
RDFM& $1225$ rectangles & $126$ &non-conforming & 5.9e-03 & 6.2e-03 \\ \hline
RDFM& $64$ rectangles & $28$ &conforming & 3.4e-03 &2.7e-02  \\ \hline
RDFM& $256$ rectangles & $56$ &conforming & 1.5e-03 & 1.4e-02 \\ \hline
\end{tabular}
\end{table}

\begin{table}[!htbp]
\centering
\caption{\label{tab:ex3_CaseB} Evaluation data of different methods of case (b) of Example \ref{ex3}}
\begin{tabular}{|c|c|c|c|c|c|}
\hline
method & matrix elements & fracture elements & mesh & err$_m$ &err$_f$ \\\hline
Box-DFM&$2691$ triangles & $130$ & conforming       &4.1e-01 & 3.8e-01  \\\hline
CC-DFM&$1386$ triangles & $95$ & conforming & 5.7e-03 & 4.6e-03 \\\hline
EDFM&$1369$ rectangles & $132$ & non-conforming & 2.9e-01 & 3.2e-01 \\\hline
Mortar-DFM&$1280$ triangles & $75$ & conforming  & 4.5e-03 & 4.9e-03 \\ \hline
P-XFEM&$961$ rectangles & $318$ & non-conforming  & 2.9e-03 & 2.2e-02 \\ \hline
D-XFEM&$1250$ triangles & $126$ & non-conforming  & 1.0e-02 & 1.9e-02 \\ \hline
RDFM& $625$ rectangles & $90$ &non-conforming & 3.1e-02 & 4.0e-02 \\ \hline
RDFM& $1225$ rectangles & $126$ &non-conforming & 2.6e-02 & 4.2e-02 \\ \hline
RDFM& $64$ rectangles & $28$ &conforming & 2.6e-03 & 1.2e-01 \\ \hline
RDFM& $256$ rectangles & $56$ &conforming & 5.9e-04 & 1.1e-01 \\ \hline
\end{tabular}
\end{table}

The reference solutions are provided by \cite{Benchmark} based on the mimetic finite difference (MFD) method on a very fine mesh containing $1136456$ matrix elements and $38600$ fracture elements.
The methods participate in the comparison are the vertex-centered control volume discrete fracture model (Box-DFM), the cell-centered two
point flux approximation control volume discrete fracture model (CC-DFM), the embedded discrete fracture model (EDFM), the mortar-flux discrete fracture model (Mortar-DFM), the primal extended finite element method (P-XFEM), and the dual extended finite element method (D-XFEM), refer to \cite{Benchmark} for a detailed introduction.

We first test our model on non-conforming meshes.
To demonstrate its performance comprehensively, we conduct the computation on two rectangular meshes with different grid sizes.
The coarse mesh is $25\times25$ while the fine mesh is $35\times35$.
We plot the contours of pressure for both cases in Figure \ref{fig:ex3contour}.
One can find the contours of reference solutions in \cite{Benchmark} for comparison.
We also slice the profiles of the pressure along the line $y=0.7$ and $x=0.5$ for case (a) and the line $(0.0, 0.1)-(0.9, 1.0)$ for case (b), and draw the comparison with the reference solutions and other models in Figure \ref{fig:ex3SlicesCaseA} and Figure \ref{fig:ex3SlicesCaseB}, respectively.

Furthermore, we compute the model on $8\times8$ and $16\times16$ uniform rectangular meshes. Note that the meshes are conforming in this case.
We plot the contours of pressure for both cases in Figure \ref{fig:ex3contour2}.
The slices of the solutions along the line $y=0.7$ and $x=0.5$ for case (a) and the line $(0.0, 0.1)-(0.9, 1.0)$ for case (b) are shown in Figure \ref{fig:ex3SlicesCaseAConforming} and Figure \ref{fig:ex3SlicesCaseBConforming} respectively, together with those of reference solutions.
We don't draw slices of other methods again since they are already shown in the previous comparisons.


A quantitative comparison of relative errors, together with other important aspects of different models, are summarized in the Table \ref{tab:ex3_CaseA} and Table \ref{tab:ex3_CaseB} for case (a) and case (b), respectively.
The relative errors on matrix and fractures/barriers of a solution are defined as:
\begin{equation}\label{ErrmFormula}
    \text{err}^{2}_\text{m} = \frac{1}{|\Omega|(\Delta p_\text{ref})^{2}} \sum_{i,j} |T^{i}_\text{ref} \cap T^{j}_\text{m}|\left(p_\text{ref}|_{T^{i}_\text{ref}}-p_\text{m}  |_{T^{j}_\text{m}}\right)^{2},
\end{equation}
\begin{equation}\label{ErrfFormula}
    \text{err}^{2}_\text{f} = \frac{1}{|\Gamma|(\Delta p_\text{ref})^{2}} \sum_{i,l} |T^{i}_\text{ref} \cap T^{l}_\text{f}|\left(p_\text{ref}|_{T^{i}_\text{ref}}-p_\text{f}  |_{T^{l}_\text{f}}\right)^{2},
\end{equation}
where $|\Omega|$ and $|\Gamma|$ are measures of the matrix and fractures respectively, $p_\text{ref}$ is the reference solution, $\Delta p_\text{ref}=\max{p_\text{ref}}-\min{p_\text{ref}}$ is the range of the reference solution, $T^{i}_\text{ref} , i=1,2,\ldots,I,$ are cells employed for the reference solution, $T^{j}_\text{m}, j=1,2,\ldots,J,$ and
$T^{l}_\text{f}, l=1,2,\ldots,L,$ are the matrix elements and fracture elements employed for the method to be evaluated, respectively, see \cite{Benchmark} for more details about the definition.

We want to point out that, due to the discontinuity of pressure across barriers, the relative error on barriers is not a good indicator of accuracy.
This concept even doesn't make much sense in our model since we use delta function to represent the barrier.
That's why the relative error of pressure on barriers fails to converge in our method, but this doesn't mean the model is failed.
In fact, as one can see from the slices of solutions along $(0.0, 0.1)-(0.9, 1.0)$ and the relative errors on the porous matrix, our method provides good approximations in barriers network case, especially on the conforming meshes.

\begin{figure}[!htbp]
\subfigure[Case (a): flow from top to bottom]{\includegraphics[width = 3in]{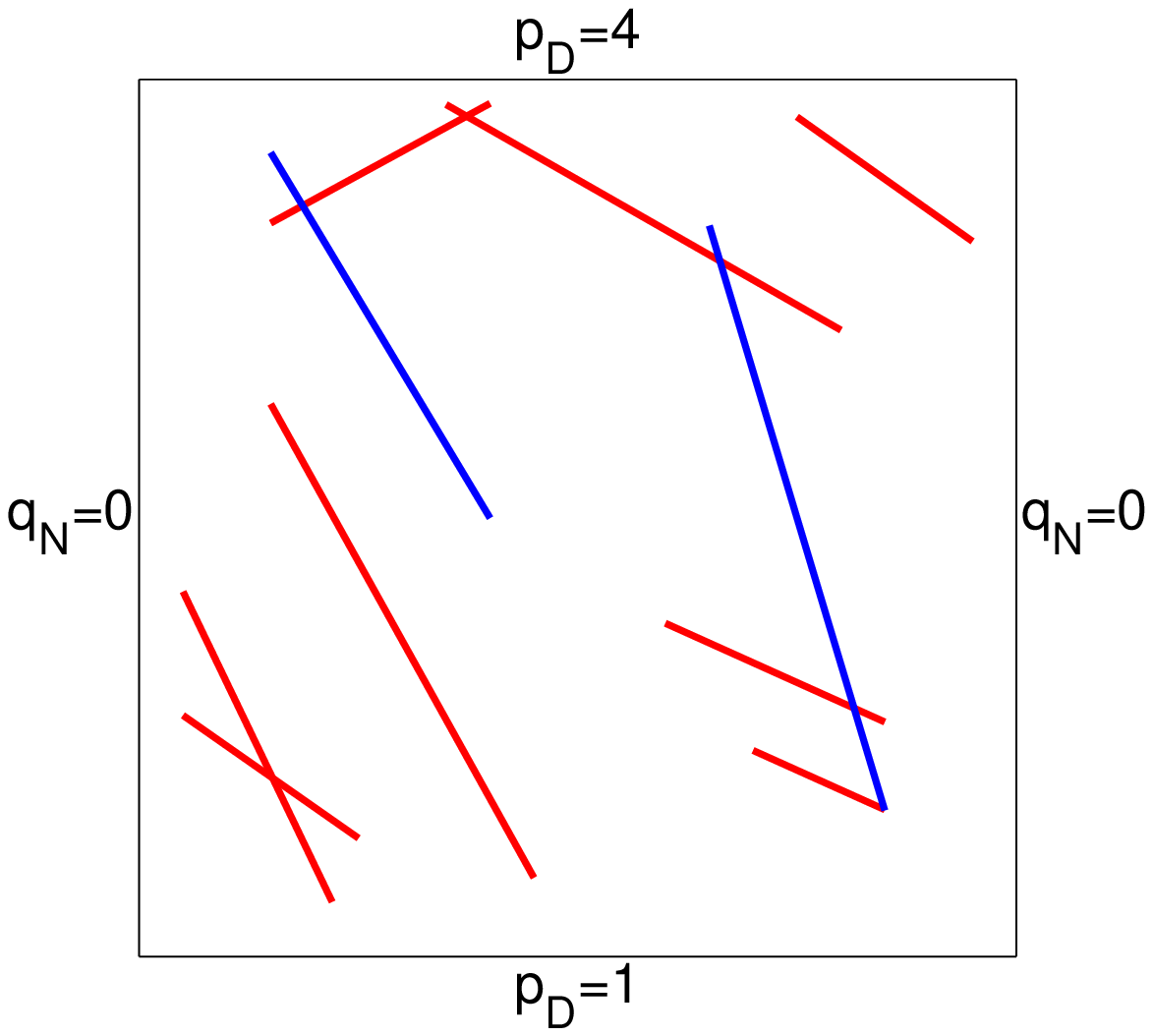}}
\subfigure[Case (b): flow from left to right]{\includegraphics[width = 3in]{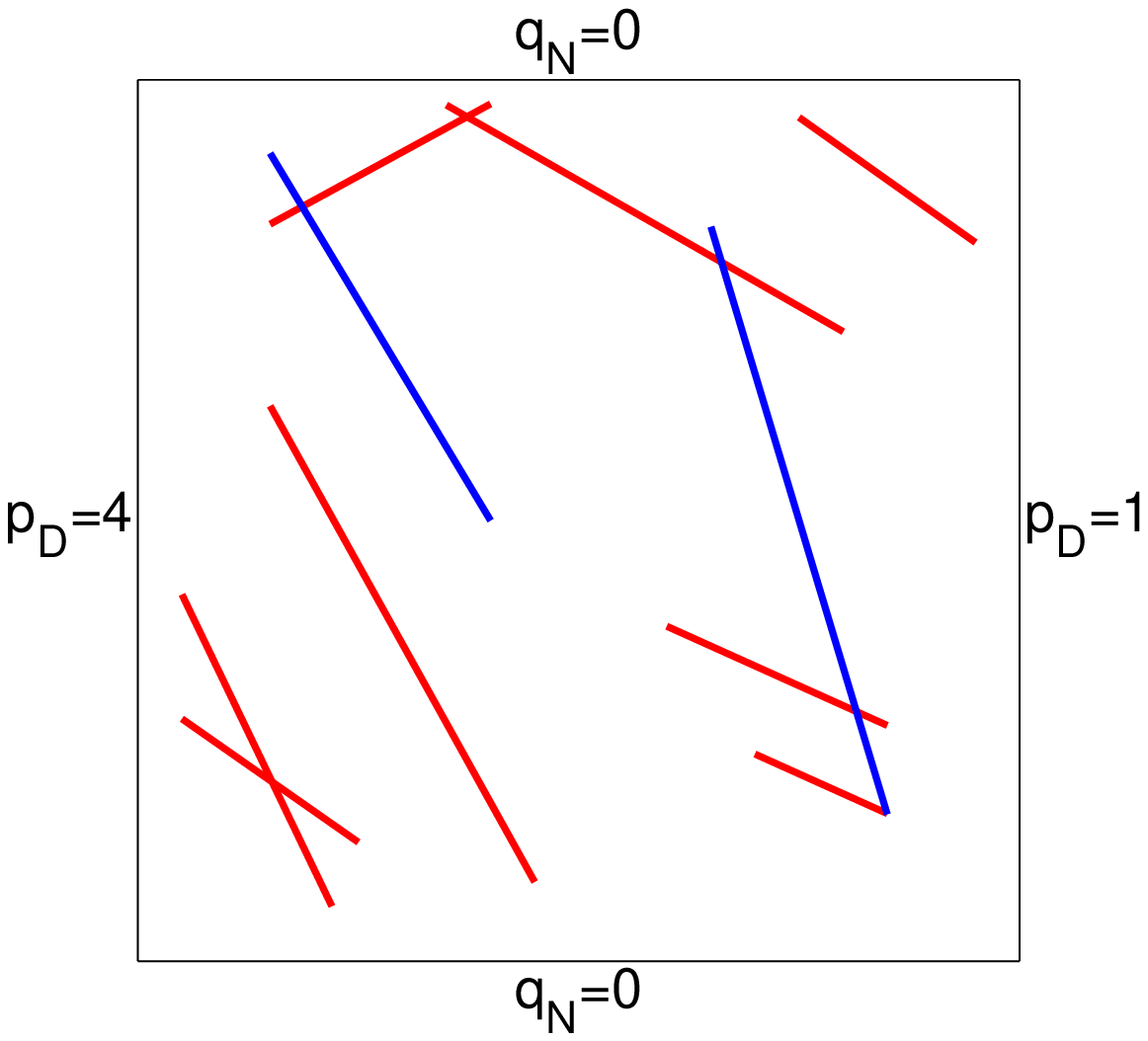}} \\
\caption{Domain and boundary conditions of Example \ref{ex4}}\label{fig:ex4settings}
\end{figure}

\begin{ex} \label{ex4}
\textbf{Complex networks}

In this example, we test complex networks containing $8$ fractures and $2$ barriers.
This example is taken from the benchmark $4.3$ in \cite{Benchmark}.
The distribution of fractures and barriers are shown in Figure \ref{fig:ex4settings}, in which fractures are colored red and barriers are colored blue.
The computational domain is the unit square $\Omega=[0,1]\times[0,1]$ with permeability ${\bf K}_{m}={\bf I}$.
The exact coordinates of the fractures and barriers are attached in the appendix \ref{appd:ex4}.
All fractures and barriers have the uniform thickness $10^{-4}$ with permeability $k_f=10^{4}$ and $k_{\epsilon}=10^{-4}$, respectively.

\textbf{Case (a): Flow from top to bottom} The top and bottom boundaries are Dirichlet with $p_D=4$ and $p_D=1$ respectively. The left and right boundaries are impermeable.

\textbf{Case (b): Flow from left to right} The left and right boundaries are Dirichlet with $p_D=4$ and $p_D=1$ respectively. The top and bottom boundaries are impermeable.
\end{ex}

\begin{figure}[!htbp]
\subfigure[Solution of case (a) on $30\times30$ mesh]{\includegraphics[width = 3in]{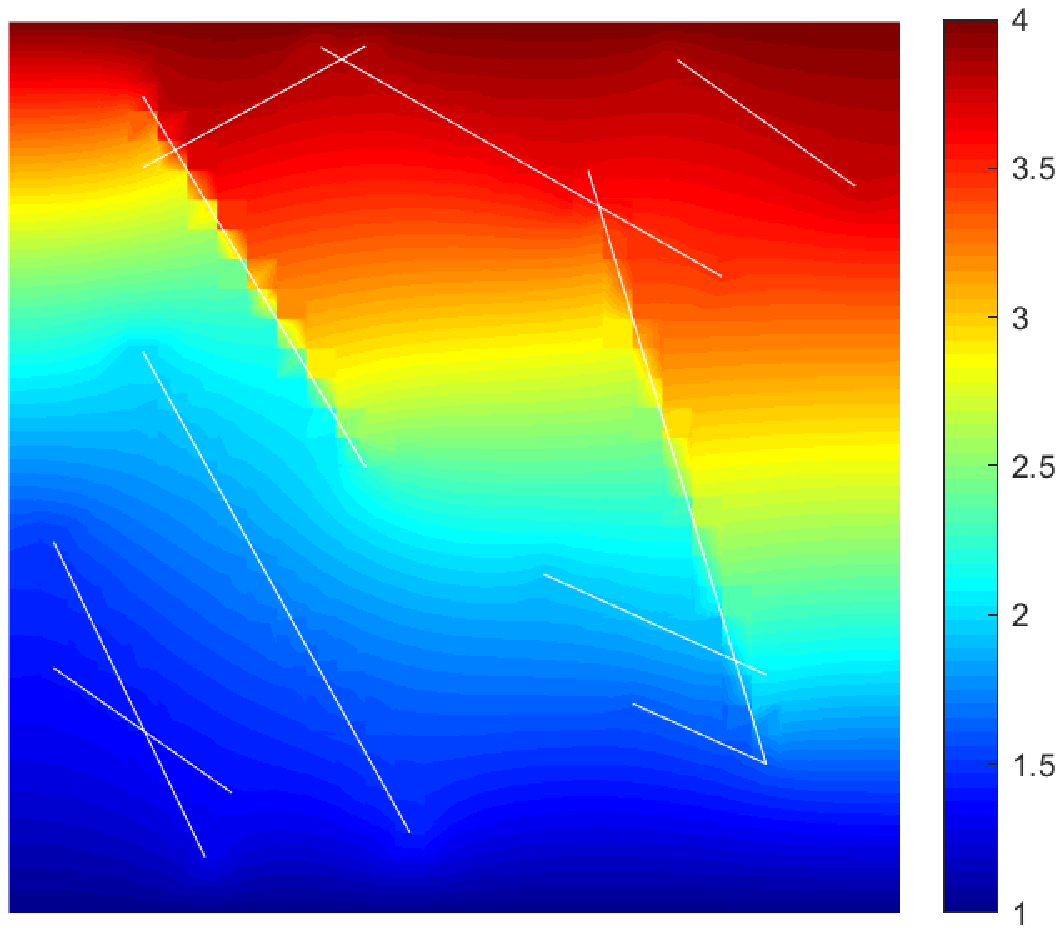}}
\subfigure[Solution of case (a) on $40\times40$ mesh]{\includegraphics[width = 3in]{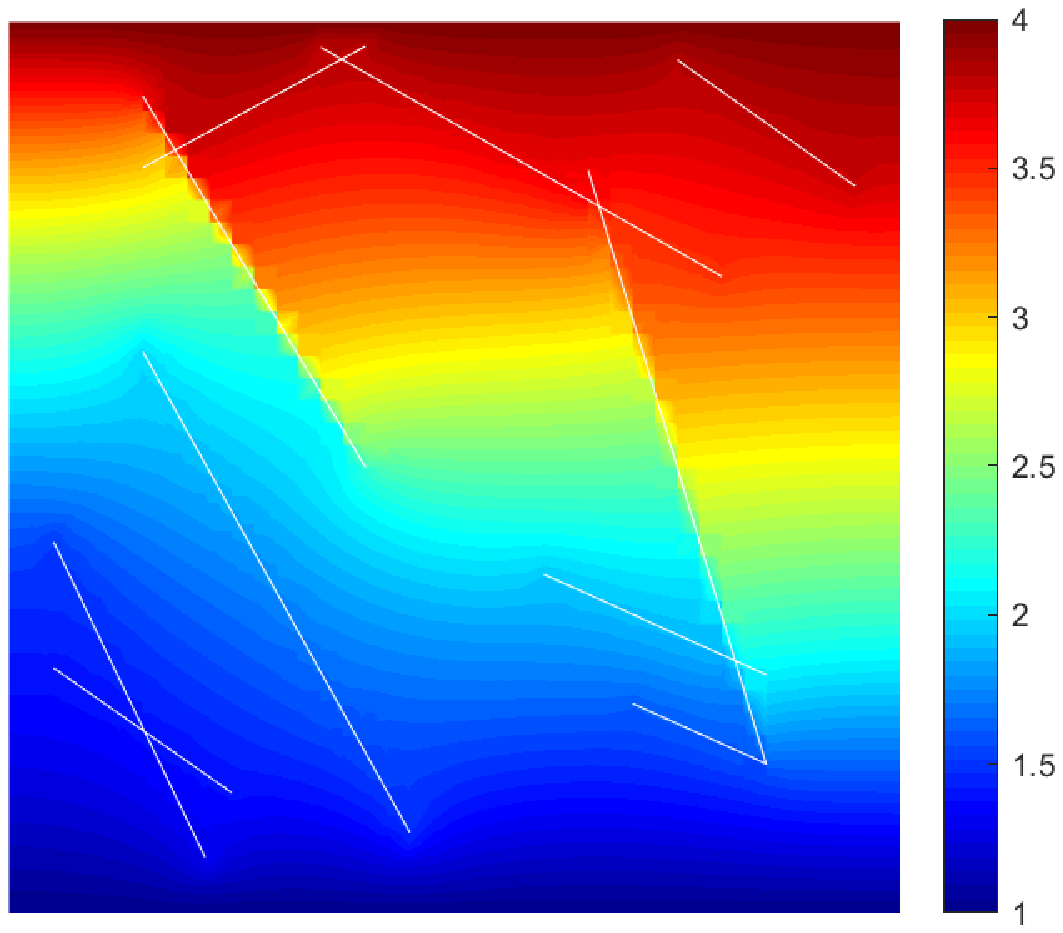}}\\
\subfigure[Solution of case (b) on $30\times30$ mesh]{\includegraphics[width = 3in]{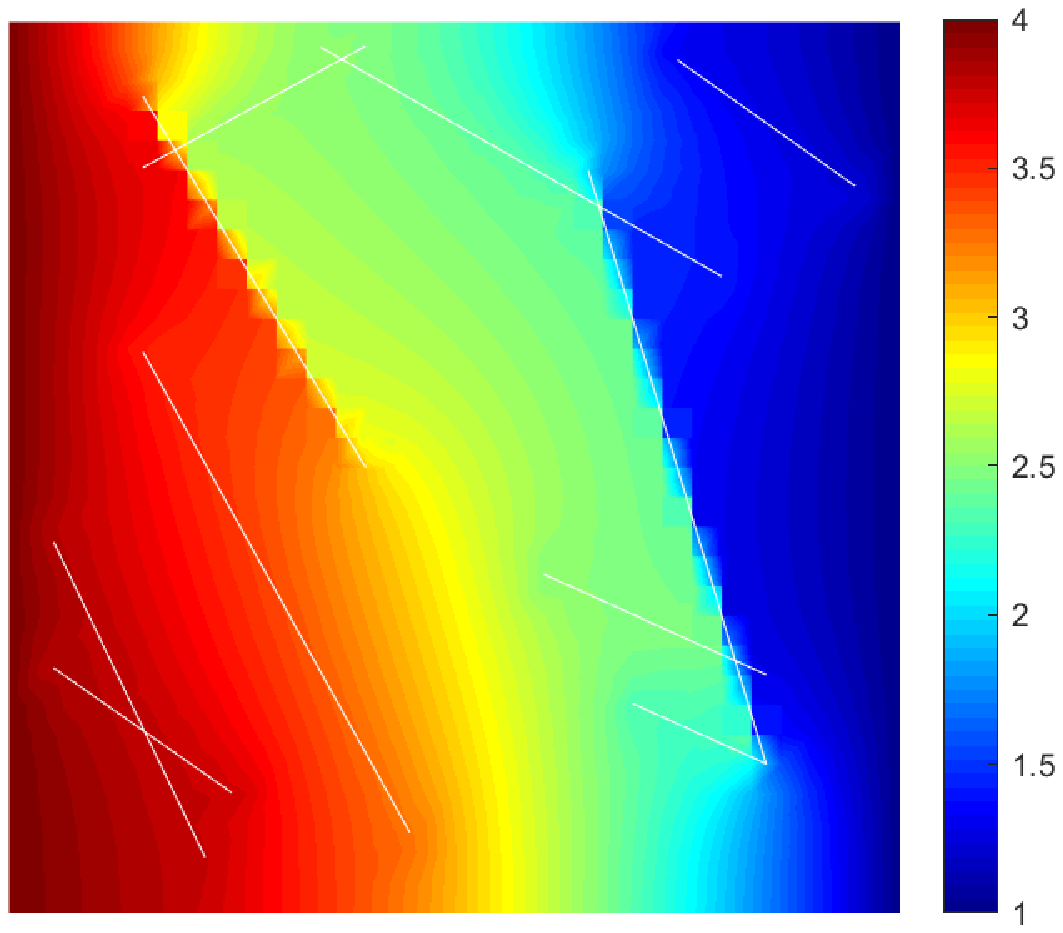}}
\subfigure[Solution of case (b) on $40\times40$ mesh]{\includegraphics[width = 3in]{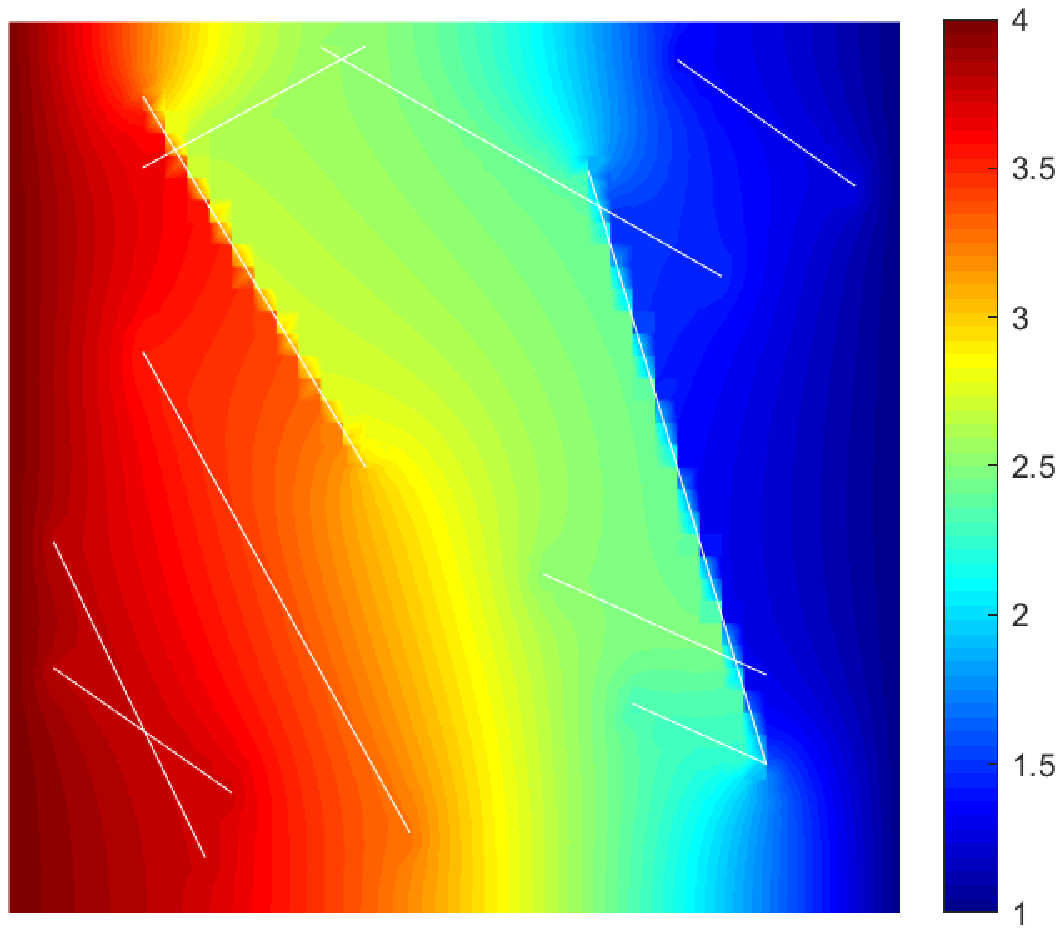}} \\
\caption{Simulation results of Example \ref{ex4} on non-conforming meshes}\label{fig:ex4contour}
\end{figure}

\begin{figure}[!tb]
\subfigure[Pressure along $(0.0, 0.5)-(1.0, 0.9)$ on $30\times30$ mesh in case (a)]{\includegraphics[width = 3in]{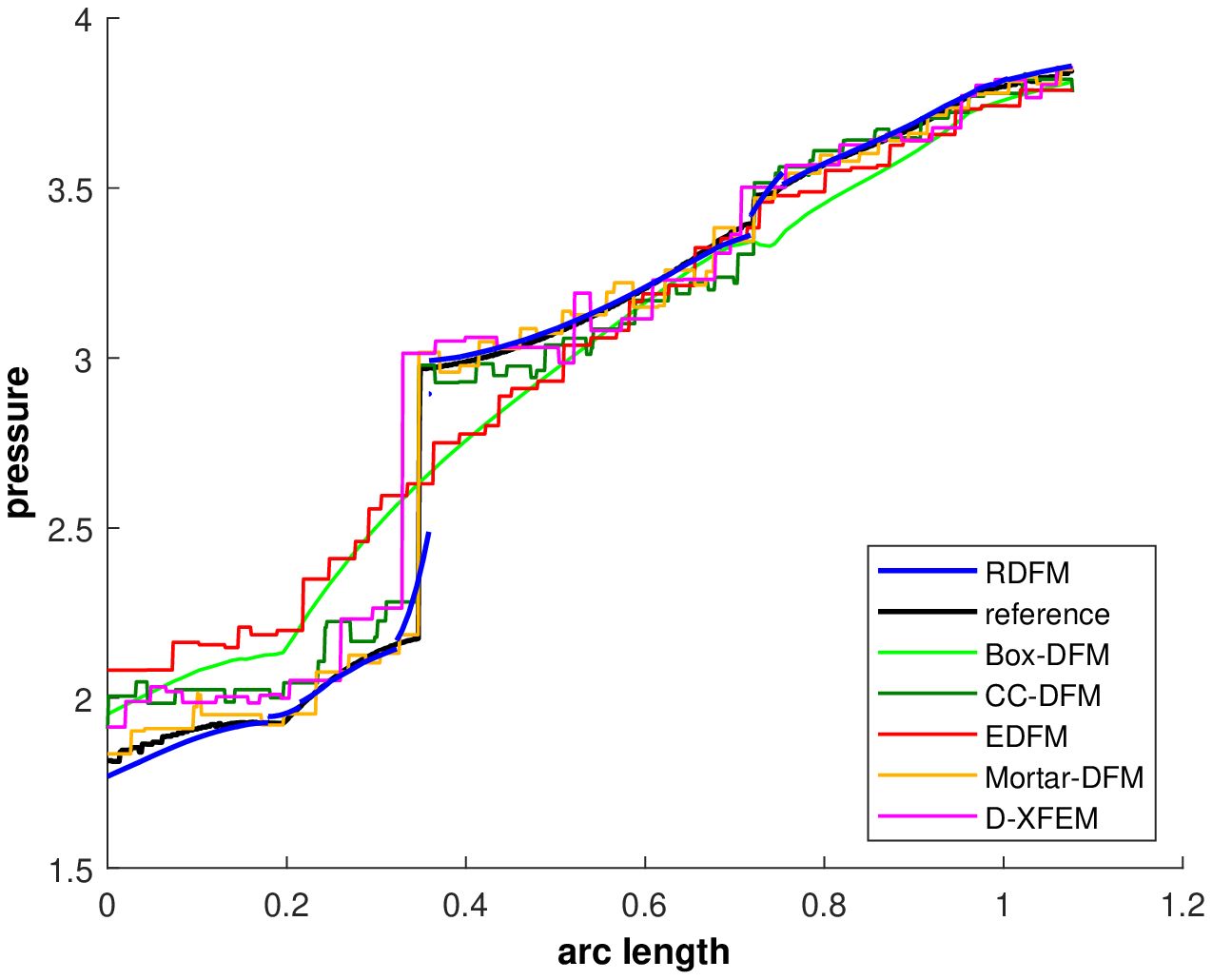}}
\subfigure[Pressure along $(0.0, 0.5)-(1.0, 0.9)$ on $40\times40$ mesh in case (a)]{\includegraphics[width = 3in]{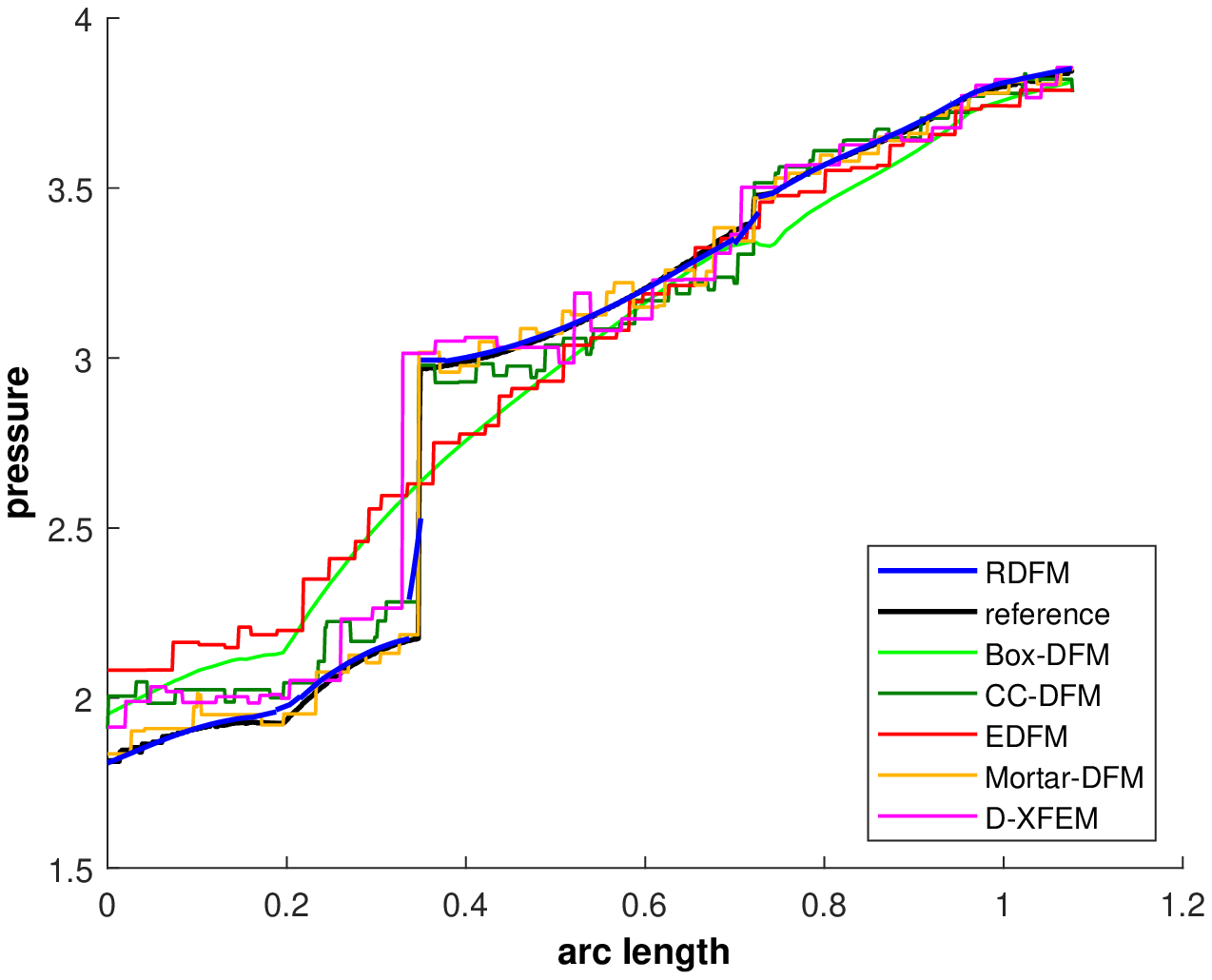}}\\
\subfigure[Pressure along $(0.0, 0.5)-(1.0, 0.9)$ on $30\times30$ mesh in case (b)]{\includegraphics[width = 3in]{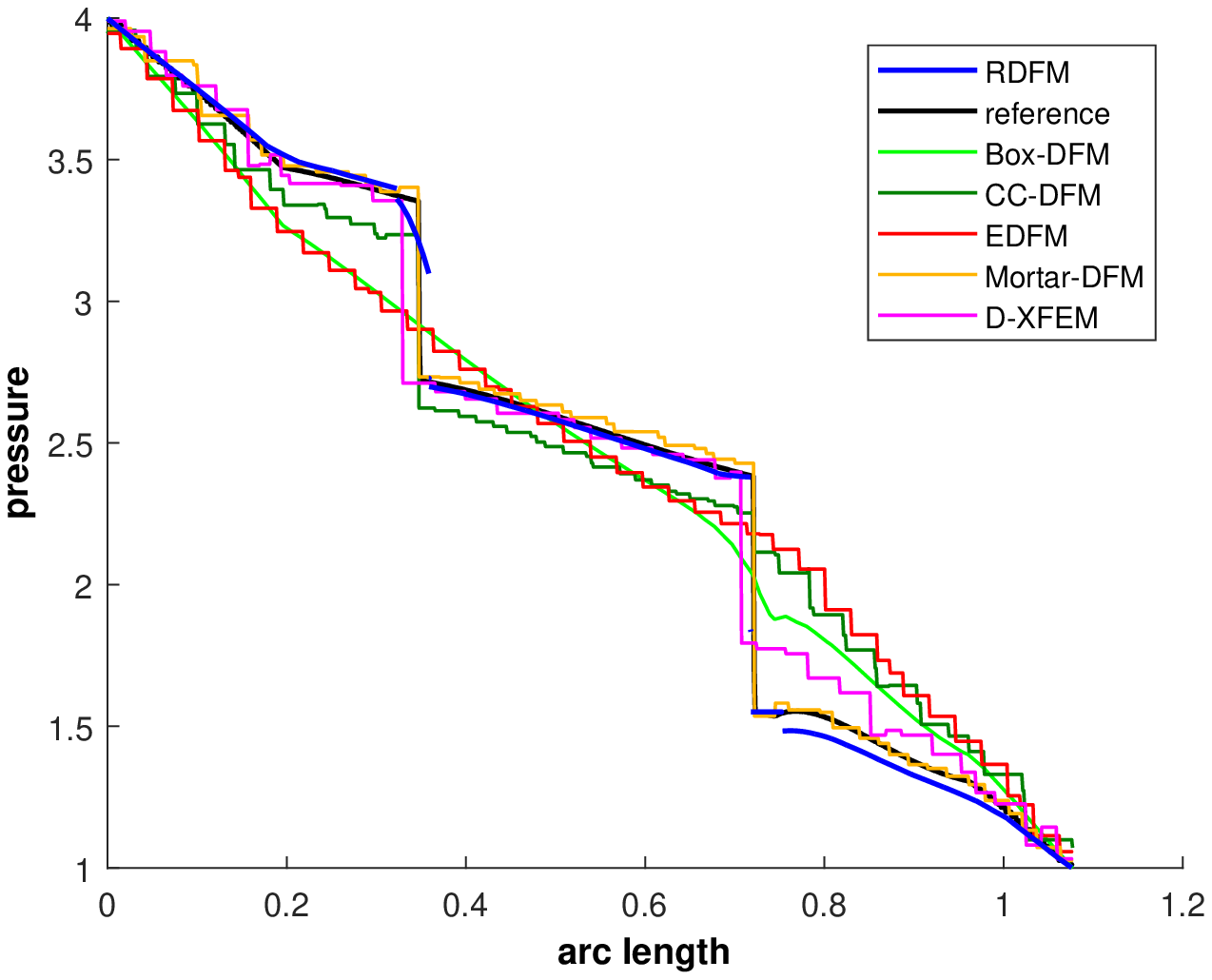}}
\subfigure[Pressure along $(0.0, 0.5)-(1.0, 0.9)$ on $40\times40$ mesh in case (b)]{\includegraphics[width = 3in]{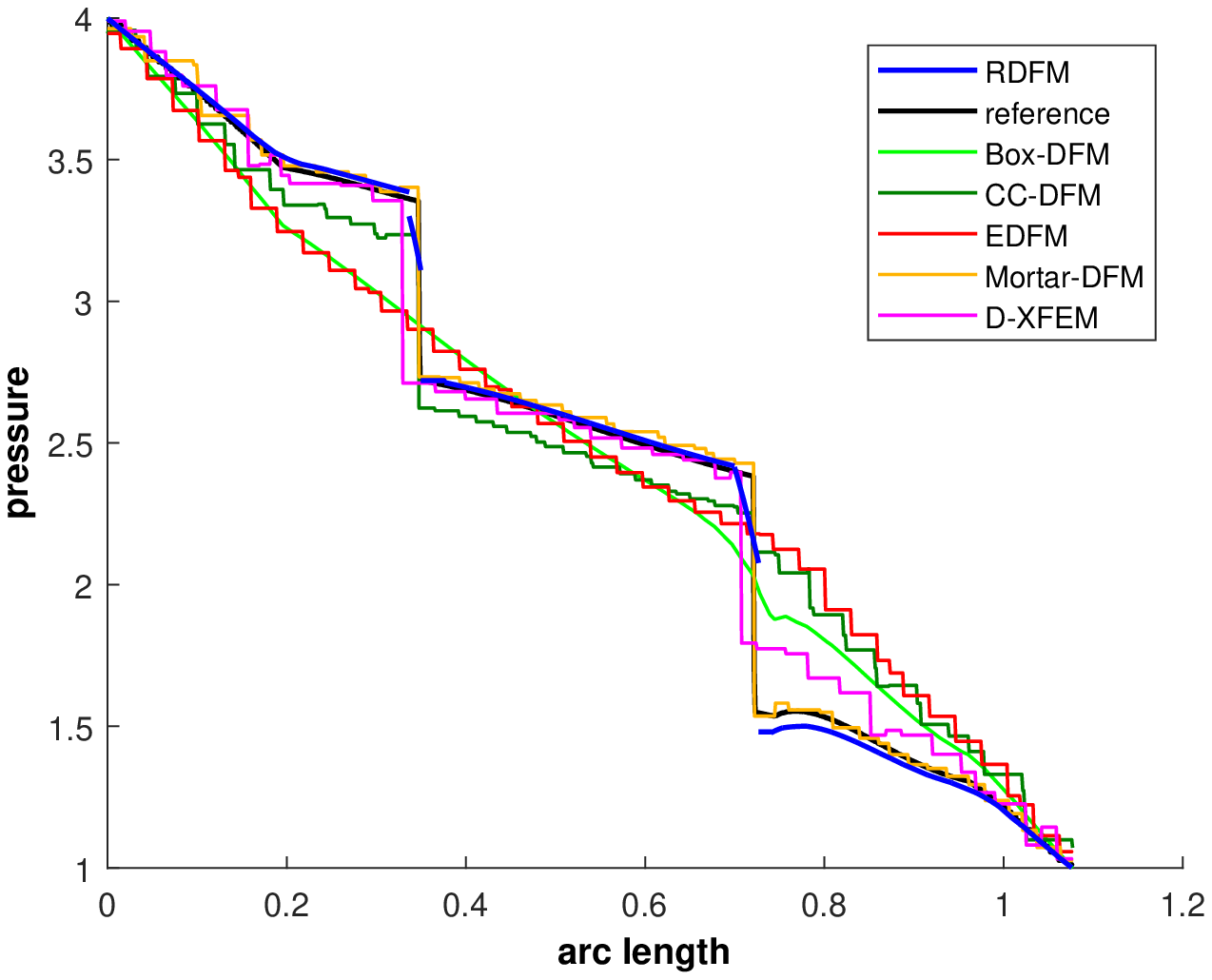}}\\
\caption{Slices of pressure of Example \ref{ex4} on non-conforming meshes}\label{fig:ex4SlicesNonconforming}
\end{figure}

\begin{figure}[!htbp]
\subfigure[Solution of case (a) on conforming mesh containing $1332$ triangles]{\includegraphics[width = 3in]{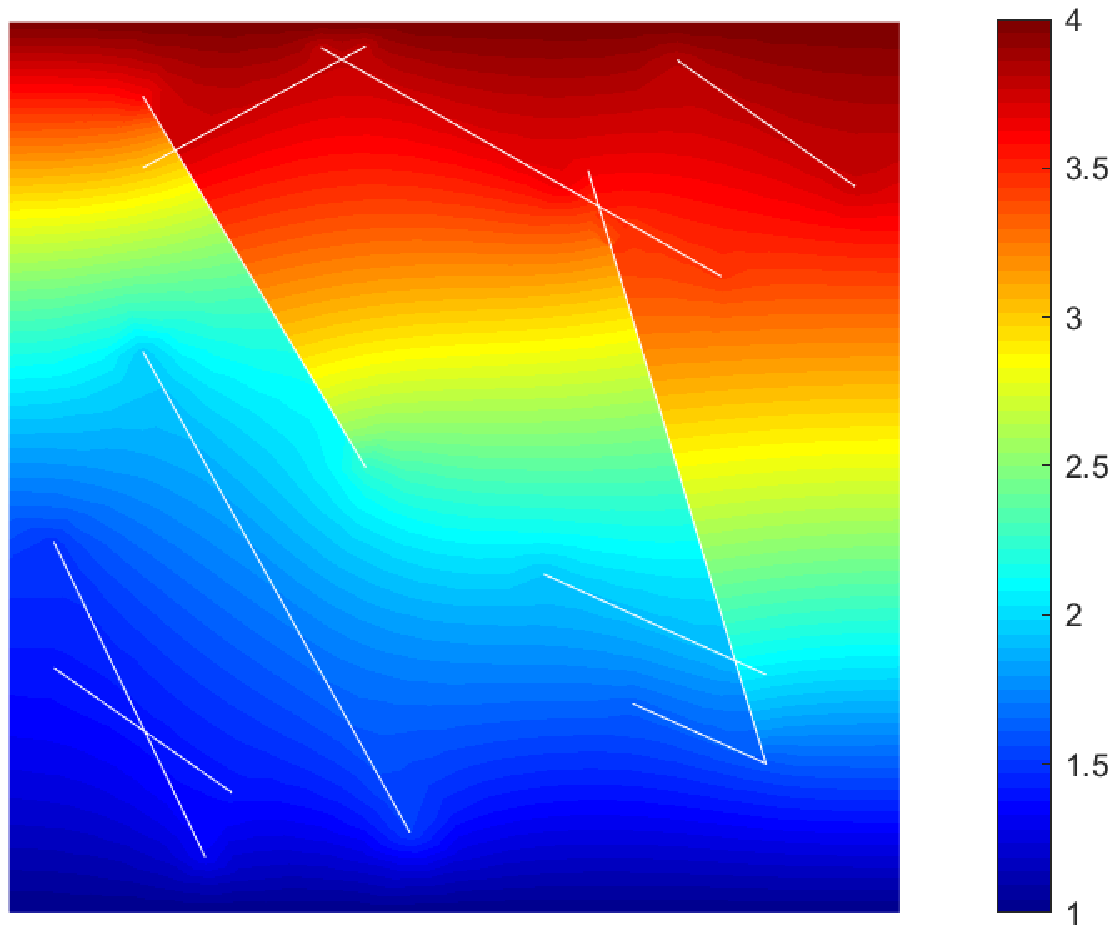}}
\subfigure[Solution of case (a) on conforming mesh containing $2664$ triangles]{\includegraphics[width = 3in]{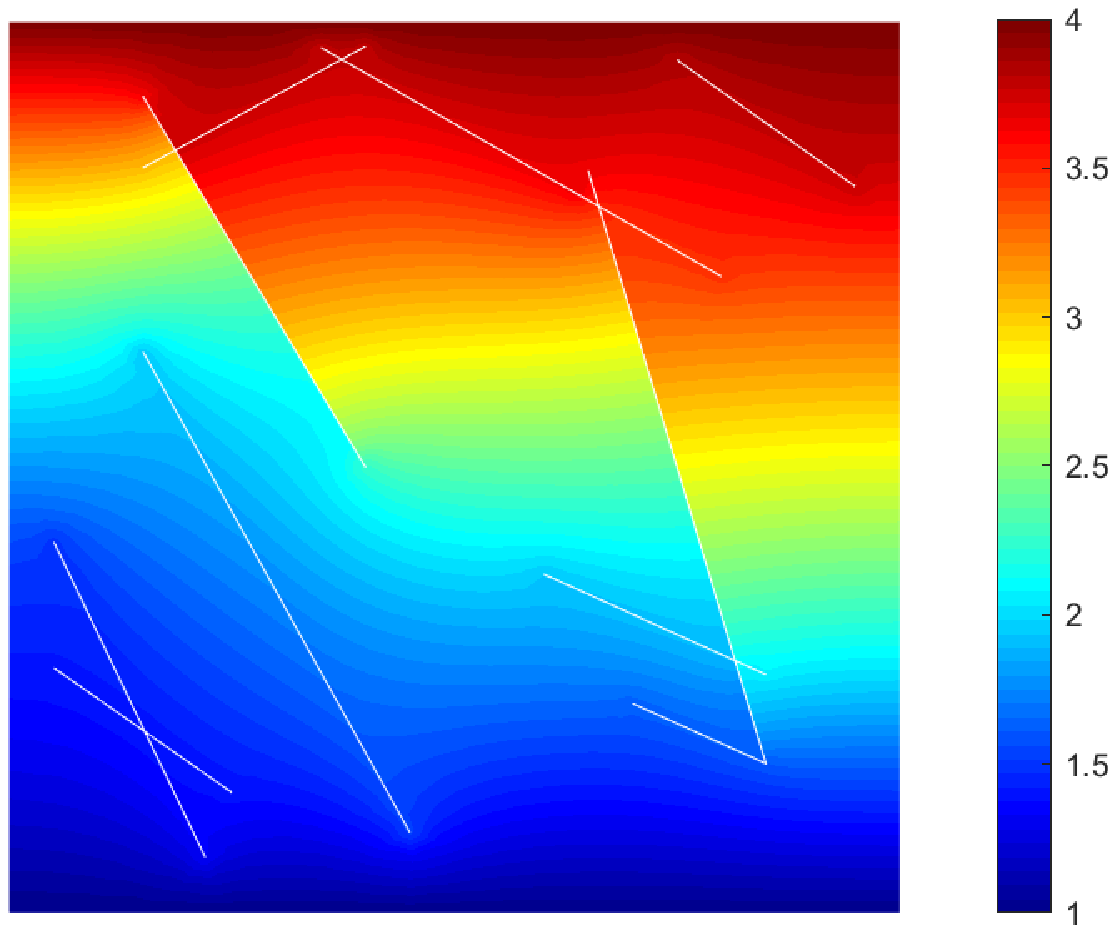}}\\
\subfigure[Solution of case (b) on conforming mesh containing $1332$ triangles]{\includegraphics[width = 3in]{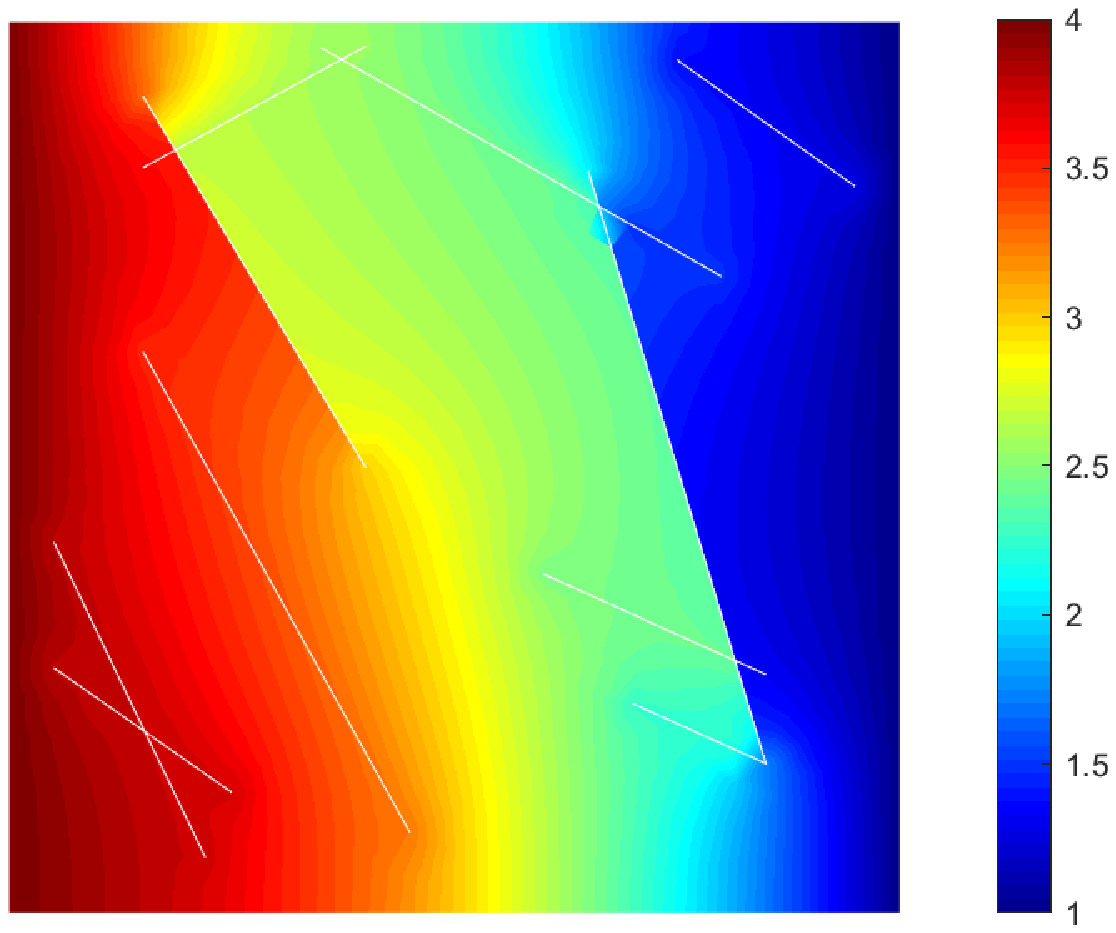}}
\subfigure[Solution of case (b) on conforming mesh containing $2664$ triangles]{\includegraphics[width = 3in]{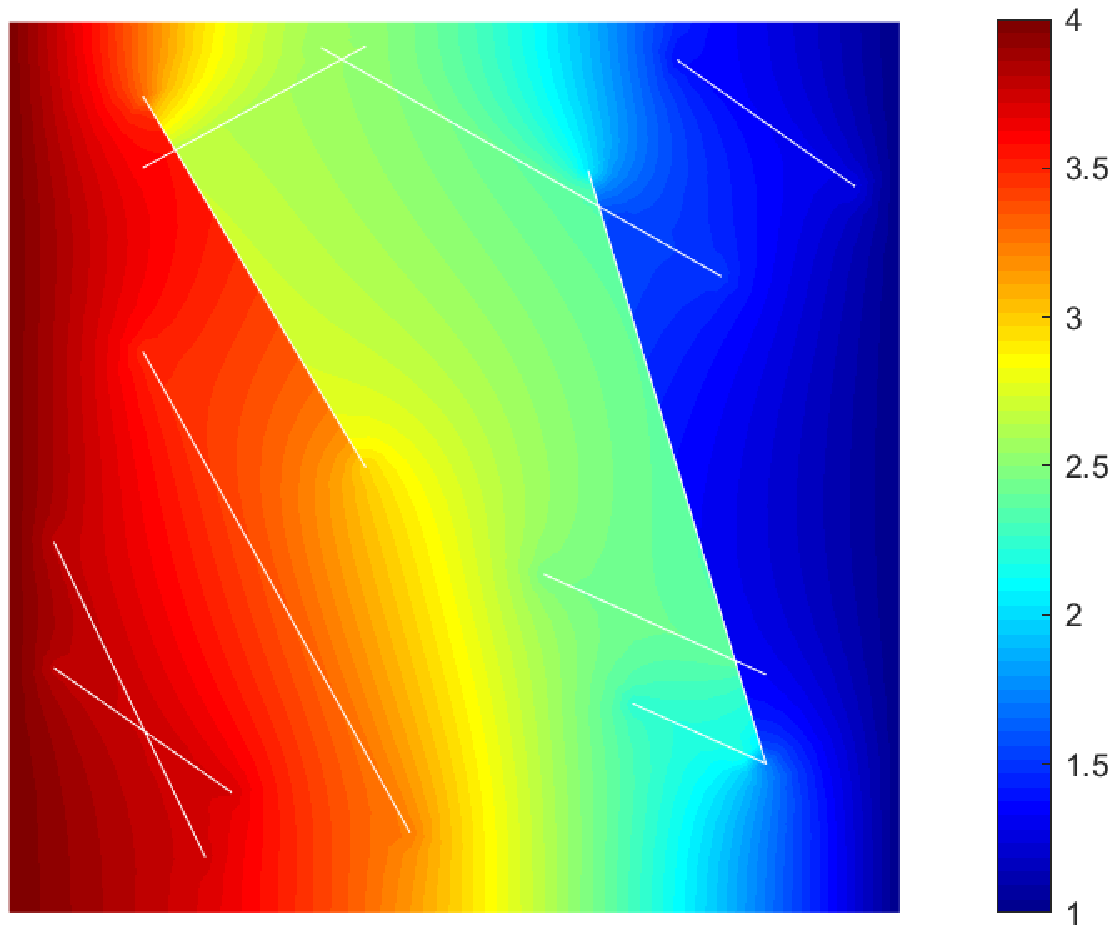}} \\
\caption{Simulation results of Example \ref{ex4} on conforming meshes}\label{fig:ex4contourConforming}
\end{figure}

\begin{figure}[!tb]
\subfigure[Pressure along $(0.0, 0.5)-(1.0, 0.9)$ on conforming mesh containing $1332$ triangles in case (a)]{\includegraphics[width = 3in]{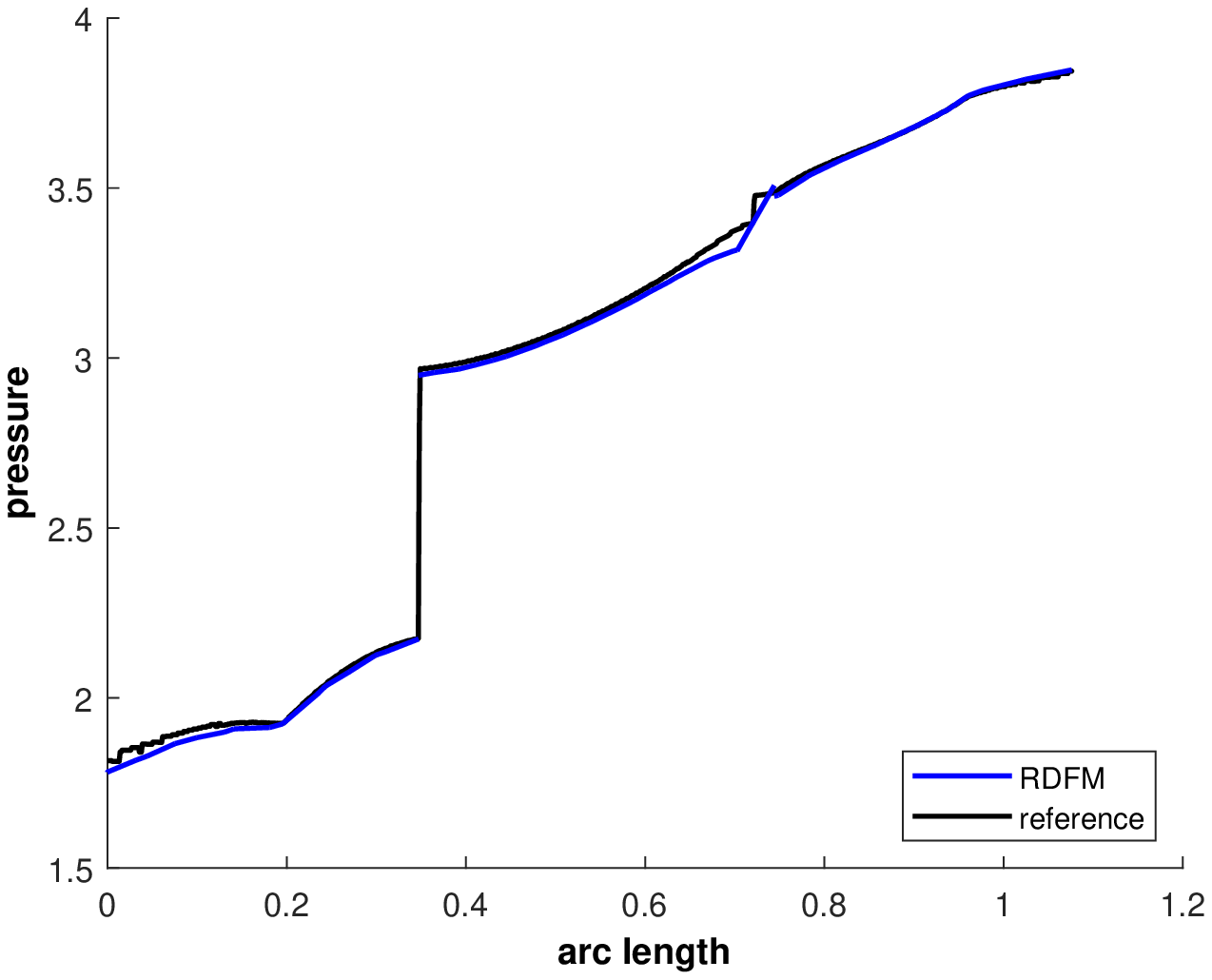}}
\subfigure[Pressure along $(0.0, 0.5)-(1.0, 0.9)$ on conforming mesh containing $2664$ triangles in case (a)]{\includegraphics[width = 3in]{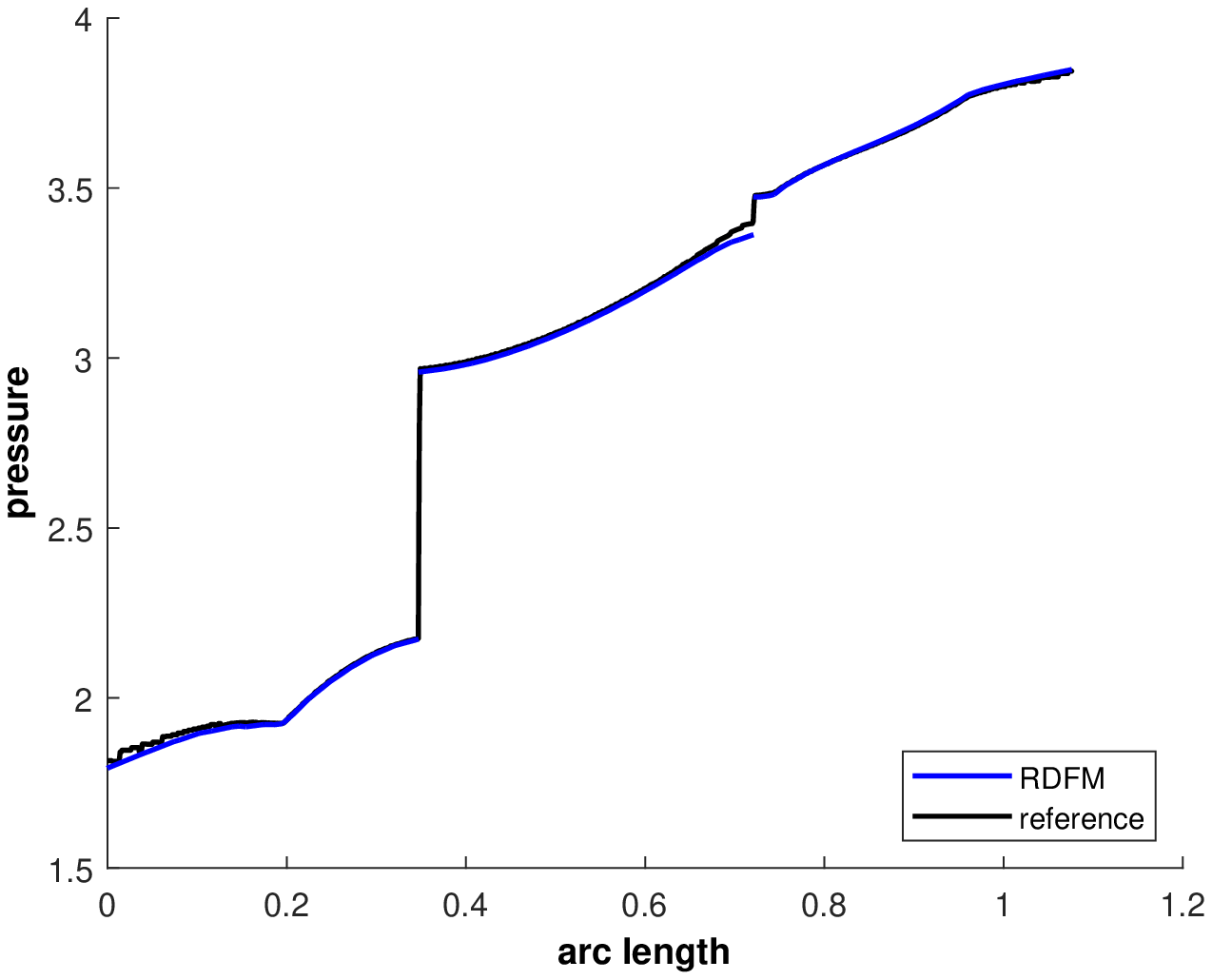}}\\
\subfigure[Pressure along $(0.0, 0.5)-(1.0, 0.9)$ on conforming mesh containing $1332$ triangles in case (b)]{\includegraphics[width = 3in]{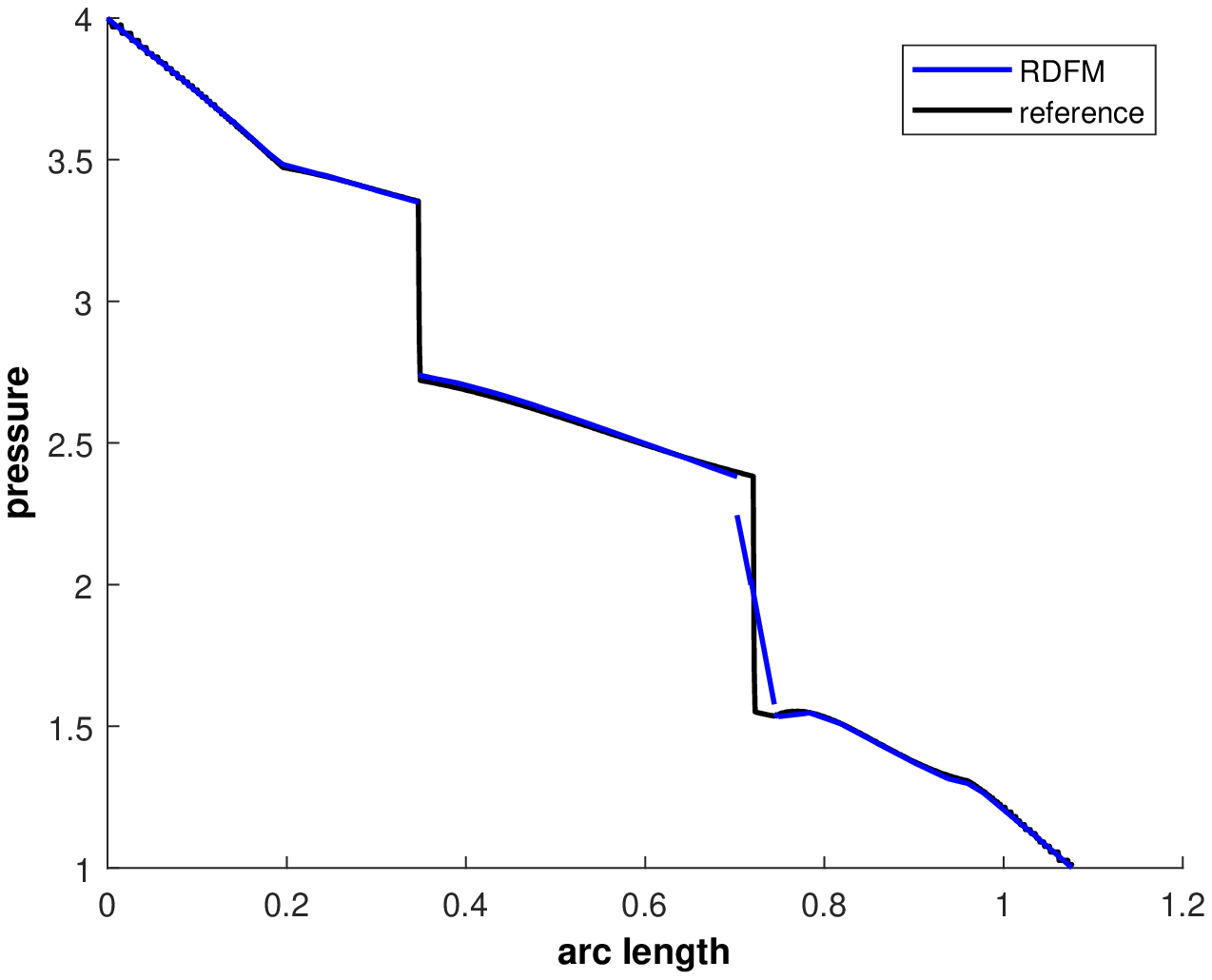}}
\subfigure[Pressure along $(0.0, 0.5)-(1.0, 0.9)$ on conforming mesh containing $2664$ triangles in case (b)]{\includegraphics[width = 3in]{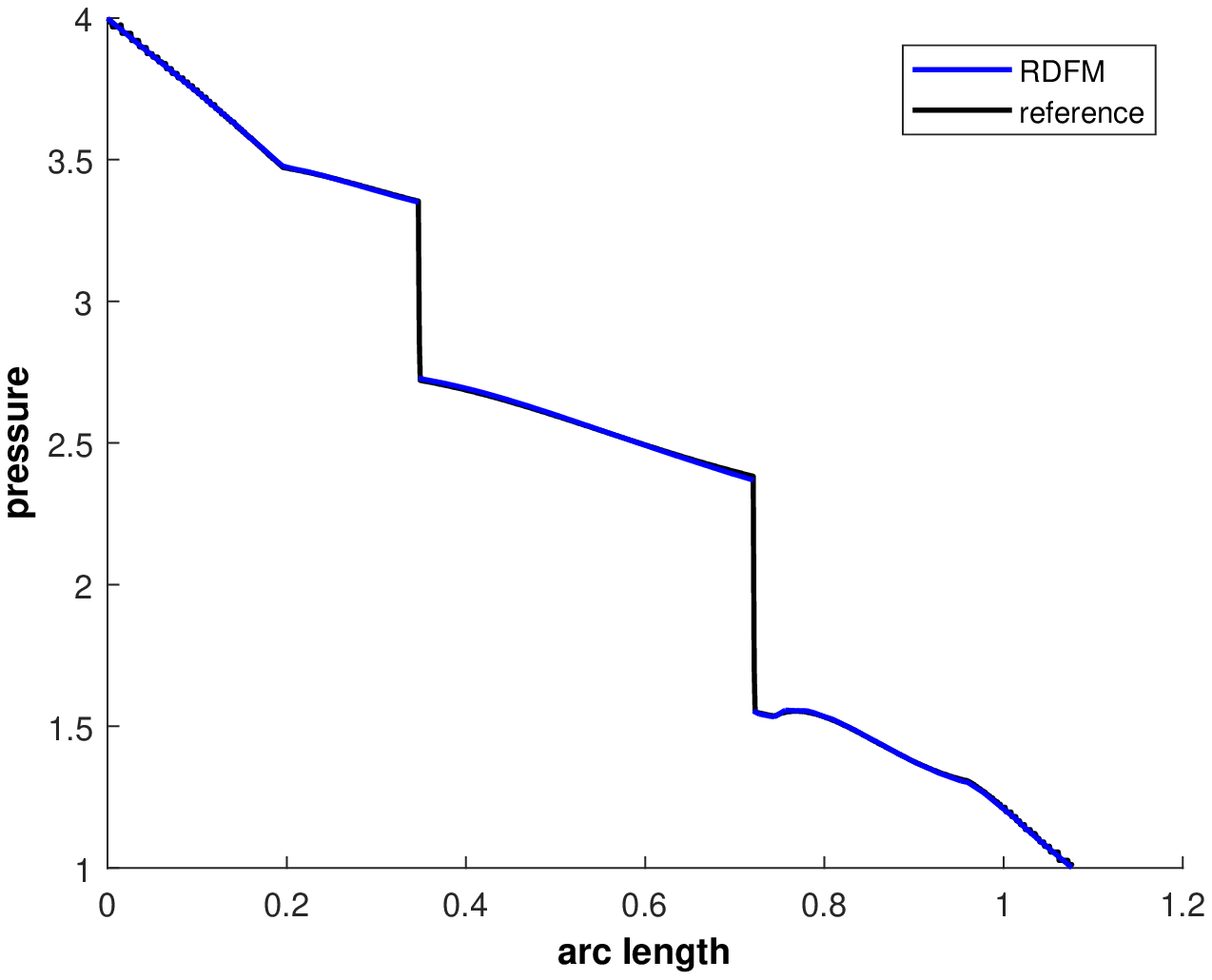}}\\
\caption{Slices of pressure of Example \ref{ex4} on conforming meshes}\label{fig:ex4SlicesConforming}
\end{figure}

\begin{table}[!htbp]
\centering
\caption{\label{tab:ex4a} Evaluation data of different methods of case (a) of Example \ref{ex4}}
\begin{tabular}{|c|c|c|c|c|c|}
  \hline
method & matrix elements & fracture elements & mesh & err$_m$ &err$_f$ \\\hline
Box-DFM&$2838$ triangles & $155$ & conforming       &4.4e-02 & 3.8e-02  \\\hline
CC-DFM&$1407$ triangles & $103$ & conforming & 2.6e-02 & 3.3e-02 \\\hline
EDFM&$1369$ rectangles & $203$ & non-conforming & 3.8e-02 & 4.5e-02 \\\hline
Mortar-DFM&$1452$ triangles & $105$ & conforming  & 1.0e-02 & 1.7e-02 \\ \hline
D-XFEM&$1922$ triangles & $199$ & non-conforming  & 1.9e-02 & 2.9e-02 \\ \hline
RDFM& $900$ rectangles & $167$ &non-conforming & 2.4e-02 & 3.5e-02 \\ \hline
RDFM& $1600$ rectangles & $219$ &non-conforming & 1.8e-02 & 3.1e-02 \\ \hline
RDFM& $1332$ triangles & $222$ &conforming & 1.0e-02 & 4.8e-02 \\ \hline
RDFM& $2664$ triangles & $341$ &conforming & 7.7e-03 & 4.8e-02 \\ \hline
\end{tabular}
\end{table}

\begin{table}[!htbp]
\centering
\caption{\label{tab:ex4b} Evaluation data of different methods of case (b) of Example \ref{ex4}}
\begin{tabular}{|c|c|c|c|c|c|}
  \hline
method & matrix elements & fracture elements & mesh & err$_m$ &err$_f$ \\\hline
Box-DFM&$2838$ triangles & $155$ & conforming       &7.5e-02 & 7.0e-02  \\\hline
CC-DFM&$1407$ triangles & $103$ & conforming & 5.2e-02 & 7.3e-02 \\\hline
CC-DFM*&$1407$ triangles & $103$ & conforming & 1.1e-02 & 2.7e-02 \\\hline
EDFM&$1369$ rectangles & $203$ & non-conforming & 5.8e-02 & 8.9e-02 \\\hline
Mortar-DFM&$1452$ triangles & $105$ & conforming  & 1.3e-02 & 2.7e-02 \\ \hline
D-XFEM&$1922$ triangles & $199$ & non-conforming  & 2.2e-02 & 3.6e-02 \\ \hline
RDFM& $900$ rectangles & $167$ &non-conforming & 2.7e-02 & 5.4e-02 \\ \hline
RDFM& $1600$ rectangles & $219$ &non-conforming & 2.3e-02 & 5.1e-02 \\ \hline
RDFM& $1332$ triangles & $222$ &conforming & 1.2e-02 & 7.3e-02 \\ \hline
RDFM& $2664$ triangles & $341$ &conforming & 7.8e-03 & 7.4e-02 \\ \hline
\end{tabular}
\end{table}

The reference solutions are computed by MFD on a very fine mesh containing $1192504$ matrix elements and $7876$ fracture elements.
The models that participate in the comparison are the same as before excepts the P-XFEM, which could not join due to its severe restriction on the geometry of fracture networks.
Moreover, another version of cell centered control volume discrete fracture model (CC-DFM*) without the elimination of intermediate fracture intersection cells is included in the comparison for case (b).

We first compute the solutions on $30\times30$ and $40\times40$ non-conforming uniform rectangular meshes for both cases.
Then, we test the model on the conforming triangular meshes with grids data provided by \cite{Benchmark}, which contains $1332$ and $2664$ triangles, respectively.

The contour plots of the simulation results on conforming and non-conforming meshes are shown in Figure \ref{fig:ex4contour} and \ref{fig:ex4contourConforming}, respectively.
One can find the contours of the reference solutions in \cite{Benchmark} for comparison.
We slice the solution along the line $(0.0, 0.5)-(1.0, 0.9)$, and draw the comparison with those of other methods in Figure \ref{fig:ex4SlicesNonconforming} and \ref{fig:ex4SlicesConforming}.

We also compare the relative errors on matrix and fractures, and other important aspects of different methods in Table \ref{tab:ex4a} and Table \ref{tab:ex4b} for case (a) and case (b), respectively.

From the tables, we can observe a relatively small matrix error for our methods, but the error on fracture-barrier network is not as small as that on matrix due to the pressure discontinuity on barriers.

\section{Applications in contaminant transportation in fractured porous media}

In this section, we study the applications of the RDFM in contaminant transportation in fractured porous media.
The flow field is governed by the equations \eqref{eq:HybridDimDarcy},\eqref{eq:ContEq} and \eqref{BVP}.
The transportation of the contaminant in flow field satisfies the convection-diffusion equation \eqref{convdiffu} with the corresponding boundary conditions \eqref{inflow}, \eqref{outflow} and initial condition \eqref{initial}:

\begin{align}
(\phi c)_t + \nabla\cdot(\bm{u}c) - \nabla\cdot\left({\bm{D}\nabla c}\right)&= f \Tilde{c},\quad (x,t)\in \Omega\times(0,T], \label{convdiffu}\\
\left(\bm{u}c-\bm{D}\nabla c\right)\cdot\bm{n}&=c_{\text{in}}\bm{u}\cdot\bm{n}, \quad (x,t)\in \Gamma_{\text{in}}\times(0,T],\label{inflow}\\
(-\bm{D}\nabla c)\cdot\bm{n}&=0, \quad (x,t)\in \Gamma_{\text{out}}\times(0,T],\label{outflow}\\
c(x,0)&=c_0(x), \quad x\in\Omega, \label{initial}
\end{align}
where $c$ is the concentration of the contaminant, $\phi$ is the porosity of porous media, $\bm{u}$ is the Darcy's velocity of the flow field, $\bm{D}$ is the diffusion coefficient of the contaminant, $f$ is the source term, $\tilde{c}=c$ if $f<0$ and $\tilde{c}=c_{\text{inject}}$ if $f>0$, $\Gamma_{\text{in}}$ is the inflow boundary and $\Gamma_{\text{out}}$ is the outflow boundary, $c_{\text{in}}$ is the concentration of the contaminant on the inflow boundary, $c_{0}$ is the initial concentration of the contaminant.

Problem \eqref{convdiffu}-\eqref{initial} is discretized by the interior penalty discontinuous Galerkin (IPDG) method in space. The semidiscrete formulation is given as follows:

Find $c\in V_h$, s.t. $\forall v\in V_h$,
\begin{equation} \label{IPDG}
\begin{split}
\left(\phi c_t, v\right)=&\left(\bm{u}c-\bm{D}\nabla c,\nabla v\right)-\int_{\Gamma_0}\hat{\bm{uc}}[v]~ds - \int_{\Gamma_{\text{in}}} \bm{u}c_{\text{in}}[v]~ds - \int_{\Gamma_{\text{out}}} \bm{u}c [v]~ds\\
&+\int_{\Gamma_0}\left(\{\bm{D}\nabla c\}\cdot[v] + \{\bm{D}\nabla v\}\cdot [c]-\tilde{\alpha}[c]\cdot[v]\right)~ds + (f\tilde{c},v),
\end{split}
\end{equation}
where $\tilde{\alpha}=O(\frac{1}{|e|})$ is the parameter of the diffusion penalty term and $\hat{\bm{uc}}|_e=\{\bm{u}c\}+\tilde{\beta}[c]$ is the Lax-Friedrichs flux on the cell interfaces $e$.

There are many choices for the time discretization. For simplicity, we use the forward Euler scheme. In the simulation, we first solve for the flow field $\bm{u}$ from \eqref{scheme:part1}, \eqref{scheme:part2} and \eqref{scheme:part3}, then simulate the time evolution of the concentration of  contaminate in the flow field by computing \eqref{IPDG}. An example is given as below.

\begin{figure}[!htbp]
\subfigure[$0.2$ PVI]{\includegraphics[width = 3in]{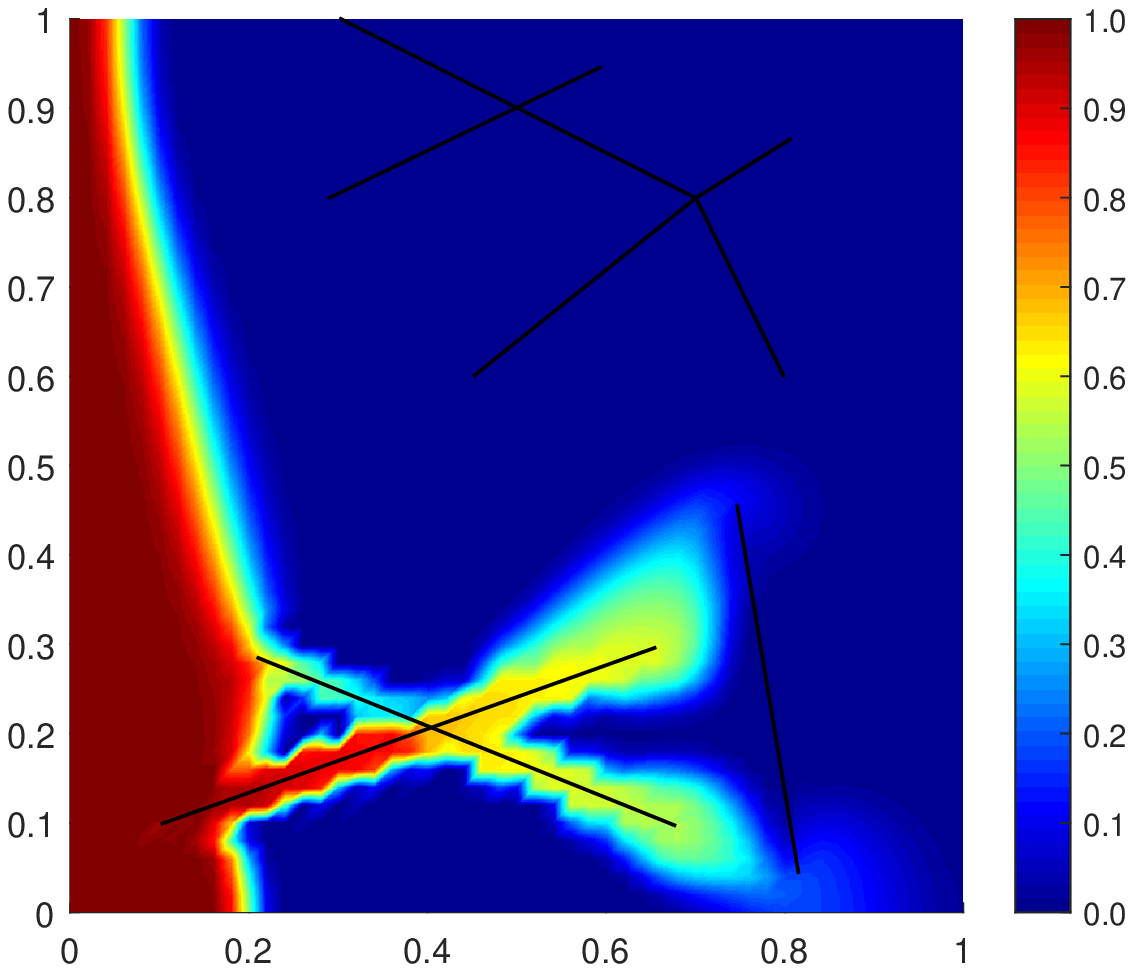}}
\subfigure[$0.4$ PVI]{\includegraphics[width = 3in]{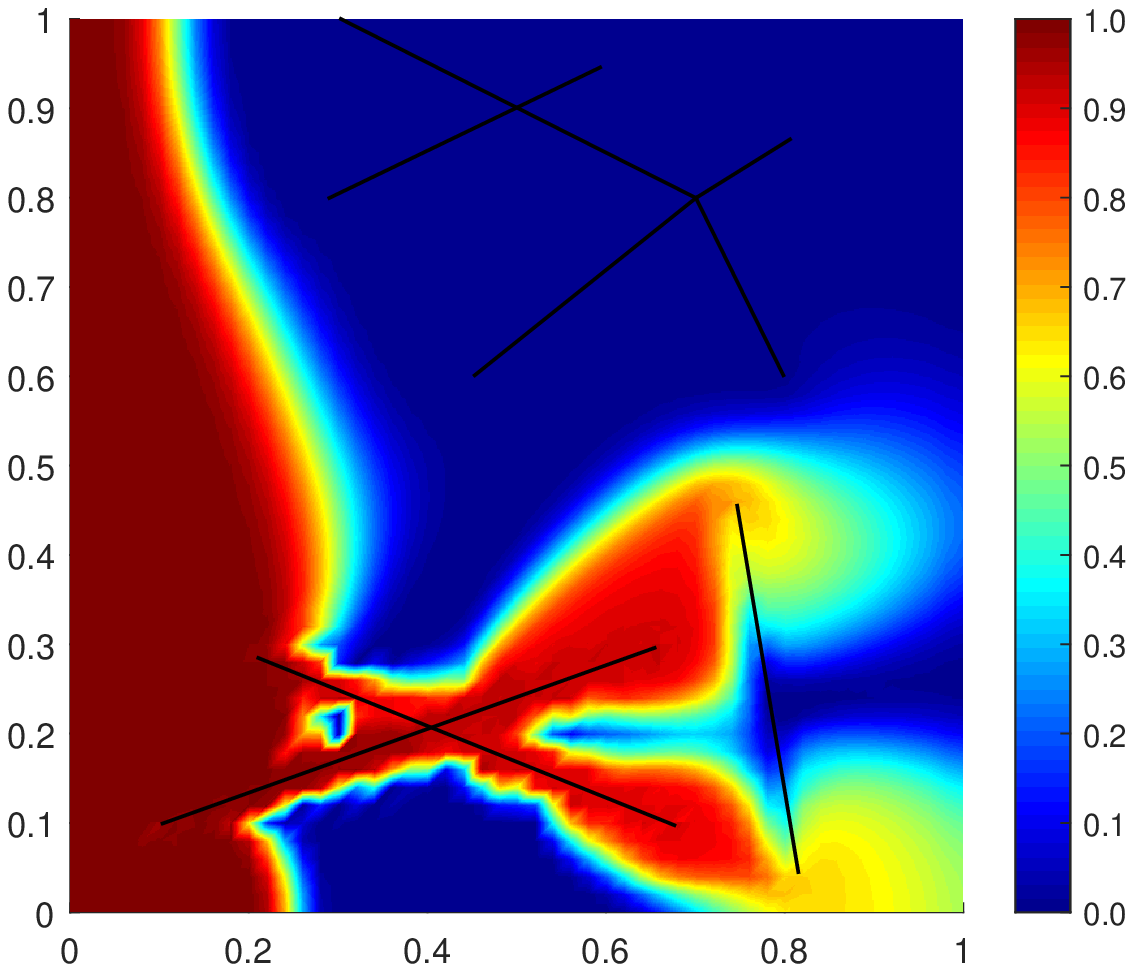}} \\
\subfigure[$0.6$ PVI]{\includegraphics[width = 3in]{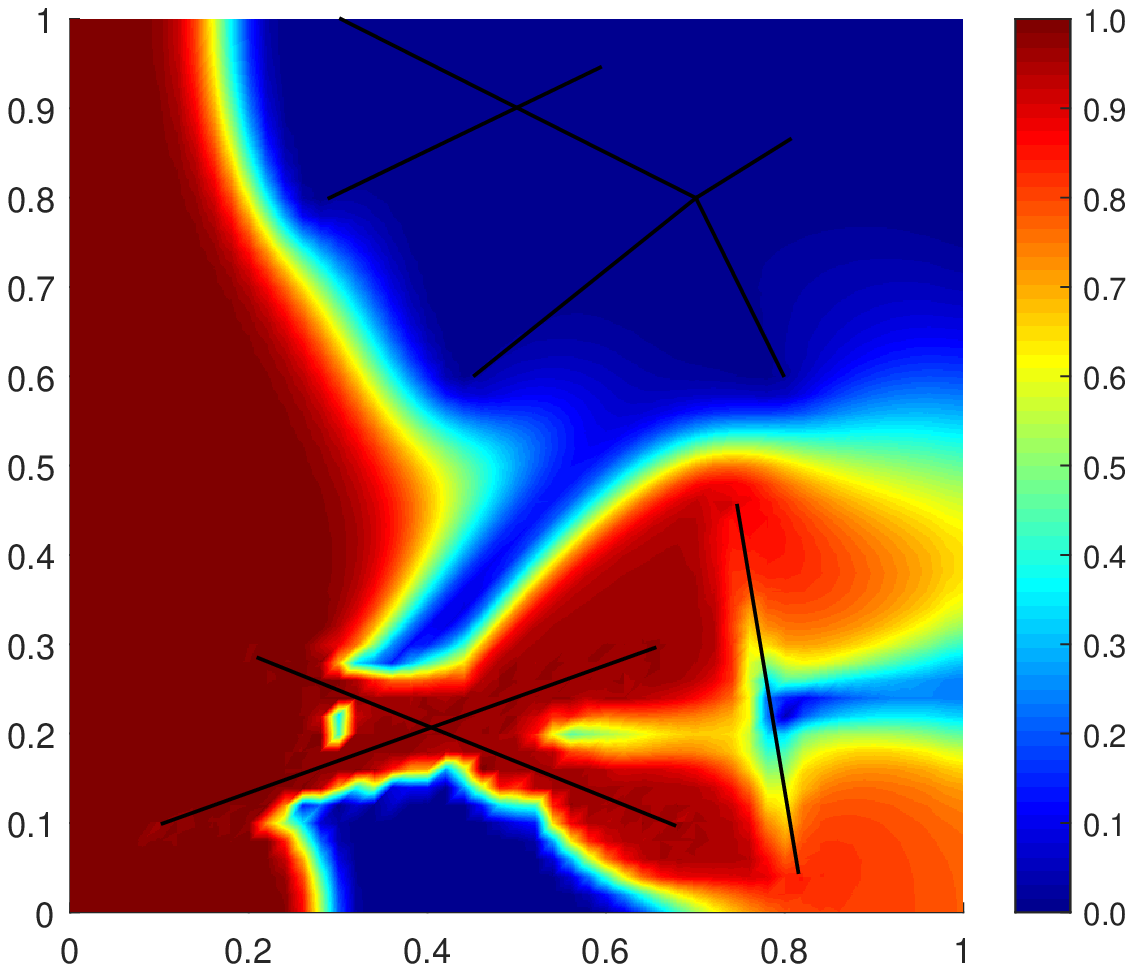}}
\subfigure[$0.8$ PVI]{\includegraphics[width = 3in]{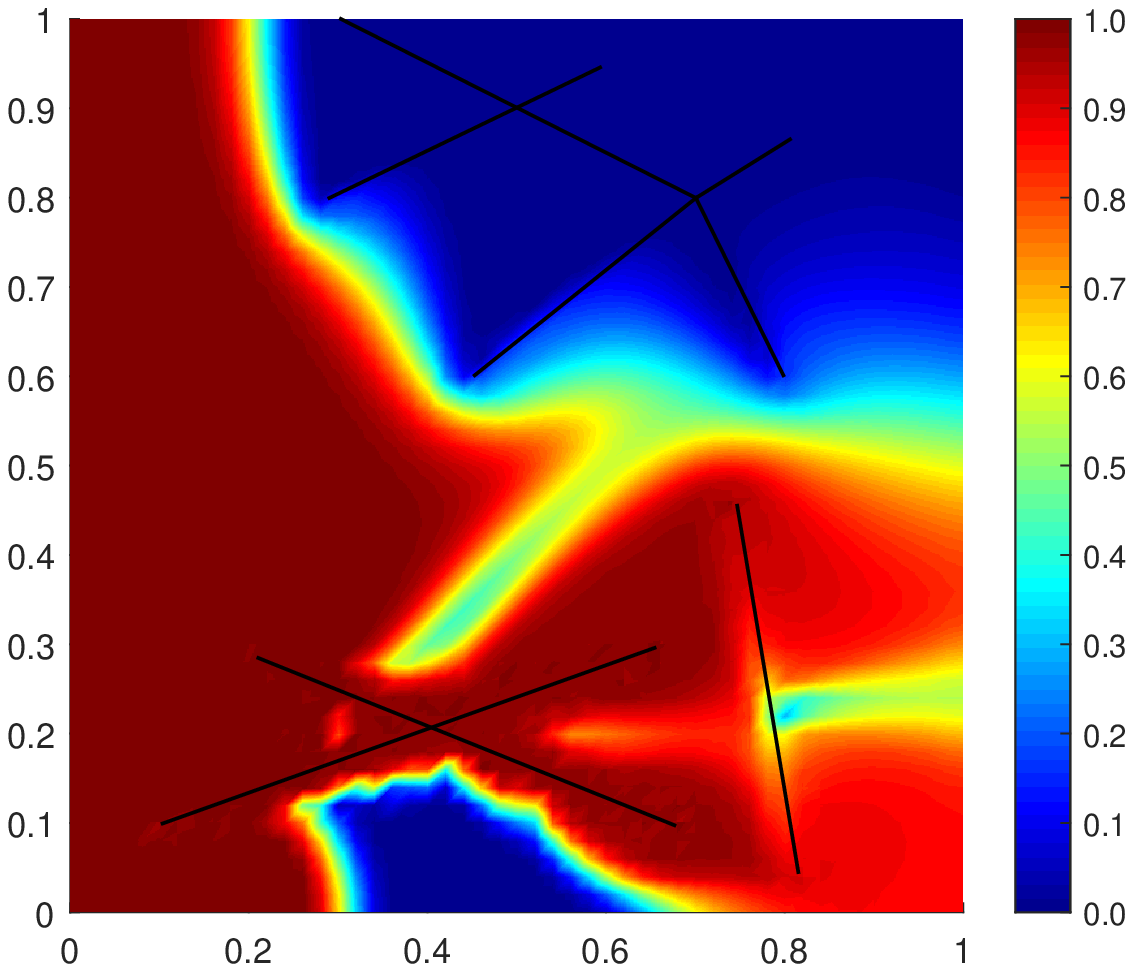}} \\
\subfigure[$1.0$ PVI]{\includegraphics[width = 3in]{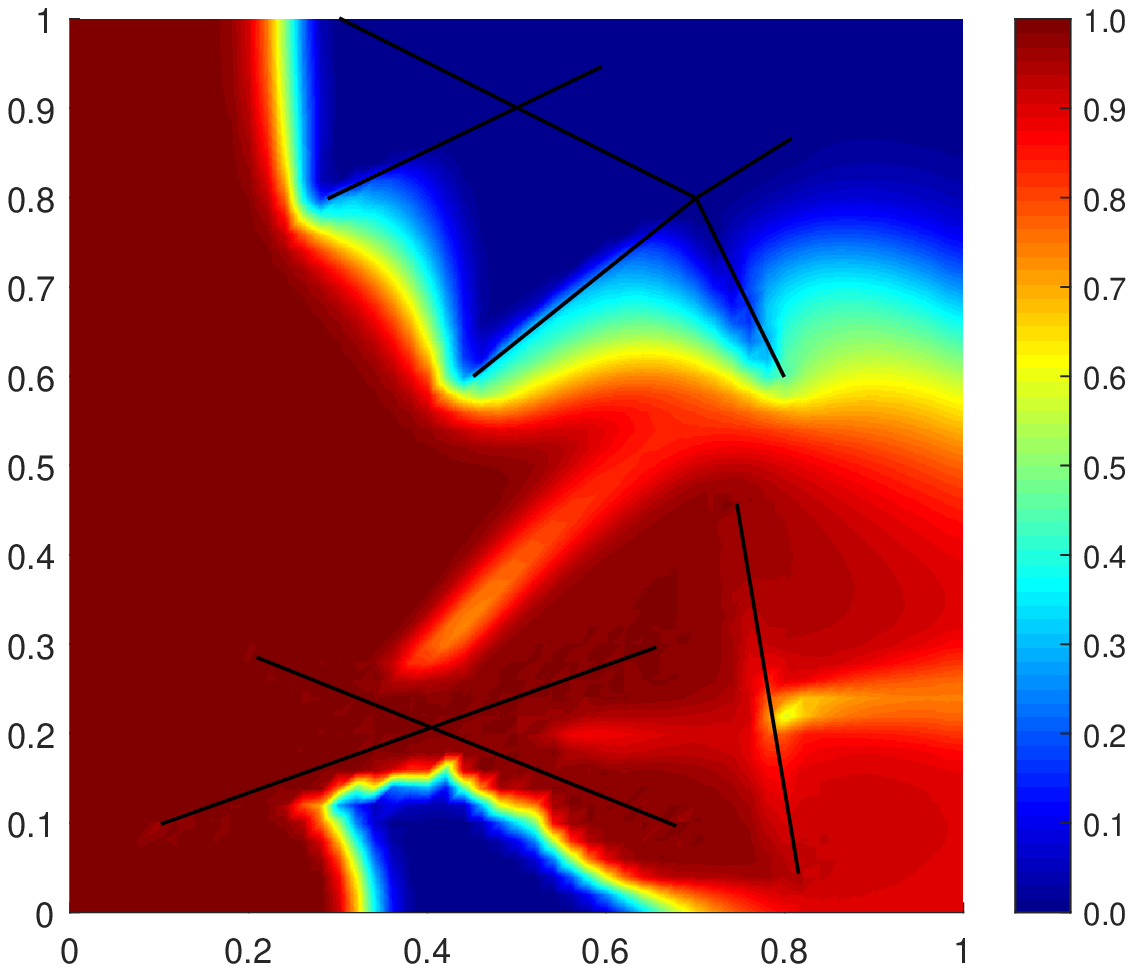}}
\subfigure[$1.2$ PVI]{\includegraphics[width = 3in]{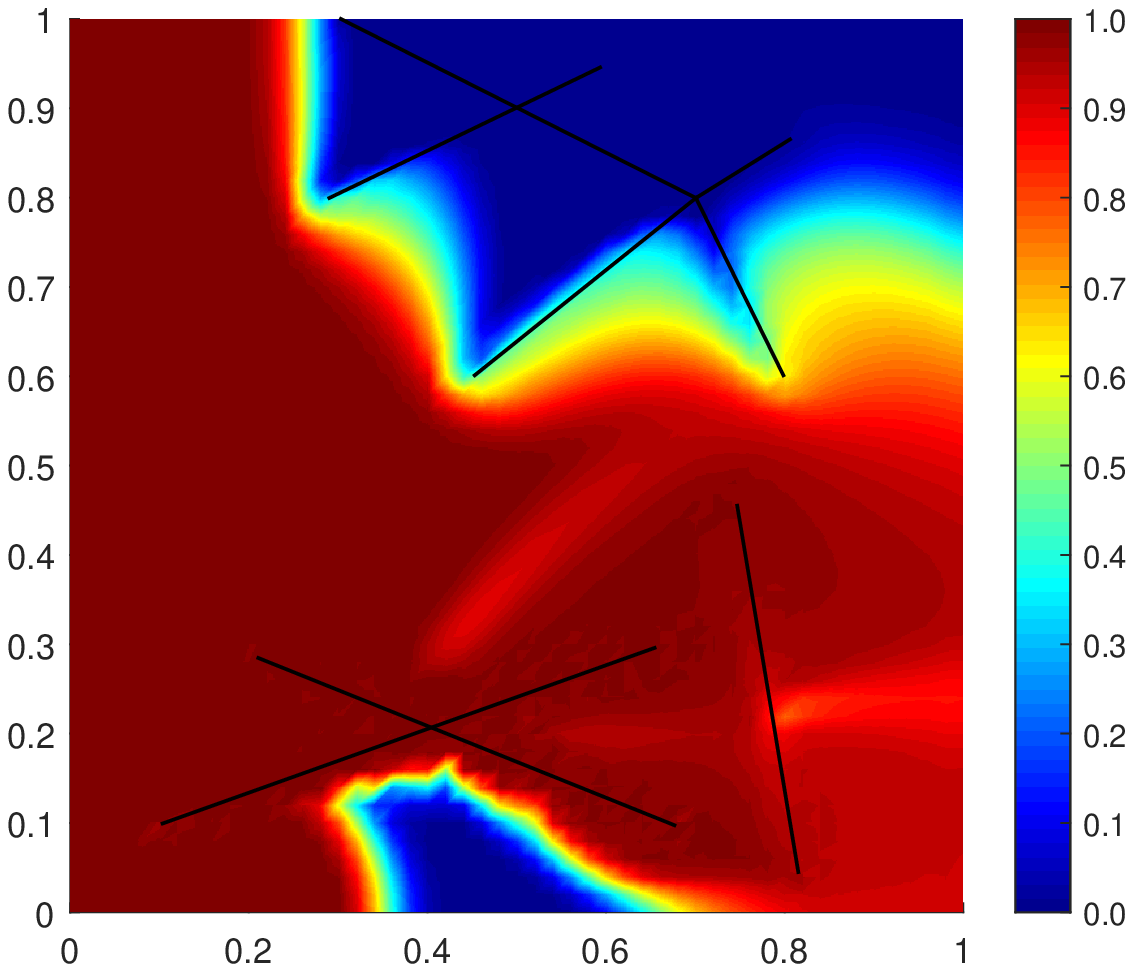}} \\
\caption{Time evolution of concentration of contaminant in Example \ref{ex6_1}}\label{fig:ex6_1concentration}
\end{figure}

\begin{ex}\label{ex6_1}
In this example, we simulate the contaminant transportation in the porous media containing fracture and barrier networks. The settings of the problem referred an example in \cite{Tene2017} but is not exactly the same. The computational domain is the unit square $\Omega=[0,1]\times[0,1]$. The permeability of the porous matrix is $k_m=1$. The porosity of the porous media is $\phi=0.2$. The porous media contains $2$ fractures and $6$ barriers, with the uniform thickness $\epsilon=1\times10^{-4}$, and permeability $k_f=1\times 10^{8}$ and $k_\epsilon=1\times10^{-8}$, respectively. The exact coordinates of the fractures and barriers are attached in the appendix \ref{appd:ex6_1}, see Figure \ref{fig:ex6_1pressure}(a) for their distribution.
For simplicity, the flow is driven by boundary conditions rather than the source term. The left and right boundaries are Dirichlet with pressure $p_D=1$ and $p_D=0$, respectively. The top and bottom boundaries are impermeable. Moreover, the left boundary is an inflow boundary with $c_{\text{in}}=1$ and the right boundary is an outflow boundary. The diffusion coefficient is set to be $\bm{D}=0.005|\bm{u}|$.
\end{ex}
The computation is conducted on a $50\times50$ non-conforming uniform rectangular mesh.
The contour plot of the pressure is shown in Figure \ref{fig:ex6_1pressure}(b).
The time evolution of the concentration of contaminant is exhibited in Figure \ref{fig:ex6_1concentration}, with different pore volume injected (PVI).

\begin{figure}[!htbp]
\subfigure[Domain and boundary conditions]{\includegraphics[width = 3in]{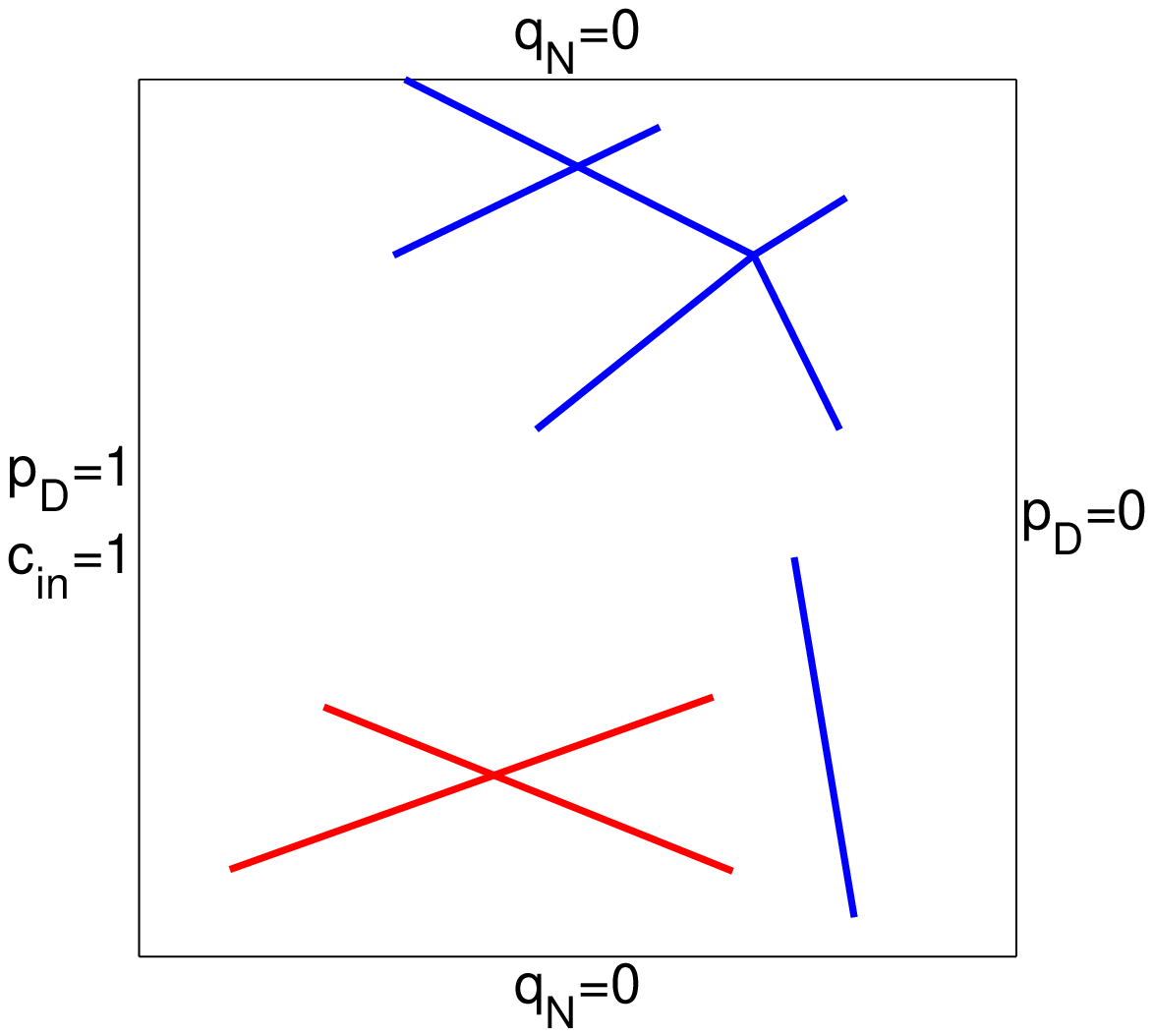}}
\subfigure[Pressure distribution on $50\times50$ mesh]{\includegraphics[width = 3in]{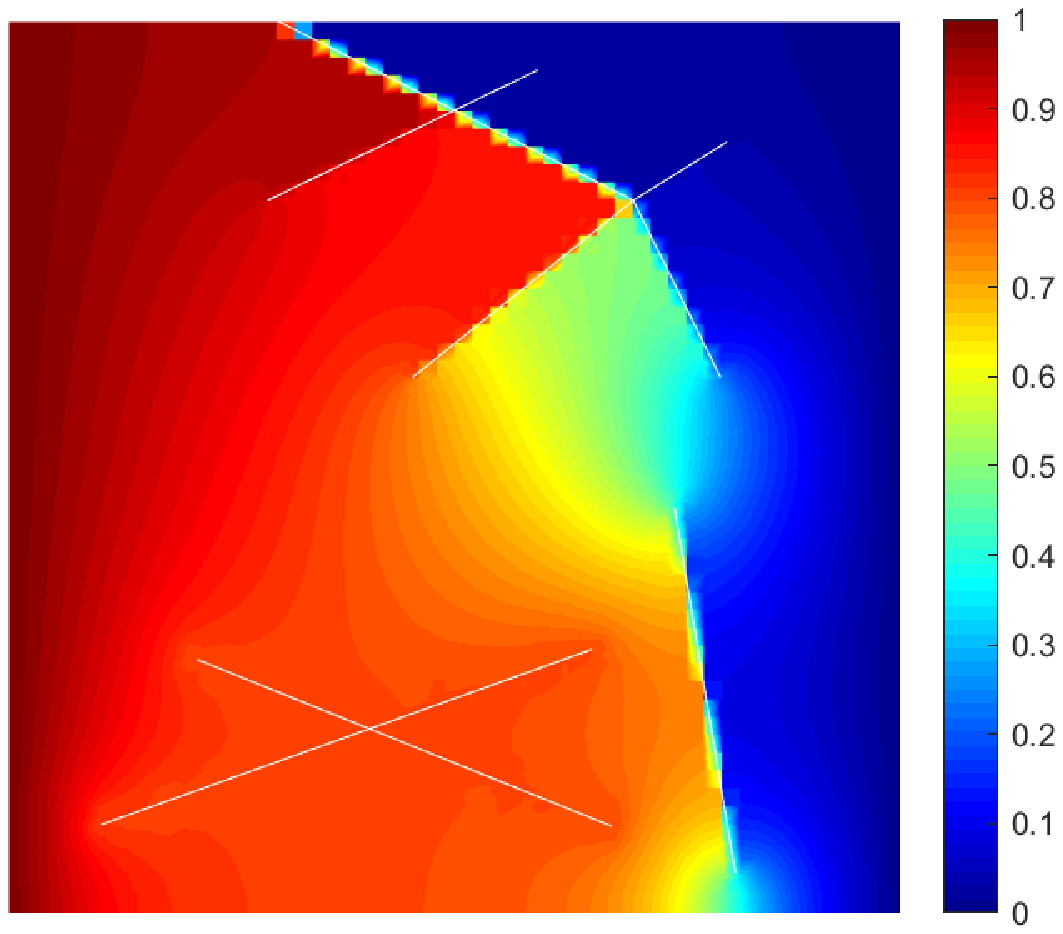}} \\
\caption{Domain setting and pressure distribution of Example \ref{ex6_1}}\label{fig:ex6_1pressure}
\end{figure}

\section{Concluding remarks}
In this paper, we have proposed the hybrid-dimensional Darcy's law and established a novel discrete fracture model by applying the local discontinuous Galerkin method to it.
Several numerical experiments showed its validity on non-conforming meshes and high accuracy on conforming meshes for porous media containing fracture and barrier networks.

Unlike the pEDFM, which are restricted on rectangular cells to account for barriers, an advantage of our model is that it does not have any requirements on the shape of elements because it describes barriers analogy to fractures. Moreover, our method is locally mass conservative thanks to the property of discontinuous Galerkin methods.

The major limitation of the method is that, if the thickness of fractures or barriers is large, our method has non-negligible model errors because we adopt Dirac-$\delta$ functions to characterize fractures and barriers.
Also, our model cannot describe the totally impermeable barriers, i.e. the case $k_\epsilon=0$.
Another shortcoming of the method is that the scheme is computationally costly if the pressure is the only variable of interest, since the unknown in the LDG scheme includes Darcy's velocity and the gradient of pressure.

There are several possible improvements and future works to be done.
First, it seems there are no essential difficulties in extending our method to curved fractures and barriers, three dimensional problems, and multi-phase flow in fractured porous media. Therefore it is worth exploring the performance of the model in these scenarios.
Another interesting scenario is the fracture being highly conductive in its tangential direction while low permeable in its normal direction. This kind of 'fracture' can be described by the hybrid-dimensional Darcy's law in principle but does not fit the current LDG scheme. We hope to establish an appropriate numerical scheme in future works to simulate this case.
Also, we expect to find more efficient numerical schemes for the PDE model \eqref{eq:HybridDimDarcy}, \eqref{eq:ContEq} and \eqref{BVP} besides the LDG methods.



\newpage
\begin{appendices}
\appendix

\section{Exact coordinates of fracture networks in Example \ref{ex4}} \label{appd:ex4}

\begin{Verbatim}
FRACTURES:
NUMBER, START_X, START_Y, END_X, END_Y 
1, 0.0500 0.4160 0.2200 0.0624;
2. 0.0500 0.2750 0.2500 0.1350;
3, 0.1500 0.6300 0.4500 0.0900;
4, 0.7000 0.2350 0.8500 0.1675;
5, 0.6000 0.3800 0.8500 0.2675;
6, 0.3500 0.9714 0.8000 0.7143;
7, 0.7500 0.9574 0.9500 0.8155;
8, 0.1500 0.8363 0.4000 0.9727;
BARRIERS
NUMBER, START_X, START_Y, END_X, END_Y 
1, 0.1500 0.9167 0.4000 0.5000;
2, 0.6500 0.8333 0.8500 0.1667;
\end{Verbatim}

\section{Exact coordinates of fracture networks in Example \ref{ex6_1}} \label{appd:ex6_1}

\begin{Verbatim}
FRACTURES:
NUMBER, START_X, START_Y, END_X, END_Y 
1, 0.2107 0.2846 0.6765 0.0974;
2. 0.1035 0.0993 0.6543 0.2959;
BARRIERS
NUMBER, START_X, START_Y, END_X, END_Y 
1, 0.3031 1.0000 0.7006 0.7996;
2, 0.7006 0.7996 0.7985 0.6011;
3, 0.2902 0.7996 0.5933 0.9457;
4, 0.4529 0.6011 0.7006 0.7996;
5, 0.7006 0.7996 0.8059 0.8652;
6, 0.7468 0.4551 0.8152 0.0449;
\end{Verbatim}

\end{appendices}

\end{document}